\documentclass[11pt,a4paper]{article}

\usepackage{theorem,enumerate}
\usepackage{amsmath,latexsym,amssymb,amsfonts}
\usepackage{eucal}
\usepackage{color}
\usepackage{comment}

\newcounter{Scounter}
\theorembodyfont{\normalfont\slshape}
\newtheorem{thm}{Theorem}[section]

\newtheorem{Thm}{Theorem}

\newtheorem{prop}[thm]{Proposition}
\newtheorem{lem}[thm]{Lemma}

\newtheorem{claim}{Claim}[section]

\numberwithin{equation}{section}

\newcommand{\proof}{\medbreak\noindent\textit{Proof.}\quad}

\newcommand{\qed}{{$\quad\square$\vs{3.6}}}

\newcommand{\vs}[1]{\vspace*{#1 mm}}

\def\GG{{ \mathcal{G}}}
\def\HH{{ \mathcal{H}}}

\bfseries\normalfont

\title{Forbidden subgraphs generating a finite set of graphs with minimum degree three and large girth}

\author{
Yoshimi Egawa$^{1}$ \and\
Michitaka Furuya$^{2}$\footnote{\texttt{e-mail:michitaka.furuya@gmail.com}}\vs{5}\\
$^{1}$\textsl{Department of Applied Mathematics,} \\
\textsl{Tokyo University of Science,}\\
\textsl{1-3 Kagurazaka, Shinjuku-ku, Tokyo 162-8601, Japan }\\
$^{2}$\textsl{College of Liberal Arts and Sciences,}\\
\textsl{Kitasato University,}\\
\textsl{1-15-1 Kitasato, Minami-ku, Sagamihara, Kanagawa 252-0373, Japan}\\
}

\date{}

\begin{document}

\maketitle

\begin{abstract}
For a family $\mathcal{H}$ of graphs, a graph $G$ is said to be {\it $\mathcal{H}$-free} if $G$ contains no member of $\mathcal{H}$ as an induced subgraph.
We let $\tilde{\mathcal{G}}_{3}(\mathcal{H})$ denote the family of connected $\mathcal{H}$-free graphs having minimum degree at least $3$.
In this paper, we characterize the non-caterpillar trees $T$ having diameter at least $7$ such that $\tilde{\mathcal{G}}_{3}(\{C_{3},C_{4},T\})$ is a finite family, where $C_{n}$ is a cycle of order $n$.
\end{abstract}



\section{Introduction}\label{sec1}

All graphs considered in this paper are finite, simple, and undirected.
Let $G$ be a graph.
We let $V(G)$ and $E(G)$ denote the {\it vertex set} and the {\it edge set} of $G$, respectively.
For a vertex $u\in V(G)$, we let $N_{G}(u)$ and $d_{G}(u)$ denote the {\it neighborhood} and the {\it degree} of $u$, respectively; thus $N_{G}(u)=\{v\in V(G): uv\in E(G)\}$ and $d_{G}(u)=|N_{G}(u)|$.
We let $\delta (G)$ and $\Delta (G)$ denote the {\it minimum degree} and the {\it maximum degree} of $G$, respectively.
For a subset $U$ of $V(G)$, we set $N_{G}(U)=\bigcup _{u\in U}N_{G}(u)$, and let $G[U]$ denote the subgraph of $G$ induced by $U$.
For two subsets $U$ and $U'$ of $V(G)$, we let $E_{G}(U,U')$ be the set of edges $uu'$ of $G$ with $u\in U$ and $u'\in U'$.
For $u,v\in V(G)$, let ${\rm dist}_{G}(u,v)$ denote the length of a shortest $u$-$v$ path of $G$.
We let ${\rm diam}(G)$ denote the maximum of ${\rm dist}_{G}(u,v)$ as $u$ and $v$ range over $V(G)$.
We let $P_{n}$ and $C_{n}$ denote the {\it path} and the {\it cycle} of order $n$, respectively.
For two positive integers $s$ and $t$, the {\it Ramsey number} $R(s,t)$ is the minimum positive integer $R$ such that any graph of order at least $R$ contains a clique of cardinality $s$ or an independent set of cardinality $t$.
For terms and symbols not defined in this paper, we refer the reader to \cite{D}.

For two graphs $G$ and $H$, we write $H\prec G$ if $G$ contains an induced copy of $H$, and $G$ is said to be {\it $H$-free} if $H\not\prec G$.
For a family $\HH$ of graphs, a graph $G$ is said to be {\it $\HH$-free} if $G$ is $H$-free for every $H\in \HH$.
In this context, the members of $\HH$ are called {\it forbidden subgraphs}.
For an integer $k\geq 1$ and a family $\HH$ of graphs, let $\tilde{\GG}_{k}(\HH)$ denote the family of connected $\HH$-free graphs with minimum degree at least $k$.
In this paper, we focus on the finiteness of $\tilde{\GG}_{3}(\{C_{3},C_{4},T\})$, where $T$ is a tree having large diameter.

\begin{figure}
\begin{center}
{\unitlength 0.1in%
\begin{picture}(50.5000,8.0500)(0.0000,-9.4000)%
%
\special{sh 1.000}%
\special{ia 200 600 50 50 0.0000000 6.2831853}%
\special{pn 8}%
\special{ar 200 600 50 50 0.0000000 6.2831853}%
%
\special{sh 1.000}%
\special{ia 400 600 50 50 0.0000000 6.2831853}%
\special{pn 8}%
\special{ar 400 600 50 50 0.0000000 6.2831853}%
%
\special{sh 1.000}%
\special{ia 600 600 50 50 0.0000000 6.2831853}%
\special{pn 8}%
\special{ar 600 600 50 50 0.0000000 6.2831853}%
%
\special{sh 1.000}%
\special{ia 800 600 50 50 0.0000000 6.2831853}%
\special{pn 8}%
\special{ar 800 600 50 50 0.0000000 6.2831853}%
%
\special{sh 1.000}%
\special{ia 1400 600 50 50 0.0000000 6.2831853}%
\special{pn 8}%
\special{ar 1400 600 50 50 0.0000000 6.2831853}%
%
\special{sh 1.000}%
\special{ia 1800 600 50 50 0.0000000 6.2831853}%
\special{pn 8}%
\special{ar 1800 600 50 50 0.0000000 6.2831853}%
%
\special{sh 1.000}%
\special{ia 2000 600 50 50 0.0000000 6.2831853}%
\special{pn 8}%
\special{ar 2000 600 50 50 0.0000000 6.2831853}%
%
\special{sh 1.000}%
\special{ia 2200 600 50 50 0.0000000 6.2831853}%
\special{pn 8}%
\special{ar 2200 600 50 50 0.0000000 6.2831853}%
%
\special{sh 1.000}%
\special{ia 2400 600 50 50 0.0000000 6.2831853}%
\special{pn 8}%
\special{ar 2400 600 50 50 0.0000000 6.2831853}%
%
\special{sh 1.000}%
\special{ia 2600 600 50 50 0.0000000 6.2831853}%
\special{pn 8}%
\special{ar 2600 600 50 50 0.0000000 6.2831853}%
%
\special{sh 1.000}%
\special{ia 3000 600 50 50 0.0000000 6.2831853}%
\special{pn 8}%
\special{ar 3000 600 50 50 0.0000000 6.2831853}%
%
\special{sh 1.000}%
\special{ia 3200 600 50 50 0.0000000 6.2831853}%
\special{pn 8}%
\special{ar 3200 600 50 50 0.0000000 6.2831853}%
%
\special{sh 1.000}%
\special{ia 2200 400 50 50 0.0000000 6.2831853}%
\special{pn 8}%
\special{ar 2200 400 50 50 0.0000000 6.2831853}%
%
\special{sh 1.000}%
\special{ia 2200 200 50 50 0.0000000 6.2831853}%
\special{pn 8}%
\special{ar 2200 200 50 50 0.0000000 6.2831853}%
%
\special{sh 1.000}%
\special{ia 2800 600 50 50 0.0000000 6.2831853}%
\special{pn 8}%
\special{ar 2800 600 50 50 0.0000000 6.2831853}%
%
\special{sh 1.000}%
\special{ia 2600 400 50 50 0.0000000 6.2831853}%
\special{pn 8}%
\special{ar 2600 400 50 50 0.0000000 6.2831853}%
%
\special{sh 1.000}%
\special{ia 2800 400 50 50 0.0000000 6.2831853}%
\special{pn 8}%
\special{ar 2800 400 50 50 0.0000000 6.2831853}%
%
\special{pn 8}%
\special{pa 1800 600}%
\special{pa 3200 600}%
\special{fp}%
%
\special{pn 8}%
\special{pa 2200 600}%
\special{pa 2200 200}%
\special{fp}%
%
\special{pn 8}%
\special{pa 2600 400}%
\special{pa 2600 600}%
\special{fp}%
%
\special{pn 8}%
\special{pa 2800 600}%
\special{pa 2800 400}%
\special{fp}%
%
\special{sh 1.000}%
\special{ia 3600 600 50 50 0.0000000 6.2831853}%
\special{pn 8}%
\special{ar 3600 600 50 50 0.0000000 6.2831853}%
%
\special{sh 1.000}%
\special{ia 3800 600 50 50 0.0000000 6.2831853}%
\special{pn 8}%
\special{ar 3800 600 50 50 0.0000000 6.2831853}%
%
\special{sh 1.000}%
\special{ia 4000 600 50 50 0.0000000 6.2831853}%
\special{pn 8}%
\special{ar 4000 600 50 50 0.0000000 6.2831853}%
%
\special{sh 1.000}%
\special{ia 4200 600 50 50 0.0000000 6.2831853}%
\special{pn 8}%
\special{ar 4200 600 50 50 0.0000000 6.2831853}%
%
\special{sh 1.000}%
\special{ia 4400 600 50 50 0.0000000 6.2831853}%
\special{pn 8}%
\special{ar 4400 600 50 50 0.0000000 6.2831853}%
%
\special{sh 1.000}%
\special{ia 4800 600 50 50 0.0000000 6.2831853}%
\special{pn 8}%
\special{ar 4800 600 50 50 0.0000000 6.2831853}%
%
\special{sh 1.000}%
\special{ia 5000 600 50 50 0.0000000 6.2831853}%
\special{pn 8}%
\special{ar 5000 600 50 50 0.0000000 6.2831853}%
%
\special{sh 1.000}%
\special{ia 4000 400 50 50 0.0000000 6.2831853}%
\special{pn 8}%
\special{ar 4000 400 50 50 0.0000000 6.2831853}%
%
\special{sh 1.000}%
\special{ia 4000 200 50 50 0.0000000 6.2831853}%
\special{pn 8}%
\special{ar 4000 200 50 50 0.0000000 6.2831853}%
%
\special{sh 1.000}%
\special{ia 4600 600 50 50 0.0000000 6.2831853}%
\special{pn 8}%
\special{ar 4600 600 50 50 0.0000000 6.2831853}%
%
\special{sh 1.000}%
\special{ia 4200 400 50 50 0.0000000 6.2831853}%
\special{pn 8}%
\special{ar 4200 400 50 50 0.0000000 6.2831853}%
%
\special{pn 8}%
\special{pa 3600 600}%
\special{pa 5000 600}%
\special{fp}%
%
\special{pn 8}%
\special{pa 4000 600}%
\special{pa 4000 200}%
\special{fp}%
%
\special{pn 8}%
\special{pa 4200 600}%
\special{pa 4200 400}%
\special{fp}%
\put(25.0000,-10.0500){\makebox(0,0){$T_{8}^{(1)}$}}%
\put(43.0000,-10.0500){\makebox(0,0){$T_{8}^{(2)}$}}%
%
\special{sh 1.000}%
\special{ia 600 400 50 50 0.0000000 6.2831853}%
\special{pn 8}%
\special{ar 600 400 50 50 0.0000000 6.2831853}%
%
\special{sh 1.000}%
\special{ia 600 200 50 50 0.0000000 6.2831853}%
\special{pn 8}%
\special{ar 600 200 50 50 0.0000000 6.2831853}%
%
\special{pn 8}%
\special{pa 600 600}%
\special{pa 600 200}%
\special{fp}%
\put(8.0000,-10.0500){\makebox(0,0){$T_{n}$}}%
%
\special{pn 8}%
\special{pa 200 600}%
\special{pa 950 600}%
\special{fp}%
%
\special{pn 8}%
\special{pa 1250 600}%
\special{pa 1400 600}%
\special{fp}%
%
\special{pn 4}%
\special{sh 1}%
\special{ar 1000 600 16 16 0 6.2831853}%
\special{sh 1}%
\special{ar 1200 600 16 16 0 6.2831853}%
\special{sh 1}%
\special{ar 1100 600 16 16 0 6.2831853}%
\special{sh 1}%
\special{ar 1100 600 16 16 0 6.2831853}%
\put(2.0000,-7.3000){\makebox(0,0){$u_{1}$}}%
\put(4.0000,-7.3000){\makebox(0,0){$u_{2}$}}%
\put(6.0000,-7.3000){\makebox(0,0){$u_{3}$}}%
\put(8.0000,-7.3000){\makebox(0,0){$u_{4}$}}%
\put(14.0000,-7.3000){\makebox(0,0){$u_{n}$}}%
\put(7.3000,-4.0000){\makebox(0,0){$v_{3}$}}%
\put(7.3000,-2.0000){\makebox(0,0){$v'_{3}$}}%
\end{picture}}%

\caption{Graphs $T_{n}$, $T_{8}^{(1)}$ and $T_{8}^{(2)}$}
\label{f1}
\end{center}
\end{figure}
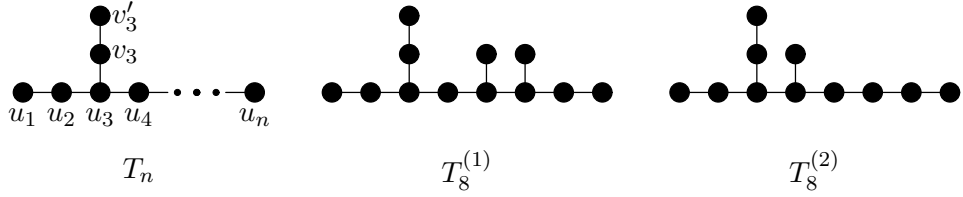

A tree $T$ is a {\it caterpillar} if there exists a path $P$ of $T$ such that $T-V(P)$ has no edges.
For an integer $n\geq 4$, let $T_{n}$ be the graph obtained from a path $u_{1}u_{2}\cdots u_{n}$ of order $n$ by adding two vertices $v_{3},v'_{3}$ and edges $u_{3}v_{3},v_{3}v'_{3}$ (see the left graph in Figure~\ref{f1}).
Let $T_{8}^{(1)}$ and $T_{8}^{(2)}$ be the graphs depicted in the center and the right in Figure~\ref{f1}, respectively.
The following is our main theorem.

\begin{thm}
\label{thm-diam7-noncat}
Let $T$ be a tree with ${\rm diam}(T)\geq 7$, and suppose that $T$ is not a caterpillar.
Then $\tilde{\GG}_{3}(\{C_{3},C_{4},T\})$ is a finite family if and only if $T$ is a subgraph of one of $T_{8}^{(1)}$, $T_{8}^{(2)}$ and $T_{9}$.
\end{thm}

As for the ``only if'' part of Theorem~\ref{thm-diam7-noncat}, we prove the following theorem in Section~\ref{sec-nec}.

\begin{thm}
\label{thm-nec-1}
Let $T$ be a tree with ${\rm diam}(T)\geq 7$, and suppose that $T$ is not a caterpillar.
If $\tilde{\GG}_{3}(\{C_{3},C_{4},T\})$ is a finite family, then $T$ is a subgraph of one of $T_{8}^{(1)}$, $T_{8}^{(2)}$ and $T_{9}$.
\end{thm}

Note that a connected graph $G$ with ${\rm diam}(G)\geq 3$ satisfies $|V(G)|\leq \Delta (G)^{{\rm diam}(G)}$ (see, for example, \cite{EFFPS}).
Thus we divide the ``if'' part of Theorem~\ref{thm-diam7-noncat} into two assertions as follows.

\begin{thm}
\label{prop-diambound}
Let $G$ be a connected $\{C_{3},C_{4}\}$-free graph with $\delta (G)\geq 3$.
Then the following hold.
\begin{enumerate}
\item[{\upshape(i)}]
If ${\rm diam}(G)\geq 20$, then $T_{8}^{(1)}\prec G$.
\item[{\upshape(ii)}]
If ${\rm diam}(G)\geq 16$, then $T_{8}^{(2)}\prec G$.
\item[{\upshape(iii)}]
If ${\rm diam}(G)\geq 12$, then $T_{9}\prec G$.
\end{enumerate}
\end{thm}

\begin{thm}
\label{prop-maxdegbound}
Let $G$ be a connected $\{C_{3},C_{4}\}$-free graph with $\delta (G)\geq 3$.
Then the following hold.
\begin{enumerate}
\item[{\upshape(i)}]
If $\Delta (G)\geq 943218$, then $T_{8}^{(1)}\prec G$.
\item[{\upshape(ii)}]
If $\Delta (G)\geq 190375$, then $T_{8}^{(2)}\prec G$.
\item[{\upshape(iii)}]
If $\Delta (G)\geq 197433$, then $T_{9}\prec G$.
\end{enumerate}
\end{thm}

Theorem~\ref{prop-diambound} is proved in Section~\ref{sec-diam}.
In Section~\ref{sec-M-cond}, we prepare some lemmas for the proof of Theorem~\ref{prop-maxdegbound}, and using them, we prove Theorem~\ref{prop-maxdegbound} in Section~\ref{sec-maxdeg}.

\subsection{Motivations and applications of Theorem~\ref{thm-diam7-noncat}}\label{sec1-1}

Our motivation of Theorem~\ref{thm-diam7-noncat} derives from several different lines of research.
The first one is characterization of families $\HH$ of connected graphs satisfying the condition that
\begin{align}
\mbox{the number of $k$-connected $\HH$-free graphs is finite.}\label{cond-intro-1}
\end{align}
If a family $\HH$ satisfies (\ref{cond-intro-1}), then for any property $P$ on graphs, although the proposition that
\begin{align}
\mbox{all $k$-connected $\HH$-free graphs satisfy $P$ with finite exceptions}\label{cond-intro-2}
\end{align}
holds, the proposition gives no information about $P$.
Thus it is important to identify families $\HH$ satisfying (\ref{cond-intro-1}) in advance.
Having such a motivation, Fujisawa, Plummer and Saito~\cite{FPS} started a study of $\HH$ satisfying (\ref{cond-intro-1}), and they determined the families $\HH$ satisfying (\ref{cond-intro-1}) for the case where $1\leq k\leq 6$ and $|\HH|\leq 2$.
In \cite{EFFPS,EF,EZ}, the research was continued by analyzing families $\HH$ satisfying (\ref{cond-intro-1}) for the case where $k=3$ and $|\HH|=3$.
However, at present, the state of research is far from complete characterization of $\HH$ for this case.
Theorem~\ref{thm-diam7-noncat} concerned with the corresponding problem for graphs with minimum degree at least $3$.
In view of the fact that connectivity conditions can often be replaced by minimum degree conditions in propositions like (\ref{cond-intro-2}), it is natural to consider connected graphs with minimum degree at least $3$ in place of $3$-connected graphs.
In fact, the authors~\cite{EF2} recently characterized the families $\HH$ with $|\HH|=3$ and $\{C_{3},C_{4}\}\not\subseteq \HH$ such that $\tilde{\GG}_{3}(\HH)$ is finite.
Thus for the finiteness problem on $\tilde{\GG}_{3}(\HH)$ with $|\HH|=3$, it remains to consider the case where $\{C_{3},C_{4}\}\subseteq \HH $.
Moreover, research done in \cite{EFFPS,EF,EZ} has revealed that caterpillars and trees with small diameter are easy to handle in the study of forbidden subgraphs.
Thus Theorem~\ref{thm-diam7-noncat} solves the most difficult case of the finiteness problem on $\tilde{\GG}_{3}(\HH)$ with $|\HH|=3$ (we plan to discuss other cases in a subsequent paper).

Our second motivation is a coloring problem.
For a graph $G$, let $\chi (G)$ denote the {\it chromatic number} of $G$, that is, the minimum number of $k$ such that there exists a function $c:V(G)\rightarrow \{1,\ldots ,k\}$ with $c(u)\neq c(v)$ for any adjacent vertices $u$ and $v$ of $G$.
In recent years, there have been remarkable progress in the study of the relationship between the chromatic number and forbidden subgraphs (see a survey \cite{SS}).
Most notably, Chudnovsky and Stacho~\cite{CS} focused on $3$-colorability of $P_{8}$-free graphs, and showed that every $\{C_{3},C_{4},P_{8}\}$-free graph is $3$-colorable.
They indeed proved the following stronger theorem which immediately assures us the $3$-colorability of $\{C_{3},C_{4},P_{8}\}$-free graphs (with induction argument).

\begin{Thm}[Chudnovsky and Stacho~\cite{CS}]
\label{Thm-P8-free}
Let $G$ be a connected $\{C_{3},C_{4},P_{8}\}$-free graph with $\delta (G)\geq 3$.
Then $G$ is the Petersen graph, the Heawood graph, or the graph obtained from the Heawood graph by contracting one edge.
In particular, $\chi (G)\leq 3$.
\end{Thm}

Since there is a polynomial-time algorithm (with respect to $|V(G)|$) to judge whether a given graph $G$ is bipartite or not, we can calculate the chromatic number of $\{C_{3},C_{4},P_{8}\}$-free graphs in polynomial time based on Theorem~\ref{Thm-P8-free}.
Here a natural question arises:
What forbidden subgraph conditions (similar to $\{C_{3},C_{4},P_{8}\}$-freeness) can yield a polynomial-time algorithm to determine the chromatic number?
In view of Theorem~\ref{Thm-P8-free}, $\{C_{3},C_{4},T\}$-freeness for a tree $T$ with large diameter appears to be a condition worth considering.
Our main result provides us with examples of such conditions which give a desired algorithm as follows.

We let $T$ be an induced subgraph of one of $T_{8}^{(1)}$, $T_{8}^{(2)}$ and $T_{9}$, and determine the chromatic number of a $\{C_{3},C_{4},T\}$-free graph $G$.
By Theorem~\ref{thm-diam7-noncat}, there exists a constant number $s=s(T)$ such that every connected $\{C_{3},C_{4},T\}$-free graph with minimum degree at least $3$ has at most $s$ vertices.
We first judge whether $G$ is bipartite or not in polynomial time.
We may assume that $G$ is judged to be a non-bipartite graph, i.e., $\chi (G)\geq 3$.
Next we find a maximal sequence $u_{1},u_{2},\ldots ,u_{m}$ of vertices of $G$ such that the degree of $u_{i}$ in $G-\{u_{1},\ldots ,u_{i-1}\}$ is at most $2$ (in polynomial time).
Since $\chi (G)\geq 3$, we have $\chi (G)=\max\{\chi (G-\{u_{1},\ldots ,u_{i}\}),3\}$ for all $i~(1\leq i\leq m)$ by induction on $i$.
In particular, $\chi (G)=\max\{\chi (G-\{u_{1},\ldots ,u_{m}\}),3\}$.
If $G-\{u_{1},\ldots ,u_{m}\}$ is an empty graph, then $\chi (G)=3$.
Thus we may assume that $G-\{u_{1},\ldots ,u_{m}\}$ is non-empty.
Let $G_{1},\ldots ,G_{p}$ be the components of $G-\{u_{1},\ldots ,u_{m}\}$.
Then $G_{i}$ is a connected $\{C_{3},C_{4},T\}$-free graph with $\delta (G_{i})\geq 3$, and hence $|V(G_{i})|\leq s$.
Since $s$ does not depend on the order of $G$, we can determine the chromatic number of $G_{i}$ in constant time (here the constant is $O(2^{s}s)$~\cite{BHK}).
Since $p\leq |V(G)|$, we calculate the value $\max\{\chi (G_{i}):1\leq i\leq p\}~(=\chi (G-\{u_{1},\ldots ,u_{m}\}))$ in linear time.
Therefore we can determine the value $\chi (G)$ in a polynomial time with respect to $|V(G)|$.
Note that the finiteness of $\tilde{\GG}_{3}(\{C_{3},C_{4},T\})$ is effectively used in the above strategy.

\section{Proof of Theorem~\ref{thm-nec-1}}\label{sec-nec}

In this section, we prove Theorem~\ref{thm-nec-1}.

\begin{figure}
\begin{center}
{\unitlength 0.1in%
\begin{picture}(55.0000,8.1000)(1.5000,-9.4000)%
%
\special{sh 1.000}%
\special{ia 200 600 50 50 0.0000000 6.2831853}%
\special{pn 8}%
\special{ar 200 600 50 50 0.0000000 6.2831853}%
%
\special{sh 1.000}%
\special{ia 400 600 50 50 0.0000000 6.2831853}%
\special{pn 8}%
\special{ar 400 600 50 50 0.0000000 6.2831853}%
%
\special{sh 1.000}%
\special{ia 600 600 50 50 0.0000000 6.2831853}%
\special{pn 8}%
\special{ar 600 600 50 50 0.0000000 6.2831853}%
%
\special{sh 1.000}%
\special{ia 800 600 50 50 0.0000000 6.2831853}%
\special{pn 8}%
\special{ar 800 600 50 50 0.0000000 6.2831853}%
%
\special{sh 1.000}%
\special{ia 1800 600 50 50 0.0000000 6.2831853}%
\special{pn 8}%
\special{ar 1800 600 50 50 0.0000000 6.2831853}%
%
\special{sh 1.000}%
\special{ia 2000 600 50 50 0.0000000 6.2831853}%
\special{pn 8}%
\special{ar 2000 600 50 50 0.0000000 6.2831853}%
%
\special{sh 1.000}%
\special{ia 600 400 50 50 0.0000000 6.2831853}%
\special{pn 8}%
\special{ar 600 400 50 50 0.0000000 6.2831853}%
%
\special{sh 1.000}%
\special{ia 600 200 50 50 0.0000000 6.2831853}%
\special{pn 8}%
\special{ar 600 200 50 50 0.0000000 6.2831853}%
%
\special{pn 8}%
\special{pa 600 600}%
\special{pa 600 200}%
\special{fp}%
\put(11.0000,-10.0500){\makebox(0,0){$T^{*}_{n}$}}%
%
\special{pn 8}%
\special{pa 200 600}%
\special{pa 950 600}%
\special{fp}%
%
\special{pn 8}%
\special{pa 1250 600}%
\special{pa 2000 600}%
\special{fp}%
%
\special{pn 4}%
\special{sh 1}%
\special{ar 1000 600 16 16 0 6.2831853}%
\special{sh 1}%
\special{ar 1200 600 16 16 0 6.2831853}%
\special{sh 1}%
\special{ar 1100 600 16 16 0 6.2831853}%
\special{sh 1}%
\special{ar 1100 600 16 16 0 6.2831853}%
%
\special{sh 1.000}%
\special{ia 1400 600 50 50 0.0000000 6.2831853}%
\special{pn 8}%
\special{ar 1400 600 50 50 0.0000000 6.2831853}%
%
\special{sh 1.000}%
\special{ia 1600 400 50 50 0.0000000 6.2831853}%
\special{pn 8}%
\special{ar 1600 400 50 50 0.0000000 6.2831853}%
%
\special{sh 1.000}%
\special{ia 1600 200 50 50 0.0000000 6.2831853}%
\special{pn 8}%
\special{ar 1600 200 50 50 0.0000000 6.2831853}%
%
\special{pn 8}%
\special{pa 1600 600}%
\special{pa 1600 200}%
\special{fp}%
%
\special{sh 1.000}%
\special{ia 1400 600 50 50 0.0000000 6.2831853}%
\special{pn 8}%
\special{ar 1400 600 50 50 0.0000000 6.2831853}%
%
\special{sh 1.000}%
\special{ia 1600 600 50 50 0.0000000 6.2831853}%
\special{pn 8}%
\special{ar 1600 600 50 50 0.0000000 6.2831853}%
%
\special{sh 1.000}%
\special{ia 1600 600 50 50 0.0000000 6.2831853}%
\special{pn 8}%
\special{ar 1600 600 50 50 0.0000000 6.2831853}%
%
\special{sh 1.000}%
\special{ia 2600 600 50 50 0.0000000 6.2831853}%
\special{pn 8}%
\special{ar 2600 600 50 50 0.0000000 6.2831853}%
%
\special{sh 1.000}%
\special{ia 2800 600 50 50 0.0000000 6.2831853}%
\special{pn 8}%
\special{ar 2800 600 50 50 0.0000000 6.2831853}%
%
\special{sh 1.000}%
\special{ia 3000 600 50 50 0.0000000 6.2831853}%
\special{pn 8}%
\special{ar 3000 600 50 50 0.0000000 6.2831853}%
%
\special{sh 1.000}%
\special{ia 3200 600 50 50 0.0000000 6.2831853}%
\special{pn 8}%
\special{ar 3200 600 50 50 0.0000000 6.2831853}%
%
\special{sh 1.000}%
\special{ia 3000 400 50 50 0.0000000 6.2831853}%
\special{pn 8}%
\special{ar 3000 400 50 50 0.0000000 6.2831853}%
%
\special{sh 1.000}%
\special{ia 3000 200 50 50 0.0000000 6.2831853}%
\special{pn 8}%
\special{ar 3000 200 50 50 0.0000000 6.2831853}%
%
\special{pn 8}%
\special{pa 3000 600}%
\special{pa 3000 200}%
\special{fp}%
%
\special{pn 8}%
\special{pa 2600 600}%
\special{pa 3800 600}%
\special{fp}%
%
\special{sh 1.000}%
\special{ia 3400 600 50 50 0.0000000 6.2831853}%
\special{pn 8}%
\special{ar 3400 600 50 50 0.0000000 6.2831853}%
%
\special{sh 1.000}%
\special{ia 3600 600 50 50 0.0000000 6.2831853}%
\special{pn 8}%
\special{ar 3600 600 50 50 0.0000000 6.2831853}%
%
\special{sh 1.000}%
\special{ia 3800 600 50 50 0.0000000 6.2831853}%
\special{pn 8}%
\special{ar 3800 600 50 50 0.0000000 6.2831853}%
%
\special{sh 1.000}%
\special{ia 3600 400 50 50 0.0000000 6.2831853}%
\special{pn 8}%
\special{ar 3600 400 50 50 0.0000000 6.2831853}%
%
\special{pn 8}%
\special{pa 3600 600}%
\special{pa 3600 400}%
\special{fp}%
\put(33.0000,-10.0500){\makebox(0,0){$S_{7}(1,0,1)$}}%
%
\special{sh 1.000}%
\special{ia 4200 600 50 50 0.0000000 6.2831853}%
\special{pn 8}%
\special{ar 4200 600 50 50 0.0000000 6.2831853}%
%
\special{sh 1.000}%
\special{ia 4400 600 50 50 0.0000000 6.2831853}%
\special{pn 8}%
\special{ar 4400 600 50 50 0.0000000 6.2831853}%
%
\special{sh 1.000}%
\special{ia 4600 600 50 50 0.0000000 6.2831853}%
\special{pn 8}%
\special{ar 4600 600 50 50 0.0000000 6.2831853}%
%
\special{sh 1.000}%
\special{ia 4800 600 50 50 0.0000000 6.2831853}%
\special{pn 8}%
\special{ar 4800 600 50 50 0.0000000 6.2831853}%
%
\special{sh 1.000}%
\special{ia 4600 400 50 50 0.0000000 6.2831853}%
\special{pn 8}%
\special{ar 4600 400 50 50 0.0000000 6.2831853}%
%
\special{sh 1.000}%
\special{ia 4600 200 50 50 0.0000000 6.2831853}%
\special{pn 8}%
\special{ar 4600 200 50 50 0.0000000 6.2831853}%
%
\special{pn 8}%
\special{pa 4600 600}%
\special{pa 4600 200}%
\special{fp}%
%
\special{pn 8}%
\special{pa 4200 600}%
\special{pa 5600 600}%
\special{fp}%
%
\special{sh 1.000}%
\special{ia 5000 600 50 50 0.0000000 6.2831853}%
\special{pn 8}%
\special{ar 5000 600 50 50 0.0000000 6.2831853}%
%
\special{sh 1.000}%
\special{ia 5200 600 50 50 0.0000000 6.2831853}%
\special{pn 8}%
\special{ar 5200 600 50 50 0.0000000 6.2831853}%
%
\special{sh 1.000}%
\special{ia 5400 600 50 50 0.0000000 6.2831853}%
\special{pn 8}%
\special{ar 5400 600 50 50 0.0000000 6.2831853}%
%
\special{sh 1.000}%
\special{ia 5600 600 50 50 0.0000000 6.2831853}%
\special{pn 8}%
\special{ar 5600 600 50 50 0.0000000 6.2831853}%
%
\special{sh 1.000}%
\special{ia 5000 400 50 50 0.0000000 6.2831853}%
\special{pn 8}%
\special{ar 5000 400 50 50 0.0000000 6.2831853}%
%
\special{pn 8}%
\special{pa 5000 600}%
\special{pa 5000 400}%
\special{fp}%
%
\special{sh 1.000}%
\special{ia 5200 400 50 50 0.0000000 6.2831853}%
\special{pn 8}%
\special{ar 5200 400 50 50 0.0000000 6.2831853}%
%
\special{pn 8}%
\special{pa 5200 600}%
\special{pa 5200 400}%
\special{fp}%
\put(49.0000,-10.0500){\makebox(0,0){$S_{8}(1,1,1,0)$}}%
%
\special{sh 1.000}%
\special{ia 3200 400 50 50 0.0000000 6.2831853}%
\special{pn 8}%
\special{ar 3200 400 50 50 0.0000000 6.2831853}%
%
\special{pn 8}%
\special{pa 3200 600}%
\special{pa 3200 400}%
\special{fp}%
%
\special{sh 1.000}%
\special{ia 4800 400 50 50 0.0000000 6.2831853}%
\special{pn 8}%
\special{ar 4800 400 50 50 0.0000000 6.2831853}%
%
\special{pn 8}%
\special{pa 4800 600}%
\special{pa 4800 400}%
\special{fp}%
\put(18.3000,-3.9500){\makebox(0,0){$w_{n-2}$}}%
\put(18.3000,-1.9500){\makebox(0,0){$w'_{n-2}$}}%
\put(32.0000,-7.2500){\makebox(0,0){$u_{4}$}}%
\put(34.0000,-7.2500){\makebox(0,0){$u_{5}$}}%
\put(36.0000,-7.2000){\makebox(0,0){$u_{6}$}}%
\put(32.0000,-2.9000){\makebox(0,0){$w_{4}$}}%
\put(36.0000,-2.9000){\makebox(0,0){$w_{6}$}}%
\put(48.0000,-2.9000){\makebox(0,0){$w_{4}$}}%
\put(50.0000,-2.9000){\makebox(0,0){$w_{5}$}}%
\put(52.0000,-2.9000){\makebox(0,0){$w_{6}$}}%
\put(50.0000,-7.2500){\makebox(0,0){$u_{5}$}}%
\put(52.0000,-7.2000){\makebox(0,0){$u_{6}$}}%
\put(54.0000,-7.2000){\makebox(0,0){$u_{7}$}}%
\put(48.0000,-7.2000){\makebox(0,0){$u_{4}$}}%
\put(16.0000,-7.3000){\makebox(0,0){$u_{n-2}$}}%
\end{picture}}%

\caption{Graphs $T^{*}_{n}$, $S_{7}(1,0,1)$ and $S_{8}(1,1,1,0)$}
\label{f-gen0}
\end{center}
\end{figure}
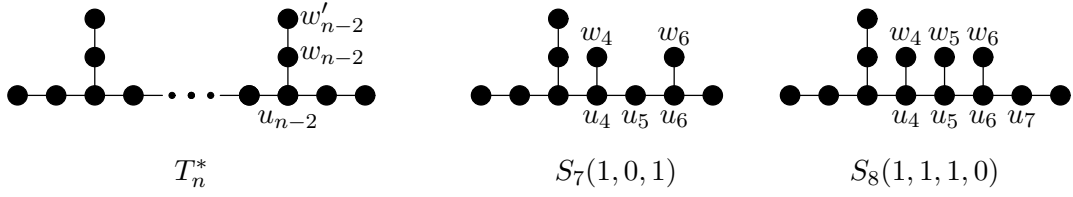

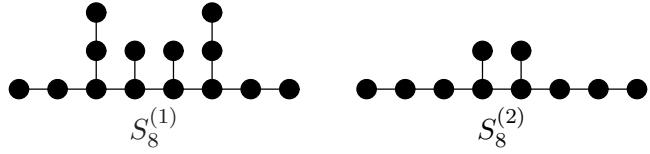
\begin{figure}
\begin{center}
{\unitlength 0.1in%
\begin{picture}(33.0000,5.8500)(1.5000,-7.3500)%
\put(9.0000,-8.0000){\makebox(0,0){$S_{8}^{(1)}$}}%
%
\special{sh 1.000}%
\special{ia 200 600 50 50 0.0000000 6.2831853}%
\special{pn 8}%
\special{ar 200 600 50 50 0.0000000 6.2831853}%
%
\special{sh 1.000}%
\special{ia 400 600 50 50 0.0000000 6.2831853}%
\special{pn 8}%
\special{ar 400 600 50 50 0.0000000 6.2831853}%
%
\special{sh 1.000}%
\special{ia 600 600 50 50 0.0000000 6.2831853}%
\special{pn 8}%
\special{ar 600 600 50 50 0.0000000 6.2831853}%
%
\special{sh 1.000}%
\special{ia 800 600 50 50 0.0000000 6.2831853}%
\special{pn 8}%
\special{ar 800 600 50 50 0.0000000 6.2831853}%
%
\special{sh 1.000}%
\special{ia 1000 600 50 50 0.0000000 6.2831853}%
\special{pn 8}%
\special{ar 1000 600 50 50 0.0000000 6.2831853}%
%
\special{sh 1.000}%
\special{ia 1200 600 50 50 0.0000000 6.2831853}%
\special{pn 8}%
\special{ar 1200 600 50 50 0.0000000 6.2831853}%
%
\special{sh 1.000}%
\special{ia 1400 600 50 50 0.0000000 6.2831853}%
\special{pn 8}%
\special{ar 1400 600 50 50 0.0000000 6.2831853}%
%
\special{sh 1.000}%
\special{ia 1600 600 50 50 0.0000000 6.2831853}%
\special{pn 8}%
\special{ar 1600 600 50 50 0.0000000 6.2831853}%
%
\special{pn 8}%
\special{pa 200 600}%
\special{pa 1600 600}%
\special{fp}%
%
\special{sh 1.000}%
\special{ia 2000 600 50 50 0.0000000 6.2831853}%
\special{pn 8}%
\special{ar 2000 600 50 50 0.0000000 6.2831853}%
%
\special{sh 1.000}%
\special{ia 2200 600 50 50 0.0000000 6.2831853}%
\special{pn 8}%
\special{ar 2200 600 50 50 0.0000000 6.2831853}%
%
\special{sh 1.000}%
\special{ia 2400 600 50 50 0.0000000 6.2831853}%
\special{pn 8}%
\special{ar 2400 600 50 50 0.0000000 6.2831853}%
%
\special{sh 1.000}%
\special{ia 2600 600 50 50 0.0000000 6.2831853}%
\special{pn 8}%
\special{ar 2600 600 50 50 0.0000000 6.2831853}%
%
\special{sh 1.000}%
\special{ia 2800 600 50 50 0.0000000 6.2831853}%
\special{pn 8}%
\special{ar 2800 600 50 50 0.0000000 6.2831853}%
%
\special{sh 1.000}%
\special{ia 3200 600 50 50 0.0000000 6.2831853}%
\special{pn 8}%
\special{ar 3200 600 50 50 0.0000000 6.2831853}%
%
\special{sh 1.000}%
\special{ia 3400 600 50 50 0.0000000 6.2831853}%
\special{pn 8}%
\special{ar 3400 600 50 50 0.0000000 6.2831853}%
%
\special{sh 1.000}%
\special{ia 3000 600 50 50 0.0000000 6.2831853}%
\special{pn 8}%
\special{ar 3000 600 50 50 0.0000000 6.2831853}%
%
\special{sh 1.000}%
\special{ia 2600 400 50 50 0.0000000 6.2831853}%
\special{pn 8}%
\special{ar 2600 400 50 50 0.0000000 6.2831853}%
%
\special{sh 1.000}%
\special{ia 2800 400 50 50 0.0000000 6.2831853}%
\special{pn 8}%
\special{ar 2800 400 50 50 0.0000000 6.2831853}%
%
\special{pn 8}%
\special{pa 2000 600}%
\special{pa 3400 600}%
\special{fp}%
%
\special{pn 8}%
\special{pa 2600 400}%
\special{pa 2600 600}%
\special{fp}%
%
\special{pn 8}%
\special{pa 2800 600}%
\special{pa 2800 400}%
\special{fp}%
\put(27.0000,-8.0000){\makebox(0,0){$S_{8}^{(2)}$}}%
%
\special{sh 1.000}%
\special{ia 600 400 50 50 0.0000000 6.2831853}%
\special{pn 8}%
\special{ar 600 400 50 50 0.0000000 6.2831853}%
%
\special{sh 1.000}%
\special{ia 600 200 50 50 0.0000000 6.2831853}%
\special{pn 8}%
\special{ar 600 200 50 50 0.0000000 6.2831853}%
%
\special{pn 8}%
\special{pa 600 600}%
\special{pa 600 200}%
\special{fp}%
%
\special{sh 1.000}%
\special{ia 1200 400 50 50 0.0000000 6.2831853}%
\special{pn 8}%
\special{ar 1200 400 50 50 0.0000000 6.2831853}%
%
\special{sh 1.000}%
\special{ia 1200 200 50 50 0.0000000 6.2831853}%
\special{pn 8}%
\special{ar 1200 200 50 50 0.0000000 6.2831853}%
%
\special{pn 8}%
\special{pa 1200 600}%
\special{pa 1200 200}%
\special{fp}%
%
\special{sh 1.000}%
\special{ia 800 400 50 50 0.0000000 6.2831853}%
\special{pn 8}%
\special{ar 800 400 50 50 0.0000000 6.2831853}%
%
\special{pn 8}%
\special{pa 800 400}%
\special{pa 800 600}%
\special{fp}%
%
\special{sh 1.000}%
\special{ia 1000 400 50 50 0.0000000 6.2831853}%
\special{pn 8}%
\special{ar 1000 400 50 50 0.0000000 6.2831853}%
%
\special{pn 8}%
\special{pa 1000 400}%
\special{pa 1000 600}%
\special{fp}%
\end{picture}}%

\caption{Graphs $S_{8}^{(1)}$ and $S_{8}^{(2)}$}
\label{f-gen1}
\end{center}
\end{figure}

For an integer $n\geq 6$, let $T^{*}_{n}$ be the graph obtained from $T_{n}$ by adding two vertices $w_{n-2},w'_{n-2}$ and edges $u_{n-2}w_{n-2},w_{n-2}w'_{n-2}$ (see the left graph in Figure~\ref{f-gen0}).
For integers $n\geq 5$ and $p_{4},\ldots ,p_{n-1}\in \{0,1\}$, let $S_{n}(p_{4},\ldots ,p_{n-1})$ be the graph obtained from $T_{n}$ by adding $\sum _{4\leq i\leq n-1}p_{i}$ vertices $w_{i}$ and edges $u_{i}w_{i}$ with $4\leq i\leq n-1$ and $p_{i}=1$ (see the center graph and the right graph in Figure~\ref{f-gen0}).
Note that $T_{n}(0,\ldots ,0)=T_{n}$, $S_{8}(0,1,1,0)\simeq T_{8}^{(1)}$ and $S_{8}(1,0,0,0)\simeq T_{8}^{(2)}$.
Let $S_{8}^{(1)}$ and $S_{8}^{(2)}$ be the graphs depicted in Figure~\ref{f-gen1}.

We start with several lemmas.
The following lemma was proved in \cite{K}.

\begin{lem}[Kochol~\cite{K}]
\label{lem-nec-1}
For every integer $g\geq 3$, there exist infinitely many $3$-connected $3$-regular graphs with girth at least $g$.
\end{lem}

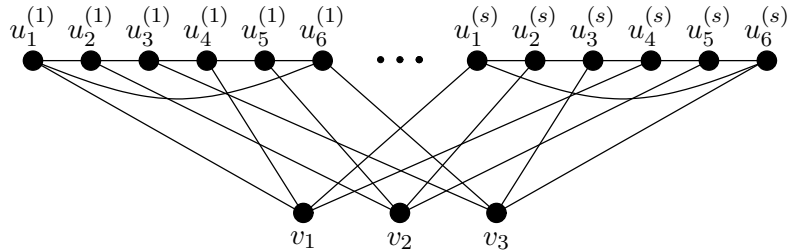
\begin{figure}
\begin{center}
{\unitlength 0.1in%
\begin{picture}(42.1000,11.2000)(2.4000,-12.8000)%
%
\special{sh 1.000}%
\special{ia 600 405 50 50 0.0000000 6.2831853}%
\special{pn 8}%
\special{ar 600 405 50 50 0.0000000 6.2831853}%
%
\special{sh 1.000}%
\special{ia 900 405 50 50 0.0000000 6.2831853}%
\special{pn 8}%
\special{ar 900 405 50 50 0.0000000 6.2831853}%
%
\special{sh 1.000}%
\special{ia 1500 405 50 50 0.0000000 6.2831853}%
\special{pn 8}%
\special{ar 1500 405 50 50 0.0000000 6.2831853}%
%
\special{sh 1.000}%
\special{ia 1200 405 50 50 0.0000000 6.2831853}%
\special{pn 8}%
\special{ar 1200 405 50 50 0.0000000 6.2831853}%
%
\special{sh 1.000}%
\special{ia 2100 405 50 50 0.0000000 6.2831853}%
\special{pn 8}%
\special{ar 2100 405 50 50 0.0000000 6.2831853}%
%
\special{sh 1.000}%
\special{ia 1800 405 50 50 0.0000000 6.2831853}%
\special{pn 8}%
\special{ar 1800 405 50 50 0.0000000 6.2831853}%
\put(6.0000,-2.2500){\makebox(0,0){$u^{(1)}_{1}$}}%
\put(12.0000,-2.2500){\makebox(0,0){$u^{(1)}_{3}$}}%
\put(18.0000,-2.2500){\makebox(0,0){$u^{(1)}_{5}$}}%
\put(21.0000,-2.2500){\makebox(0,0){$u^{(1)}_{6}$}}%
\put(15.0000,-2.2500){\makebox(0,0){$u^{(1)}_{4}$}}%
\put(9.0000,-2.2500){\makebox(0,0){$u^{(1)}_{2}$}}%
%
\special{pn 8}%
\special{pa 600 405}%
\special{pa 631 417}%
\special{pa 662 430}%
\special{pa 786 478}%
\special{pa 816 489}%
\special{pa 909 522}%
\special{pa 940 532}%
\special{pa 1033 559}%
\special{pa 1064 567}%
\special{pa 1126 581}%
\special{pa 1157 587}%
\special{pa 1187 592}%
\special{pa 1249 600}%
\special{pa 1311 604}%
\special{pa 1342 605}%
\special{pa 1373 605}%
\special{pa 1435 601}%
\special{pa 1466 598}%
\special{pa 1528 590}%
\special{pa 1559 584}%
\special{pa 1589 578}%
\special{pa 1620 571}%
\special{pa 1682 555}%
\special{pa 1744 537}%
\special{pa 1806 517}%
\special{pa 1899 484}%
\special{pa 1930 472}%
\special{pa 1960 460}%
\special{pa 2053 424}%
\special{pa 2084 411}%
\special{pa 2100 405}%
\special{fp}%
%
\special{pn 8}%
\special{pa 2100 405}%
\special{pa 600 405}%
\special{fp}%
%
\special{sh 1.000}%
\special{ia 2900 405 50 50 0.0000000 6.2831853}%
\special{pn 8}%
\special{ar 2900 405 50 50 0.0000000 6.2831853}%
%
\special{sh 1.000}%
\special{ia 3200 405 50 50 0.0000000 6.2831853}%
\special{pn 8}%
\special{ar 3200 405 50 50 0.0000000 6.2831853}%
%
\special{sh 1.000}%
\special{ia 3800 405 50 50 0.0000000 6.2831853}%
\special{pn 8}%
\special{ar 3800 405 50 50 0.0000000 6.2831853}%
%
\special{sh 1.000}%
\special{ia 3500 405 50 50 0.0000000 6.2831853}%
\special{pn 8}%
\special{ar 3500 405 50 50 0.0000000 6.2831853}%
%
\special{sh 1.000}%
\special{ia 4400 405 50 50 0.0000000 6.2831853}%
\special{pn 8}%
\special{ar 4400 405 50 50 0.0000000 6.2831853}%
%
\special{sh 1.000}%
\special{ia 4100 405 50 50 0.0000000 6.2831853}%
\special{pn 8}%
\special{ar 4100 405 50 50 0.0000000 6.2831853}%
\put(29.0000,-2.2500){\makebox(0,0){$u^{(s)}_{1}$}}%
\put(35.0000,-2.2500){\makebox(0,0){$u^{(s)}_{3}$}}%
\put(41.0000,-2.2500){\makebox(0,0){$u^{(s)}_{5}$}}%
\put(44.0000,-2.2500){\makebox(0,0){$u^{(s)}_{6}$}}%
\put(38.0000,-2.2500){\makebox(0,0){$u^{(s)}_{4}$}}%
\put(32.0000,-2.2500){\makebox(0,0){$u^{(s)}_{2}$}}%
%
\special{pn 8}%
\special{pa 2900 405}%
\special{pa 2931 417}%
\special{pa 2962 430}%
\special{pa 3086 478}%
\special{pa 3116 489}%
\special{pa 3209 522}%
\special{pa 3240 532}%
\special{pa 3333 559}%
\special{pa 3364 567}%
\special{pa 3426 581}%
\special{pa 3457 587}%
\special{pa 3487 592}%
\special{pa 3549 600}%
\special{pa 3611 604}%
\special{pa 3642 605}%
\special{pa 3673 605}%
\special{pa 3735 601}%
\special{pa 3766 598}%
\special{pa 3828 590}%
\special{pa 3859 584}%
\special{pa 3889 578}%
\special{pa 3920 571}%
\special{pa 3982 555}%
\special{pa 4044 537}%
\special{pa 4106 517}%
\special{pa 4199 484}%
\special{pa 4230 472}%
\special{pa 4260 460}%
\special{pa 4353 424}%
\special{pa 4384 411}%
\special{pa 4400 405}%
\special{fp}%
%
\special{pn 8}%
\special{pa 4400 405}%
\special{pa 2900 405}%
\special{fp}%
%
\special{pn 4}%
\special{sh 1}%
\special{ar 2500 400 16 16 0 6.2831853}%
\special{sh 1}%
\special{ar 2400 400 16 16 0 6.2831853}%
\special{sh 1}%
\special{ar 2600 400 16 16 0 6.2831853}%
\special{sh 1}%
\special{ar 2600 400 16 16 0 6.2831853}%
%
\special{sh 1.000}%
\special{ia 2500 1200 50 50 0.0000000 6.2831853}%
\special{pn 8}%
\special{ar 2500 1200 50 50 0.0000000 6.2831853}%
%
\special{sh 1.000}%
\special{ia 3000 1200 50 50 0.0000000 6.2831853}%
\special{pn 8}%
\special{ar 3000 1200 50 50 0.0000000 6.2831853}%
\put(25.0000,-13.4500){\makebox(0,0){$v_{2}$}}%
\put(30.0000,-13.4500){\makebox(0,0){$v_{3}$}}%
%
\special{sh 1.000}%
\special{ia 2000 1200 50 50 0.0000000 6.2831853}%
\special{pn 8}%
\special{ar 2000 1200 50 50 0.0000000 6.2831853}%
\put(20.0000,-13.4000){\makebox(0,0){$v_{1}$}}%
%
\special{pn 8}%
\special{pa 2000 1205}%
\special{pa 600 405}%
\special{fp}%
\special{pa 1500 405}%
\special{pa 2000 1205}%
\special{fp}%
\special{pa 2000 1205}%
\special{pa 2900 405}%
\special{fp}%
\special{pa 3800 405}%
\special{pa 2000 1205}%
\special{fp}%
%
\special{pn 8}%
\special{pa 2500 1205}%
\special{pa 900 405}%
\special{fp}%
\special{pa 4100 405}%
\special{pa 2500 1205}%
\special{fp}%
\special{pa 2500 1205}%
\special{pa 3200 405}%
\special{fp}%
\special{pa 1800 405}%
\special{pa 2500 1205}%
\special{fp}%
%
\special{pn 8}%
\special{pa 3000 1200}%
\special{pa 4400 400}%
\special{fp}%
\special{pa 3500 400}%
\special{pa 3000 1200}%
\special{fp}%
\special{pa 3000 1200}%
\special{pa 2100 400}%
\special{fp}%
\special{pa 1200 400}%
\special{pa 3000 1200}%
\special{fp}%
\end{picture}}%

\caption{Graph $H^{(1)}_{s}$}
\label{fH1}
\end{center}
\end{figure}

For an integer $s\geq 1$, we let $H^{(1)}_{s}$ be the graph obtained from $s$ pairwise vertex-disjoint cycles $D_{i}=u^{(i)}_{1}u^{(i)}_{2}\cdots u^{(i)}_{6}u^{(i)}_{1}~(1\leq i\leq s)$ by adding three vertices $v_{1},v_{2},v_{3}$ and the set of edges
$$
\{u^{(i)}_{j}v_{h}:1\leq i\leq s,~1\leq j\leq 6,~1\leq h\leq 3,~j\equiv h~({\rm mod~}3)\}
$$
(see Figure~\ref{fH1}).

\begin{lem}
\label{lem-nec-H^1-1}
Let $s\geq 5$ be an integer, and let $T$ be an induced subtree of $H^{(1)}_{s}$ with $\Delta (T)\leq 3$ which is not a caterpillar.
Then the following holds.
\begin{enumerate}
\item[{\upshape(i)}]
The graph $H^{(1)}_{s}$ is $P_{10}$-free.
In particular, ${\rm diam}(T)\leq 8$.
\item[{\upshape(ii)}]
If ${\rm diam}(T)=8$, then $T$ is a subgraph of $T^{*}_{9}$.
\item[{\upshape(iii)}]
If ${\rm diam}(T)=7$, then $T$ is a subgraph of $S_{8}(0,0,0,1)$ or $S_{8}^{(1)}$.
\end{enumerate}
\end{lem}
\proof
Let $P$ be an induced subpath of $H^{(1)}_{s}$ having $9$ vertices.
Then we see that either $P$ contains all of $v_{1},v_{2},v_{3}$ or, by the symmetry of vertices of $H^{(1)}_{s}$, we may assume that $P$ is equal to
$$
Q_{1}=u^{(1)}_{6}u^{(1)}_{1}v_{1}u^{(2)}_{1}u^{(2)}_{6}u^{(2)}_{5}v_{2}u^{(3)}_{2}u^{(3)}_{3}.
$$
In either case, $P$ is a maximal induced subpath of $H^{(1)}_{s}$, which implies that $H^{(1)}_{s}$ contains no induced copy of $P_{10}$.
This proves (i).

Suppose that ${\rm diam}(T)=8$, and let $P'$ be a path of $T$ having $9$ vertices.
If $P'$ contains all of $v_{1},v_{2},v_{3}$, then any induced subtree of $H^{(1)}_{s}$ containing all vertices of $P'$ is a caterpillar, which is a contradiction.
Thus we may assume that $P'=Q_{1}$.
Then $T$ does not contain $v_{3}$, and contains no three consecutive vertices of $D_{j}~(1\leq j\leq s)$ except for $u_{1}^{(2)}u_{6}^{(2)}u_{5}^{(2)}$.
This together with the assumption that $\Delta (T)\leq 3$ implies that $T$ is an induced subgraph of $T^{*}_{9}$ (for example, the subgraph of $H^{(1)}_{s}$ induced by  $\{u^{(1)}_{6},u^{(1)}_{1},v_{1},u^{(4)}_{6},u^{(4)}_{1},u^{(2)}_{1},u^{(2)}_{6},u^{(2)}_{5},v_{2},u^{(5)}_{3},u^{(5)}_{2},u^{(3)}_{2},u^{(3)}_{3}\}$ is isomorphic to $T^{*}_{9}$), which proves (ii).

Suppose that ${\rm diam}(T)=7$, and let $P''$ be a path of $T$ having $8$ vertices.
By the symmetry of vertices of $H^{(1)}_{s}$, we may assume that one of the following holds.
\begin{enumerate}
\item[{$\bullet $}]
$P''$ contains all of $v_{1},v_{2},v_{3}$;
\item[{$\bullet $}]
$P''$ is equal to $Q_{1}-u^{(3)}_{3}$, which is a subpath of $Q_{1}$; or
\item[{$\bullet $}]
$P''$ is equal to $Q_{2}=u^{(1)}_{6}u^{(1)}_{1}v_{1}u^{(2)}_{1}u^{(2)}_{2}v_{2}u^{(3)}_{2}u^{(3)}_{3}$.
\end{enumerate}
Now we focus on an induced tree $T'$ containing all vertices of $P''$ which is maximal under the condition that ${\rm diam}(T')=7$ and $\Delta (T')\leq 3$.
Then we see that
\begin{enumerate}
\item[{$\bullet $}]
if $P''$ contains all of $v_{1},v_{2},v_{3}$, then $T'$ is a caterpillar;
\item[{$\bullet $}]
if $P''=Q_{1}-u^{(3)}_{3}$ or $P''=Q_{2}$, then $T'$ does not contain $v_{3}$;
\item[{$\bullet $}]
if $P''=Q_{1}-u^{(3)}_{3}$ or $P''=Q_{2}$, then $T'$ contains no three consecutive vertices of $D_{j}~(1\leq j\leq s)$ except for vertices of $D_{2}$; and
\item[{$\bullet $}]
from the above facts, it follows that the following hold:
\begin{enumerate}
\item[{$\circ $}]
if $P''=Q_{1}-u^{(3)}_{3}$, then $T'$ is an induced copy of $S_{8}(0,0,0,1)$ (for example, the subgraph of $H^{(1)}_{s}$ induced by $\{u^{(1)}_{6},u^{(1)}_{1},v_{1},u^{(4)}_{6},u^{(4)}_{1},u^{(2)}_{1},u^{(2)}_{6},u^{(2)}_{5},v_{2},u^{(3)}_{2},u^{(3)}_{5}\}$ is isomorphic to $S_{8}(0,0,0,1)$); and
\item[{$\circ $}]
if $P''=Q_{2}$, then $T'$ is an induced copy of $S_{8}^{(1)}$ (for example, the subgraph of $H^{(1)}_{s}$ induced by $\{u^{(1)}_{6},u^{(1)}_{1},v_{1},u^{(4)}_{6},u^{(4)}_{1},u^{(2)}_{1},u^{(2)}_{6},u^{(2)}_{2},u^{(2)}_{3},v_{2},u^{(3)}_{3},u^{(3)}_{2},u^{(5)}_{2},u^{(5)}_{3}\}$ is isomorphic to $S_{8}^{(1)}$).
\end{enumerate}
\end{enumerate}
Since $T$ is a tree with ${\rm diam}(T)=7$ and $\Delta (T)\leq 3$ and  not a caterpillar, this implies that $T\prec S_{8}(0,0,0,1)$ or $T\prec S_{8}^{(1)}$, which proves (iii).
\qed

\begin{figure}
\begin{center}
{\unitlength 0.1in%
\begin{picture}(66.1000,16.9000)(2.4000,-18.5000)%
%
\special{sh 1.000}%
\special{ia 600 405 50 50 0.0000000 6.2831853}%
\special{pn 8}%
\special{ar 600 405 50 50 0.0000000 6.2831853}%
%
\special{sh 1.000}%
\special{ia 900 405 50 50 0.0000000 6.2831853}%
\special{pn 8}%
\special{ar 900 405 50 50 0.0000000 6.2831853}%
%
\special{sh 1.000}%
\special{ia 1500 405 50 50 0.0000000 6.2831853}%
\special{pn 8}%
\special{ar 1500 405 50 50 0.0000000 6.2831853}%
%
\special{sh 1.000}%
\special{ia 1200 405 50 50 0.0000000 6.2831853}%
\special{pn 8}%
\special{ar 1200 405 50 50 0.0000000 6.2831853}%
%
\special{sh 1.000}%
\special{ia 1800 405 50 50 0.0000000 6.2831853}%
\special{pn 8}%
\special{ar 1800 405 50 50 0.0000000 6.2831853}%
\put(6.0000,-2.2500){\makebox(0,0){$u^{(1)}_{1}$}}%
\put(12.0000,-2.2500){\makebox(0,0){$u^{(1)}_{3}$}}%
\put(18.0000,-2.2500){\makebox(0,0){$u^{(1)}_{5}$}}%
\put(15.0000,-2.2500){\makebox(0,0){$u^{(1)}_{4}$}}%
\put(9.0000,-2.2500){\makebox(0,0){$u^{(1)}_{2}$}}%
%
\special{pn 8}%
\special{pa 1800 400}%
\special{pa 600 400}%
\special{fp}%
%
\special{pn 4}%
\special{sh 1}%
\special{ar 3700 800 16 16 0 6.2831853}%
\special{sh 1}%
\special{ar 3600 800 16 16 0 6.2831853}%
\special{sh 1}%
\special{ar 3800 800 16 16 0 6.2831853}%
\special{sh 1}%
\special{ar 3800 800 16 16 0 6.2831853}%
%
\special{pn 8}%
\special{pa 600 400}%
\special{pa 630 415}%
\special{pa 661 430}%
\special{pa 721 460}%
\special{pa 752 474}%
\special{pa 782 488}%
\special{pa 813 502}%
\special{pa 843 515}%
\special{pa 873 527}%
\special{pa 904 539}%
\special{pa 934 550}%
\special{pa 964 560}%
\special{pa 995 569}%
\special{pa 1025 577}%
\special{pa 1055 584}%
\special{pa 1086 590}%
\special{pa 1146 598}%
\special{pa 1177 600}%
\special{pa 1207 600}%
\special{pa 1238 599}%
\special{pa 1268 596}%
\special{pa 1298 592}%
\special{pa 1329 587}%
\special{pa 1359 581}%
\special{pa 1389 573}%
\special{pa 1420 565}%
\special{pa 1480 545}%
\special{pa 1511 533}%
\special{pa 1571 509}%
\special{pa 1602 495}%
\special{pa 1632 482}%
\special{pa 1663 468}%
\special{pa 1723 438}%
\special{pa 1754 423}%
\special{pa 1800 400}%
\special{fp}%
%
\special{sh 1.000}%
\special{ia 2200 405 50 50 0.0000000 6.2831853}%
\special{pn 8}%
\special{ar 2200 405 50 50 0.0000000 6.2831853}%
%
\special{sh 1.000}%
\special{ia 2500 405 50 50 0.0000000 6.2831853}%
\special{pn 8}%
\special{ar 2500 405 50 50 0.0000000 6.2831853}%
%
\special{sh 1.000}%
\special{ia 3100 405 50 50 0.0000000 6.2831853}%
\special{pn 8}%
\special{ar 3100 405 50 50 0.0000000 6.2831853}%
%
\special{sh 1.000}%
\special{ia 2800 405 50 50 0.0000000 6.2831853}%
\special{pn 8}%
\special{ar 2800 405 50 50 0.0000000 6.2831853}%
%
\special{sh 1.000}%
\special{ia 3400 405 50 50 0.0000000 6.2831853}%
\special{pn 8}%
\special{ar 3400 405 50 50 0.0000000 6.2831853}%
\put(22.0000,-2.2500){\makebox(0,0){$v^{(1)}_{1}$}}%
\put(28.0000,-2.2500){\makebox(0,0){$v^{(1)}_{3}$}}%
\put(34.0000,-2.2500){\makebox(0,0){$v^{(1)}_{5}$}}%
\put(31.0000,-2.2500){\makebox(0,0){$v^{(1)}_{4}$}}%
\put(25.0000,-2.2500){\makebox(0,0){$v^{(1)}_{2}$}}%
%
\special{pn 8}%
\special{pa 3400 400}%
\special{pa 2200 400}%
\special{fp}%
%
\special{pn 8}%
\special{pa 2200 400}%
\special{pa 2230 415}%
\special{pa 2261 430}%
\special{pa 2321 460}%
\special{pa 2352 474}%
\special{pa 2382 488}%
\special{pa 2413 502}%
\special{pa 2443 515}%
\special{pa 2473 527}%
\special{pa 2504 539}%
\special{pa 2534 550}%
\special{pa 2564 560}%
\special{pa 2595 569}%
\special{pa 2625 577}%
\special{pa 2655 584}%
\special{pa 2686 590}%
\special{pa 2746 598}%
\special{pa 2777 600}%
\special{pa 2807 600}%
\special{pa 2838 599}%
\special{pa 2868 596}%
\special{pa 2898 592}%
\special{pa 2929 587}%
\special{pa 2959 581}%
\special{pa 2989 573}%
\special{pa 3020 565}%
\special{pa 3080 545}%
\special{pa 3111 533}%
\special{pa 3171 509}%
\special{pa 3202 495}%
\special{pa 3232 482}%
\special{pa 3263 468}%
\special{pa 3323 438}%
\special{pa 3354 423}%
\special{pa 3400 400}%
\special{fp}%
%
\special{sh 1.000}%
\special{ia 2000 1200 50 50 0.0000000 6.2831853}%
\special{pn 8}%
\special{ar 2000 1200 50 50 0.0000000 6.2831853}%
%
\special{sh 1.000}%
\special{ia 1700 1200 50 50 0.0000000 6.2831853}%
\special{pn 8}%
\special{ar 1700 1200 50 50 0.0000000 6.2831853}%
%
\special{sh 1.000}%
\special{ia 1400 1200 50 50 0.0000000 6.2831853}%
\special{pn 8}%
\special{ar 1400 1200 50 50 0.0000000 6.2831853}%
%
\special{sh 1.000}%
\special{ia 2600 1200 50 50 0.0000000 6.2831853}%
\special{pn 8}%
\special{ar 2600 1200 50 50 0.0000000 6.2831853}%
%
\special{sh 1.000}%
\special{ia 2300 1200 50 50 0.0000000 6.2831853}%
\special{pn 8}%
\special{ar 2300 1200 50 50 0.0000000 6.2831853}%
%
\special{sh 1.000}%
\special{ia 4000 405 50 50 0.0000000 6.2831853}%
\special{pn 8}%
\special{ar 4000 405 50 50 0.0000000 6.2831853}%
%
\special{sh 1.000}%
\special{ia 4300 405 50 50 0.0000000 6.2831853}%
\special{pn 8}%
\special{ar 4300 405 50 50 0.0000000 6.2831853}%
%
\special{sh 1.000}%
\special{ia 4900 405 50 50 0.0000000 6.2831853}%
\special{pn 8}%
\special{ar 4900 405 50 50 0.0000000 6.2831853}%
%
\special{sh 1.000}%
\special{ia 4600 405 50 50 0.0000000 6.2831853}%
\special{pn 8}%
\special{ar 4600 405 50 50 0.0000000 6.2831853}%
%
\special{sh 1.000}%
\special{ia 5200 405 50 50 0.0000000 6.2831853}%
\special{pn 8}%
\special{ar 5200 405 50 50 0.0000000 6.2831853}%
\put(40.0000,-2.2500){\makebox(0,0){$u^{(s)}_{1}$}}%
\put(46.0000,-2.2500){\makebox(0,0){$u^{(s)}_{3}$}}%
\put(52.0000,-2.2500){\makebox(0,0){$u^{(s)}_{5}$}}%
\put(49.0000,-2.2500){\makebox(0,0){$u^{(s)}_{4}$}}%
\put(43.0000,-2.2500){\makebox(0,0){$u^{(s)}_{2}$}}%
%
\special{pn 8}%
\special{pa 5200 400}%
\special{pa 4000 400}%
\special{fp}%
%
\special{pn 8}%
\special{pa 4000 400}%
\special{pa 4030 415}%
\special{pa 4061 430}%
\special{pa 4121 460}%
\special{pa 4152 474}%
\special{pa 4182 488}%
\special{pa 4213 502}%
\special{pa 4243 515}%
\special{pa 4273 527}%
\special{pa 4304 539}%
\special{pa 4334 550}%
\special{pa 4364 560}%
\special{pa 4395 569}%
\special{pa 4425 577}%
\special{pa 4455 584}%
\special{pa 4486 590}%
\special{pa 4546 598}%
\special{pa 4577 600}%
\special{pa 4607 600}%
\special{pa 4638 599}%
\special{pa 4668 596}%
\special{pa 4698 592}%
\special{pa 4729 587}%
\special{pa 4759 581}%
\special{pa 4789 573}%
\special{pa 4820 565}%
\special{pa 4880 545}%
\special{pa 4911 533}%
\special{pa 4971 509}%
\special{pa 5002 495}%
\special{pa 5032 482}%
\special{pa 5063 468}%
\special{pa 5123 438}%
\special{pa 5154 423}%
\special{pa 5200 400}%
\special{fp}%
%
\special{sh 1.000}%
\special{ia 5600 405 50 50 0.0000000 6.2831853}%
\special{pn 8}%
\special{ar 5600 405 50 50 0.0000000 6.2831853}%
%
\special{sh 1.000}%
\special{ia 5900 405 50 50 0.0000000 6.2831853}%
\special{pn 8}%
\special{ar 5900 405 50 50 0.0000000 6.2831853}%
%
\special{sh 1.000}%
\special{ia 6500 405 50 50 0.0000000 6.2831853}%
\special{pn 8}%
\special{ar 6500 405 50 50 0.0000000 6.2831853}%
%
\special{sh 1.000}%
\special{ia 6200 405 50 50 0.0000000 6.2831853}%
\special{pn 8}%
\special{ar 6200 405 50 50 0.0000000 6.2831853}%
%
\special{sh 1.000}%
\special{ia 6800 405 50 50 0.0000000 6.2831853}%
\special{pn 8}%
\special{ar 6800 405 50 50 0.0000000 6.2831853}%
\put(56.0000,-2.2500){\makebox(0,0){$v^{(s)}_{1}$}}%
\put(62.0000,-2.2500){\makebox(0,0){$v^{(s)}_{3}$}}%
\put(68.0000,-2.2500){\makebox(0,0){$v^{(s)}_{5}$}}%
\put(65.0000,-2.2500){\makebox(0,0){$v^{(s)}_{4}$}}%
\put(59.0000,-2.2500){\makebox(0,0){$v^{(s)}_{2}$}}%
%
\special{pn 8}%
\special{pa 6800 400}%
\special{pa 5600 400}%
\special{fp}%
%
\special{pn 8}%
\special{pa 5600 400}%
\special{pa 5630 415}%
\special{pa 5661 430}%
\special{pa 5721 460}%
\special{pa 5752 474}%
\special{pa 5782 488}%
\special{pa 5813 502}%
\special{pa 5843 515}%
\special{pa 5873 527}%
\special{pa 5904 539}%
\special{pa 5934 550}%
\special{pa 5964 560}%
\special{pa 5995 569}%
\special{pa 6025 577}%
\special{pa 6055 584}%
\special{pa 6086 590}%
\special{pa 6146 598}%
\special{pa 6177 600}%
\special{pa 6207 600}%
\special{pa 6238 599}%
\special{pa 6268 596}%
\special{pa 6298 592}%
\special{pa 6329 587}%
\special{pa 6359 581}%
\special{pa 6389 573}%
\special{pa 6420 565}%
\special{pa 6480 545}%
\special{pa 6511 533}%
\special{pa 6571 509}%
\special{pa 6602 495}%
\special{pa 6632 482}%
\special{pa 6663 468}%
\special{pa 6723 438}%
\special{pa 6754 423}%
\special{pa 6800 400}%
\special{fp}%
\put(13.9000,-10.2000){\makebox(0,0){$w^{(1)}_{1}$}}%
\put(20.2000,-10.2000){\makebox(0,0){$w^{(1)}_{3}$}}%
\put(26.3000,-10.2000){\makebox(0,0){$w^{(1)}_{5}$}}%
\put(23.2000,-10.2000){\makebox(0,0){$w^{(1)}_{4}$}}%
\put(17.2000,-10.2000){\makebox(0,0){$w^{(1)}_{2}$}}%
%
\special{pn 8}%
\special{pa 1400 1205}%
\special{pa 600 405}%
\special{fp}%
%
\special{pn 8}%
\special{pa 2000 1200}%
\special{pa 1200 400}%
\special{fp}%
%
\special{pn 8}%
\special{pa 2600 1200}%
\special{pa 1800 400}%
\special{fp}%
%
\special{pn 8}%
\special{pa 2300 1200}%
\special{pa 1500 400}%
\special{fp}%
%
\special{pn 8}%
\special{pa 1700 1200}%
\special{pa 900 400}%
\special{fp}%
%
\special{pn 8}%
\special{pa 2600 1205}%
\special{pa 3400 405}%
\special{fp}%
%
\special{pn 8}%
\special{pa 2000 1200}%
\special{pa 2800 400}%
\special{fp}%
%
\special{pn 8}%
\special{pa 1400 1200}%
\special{pa 2200 400}%
\special{fp}%
%
\special{pn 8}%
\special{pa 1700 1200}%
\special{pa 2500 400}%
\special{fp}%
%
\special{pn 8}%
\special{pa 2300 1200}%
\special{pa 3100 400}%
\special{fp}%
%
\special{sh 1.000}%
\special{ia 5400 1200 50 50 0.0000000 6.2831853}%
\special{pn 8}%
\special{ar 5400 1200 50 50 0.0000000 6.2831853}%
%
\special{sh 1.000}%
\special{ia 5100 1200 50 50 0.0000000 6.2831853}%
\special{pn 8}%
\special{ar 5100 1200 50 50 0.0000000 6.2831853}%
%
\special{sh 1.000}%
\special{ia 4800 1200 50 50 0.0000000 6.2831853}%
\special{pn 8}%
\special{ar 4800 1200 50 50 0.0000000 6.2831853}%
%
\special{sh 1.000}%
\special{ia 6000 1200 50 50 0.0000000 6.2831853}%
\special{pn 8}%
\special{ar 6000 1200 50 50 0.0000000 6.2831853}%
%
\special{sh 1.000}%
\special{ia 5700 1200 50 50 0.0000000 6.2831853}%
\special{pn 8}%
\special{ar 5700 1200 50 50 0.0000000 6.2831853}%
\put(48.0000,-10.2000){\makebox(0,0){$w^{(s)}_{1}$}}%
\put(54.2000,-10.2000){\makebox(0,0){$w^{(s)}_{3}$}}%
\put(60.5000,-10.2000){\makebox(0,0){$w^{(s)}_{5}$}}%
\put(57.2000,-10.2000){\makebox(0,0){$w^{(s)}_{4}$}}%
\put(51.2000,-10.2000){\makebox(0,0){$w^{(s)}_{2}$}}%
%
\special{pn 8}%
\special{pa 4800 1205}%
\special{pa 4000 405}%
\special{fp}%
%
\special{pn 8}%
\special{pa 5400 1200}%
\special{pa 4600 400}%
\special{fp}%
%
\special{pn 8}%
\special{pa 6000 1200}%
\special{pa 5200 400}%
\special{fp}%
%
\special{pn 8}%
\special{pa 5700 1200}%
\special{pa 4900 400}%
\special{fp}%
%
\special{pn 8}%
\special{pa 5100 1200}%
\special{pa 4300 400}%
\special{fp}%
%
\special{pn 8}%
\special{pa 6000 1205}%
\special{pa 6800 405}%
\special{fp}%
%
\special{pn 8}%
\special{pa 5400 1200}%
\special{pa 6200 400}%
\special{fp}%
%
\special{pn 8}%
\special{pa 4800 1200}%
\special{pa 5600 400}%
\special{fp}%
%
\special{pn 8}%
\special{pa 5100 1200}%
\special{pa 5900 400}%
\special{fp}%
%
\special{pn 8}%
\special{pa 5700 1200}%
\special{pa 6500 400}%
\special{fp}%
%
\special{sh 1.000}%
\special{ia 3700 1800 50 50 0.0000000 6.2831853}%
\special{pn 8}%
\special{ar 3700 1800 50 50 0.0000000 6.2831853}%
%
\special{pn 8}%
\special{pa 3700 1800}%
\special{pa 1400 1200}%
\special{fp}%
%
\special{pn 8}%
\special{pa 3700 1800}%
\special{pa 1700 1200}%
\special{fp}%
\special{pa 2300 1200}%
\special{pa 3700 1800}%
\special{fp}%
%
\special{pn 8}%
\special{pa 3700 1800}%
\special{pa 2000 1200}%
\special{fp}%
\special{pa 2600 1200}%
\special{pa 3700 1800}%
\special{fp}%
%
\special{pn 8}%
\special{pa 3700 1800}%
\special{pa 6000 1200}%
\special{fp}%
%
\special{pn 8}%
\special{pa 3700 1800}%
\special{pa 5700 1200}%
\special{fp}%
\special{pa 5100 1200}%
\special{pa 3700 1800}%
\special{fp}%
%
\special{pn 8}%
\special{pa 3700 1800}%
\special{pa 5400 1200}%
\special{fp}%
\special{pa 4800 1200}%
\special{pa 3700 1800}%
\special{fp}%
\put(37.0000,-16.7000){\makebox(0,0){$z$}}%
\end{picture}}%

\caption{Graph $H^{(2)}_{s}$}
\label{fH2}
\end{center}
\end{figure}
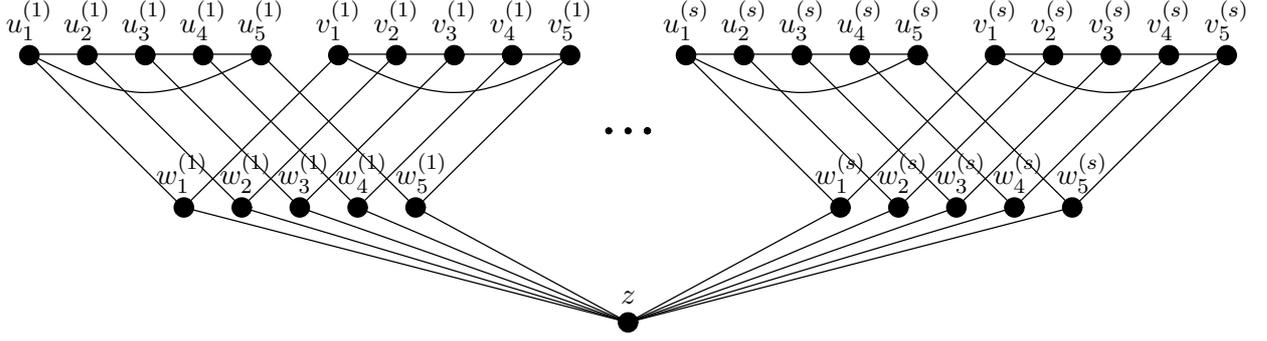

Let $A$ be the graph obtained from two vertex-disjoint cycles $u_{1}u_{2}\cdots u_{5}u_{1}$ and $v_{1}v_{2}\cdots v_{5}v_{1}$ by adding five vertices $w_{1},\ldots ,w_{5}$ and edges $u_{j}w_{j},v_{j}w_{j}~(1\leq j\leq 5)$.
For an integer $s\geq 3$, we let $H^{(2)}_{s}$ be the graph obtained from $s$ pairwise vertex-disjoint copies $A_{1},\ldots ,A_{s}$ of $A$ such that $V(A_{i})=\{u^{(i)}_{j},v^{(i)}_{j},w^{(i)}_{j}:1\leq j\leq 5\}$ where $u^{(i)}_{j}$, $v^{(i)}_{j}$ and $w^{(i)}_{j}$ respectively correspond to $u_{j}$, $v_{j}$ and $w_{j}$, by adding a new vertex $z$ and edges $zw^{(i)}_{j}~(1\leq i\leq s,~1\leq j\leq 5)$ (see Figure~\ref{fH2}).

\begin{lem}
\label{lem-nec-H^2-1}
For an integer $s\geq 3$, $H^{(2)}_{s}$ is $\{S_{8}(0,0,0,1),T^{*}_{8}\}$-free.
\end{lem}
\proof
Let $T$ be an induced subgraph of $H^{(2)}_{s}$ isomorphic to $T_{8}$.
It suffices to show that
\begin{align}
\mbox{no induced subgraph of $H^{(2)}_{s}$ isomorphic to $S_{8}(0,0,0,1)$ or $T^{*}_{8}$ contains all vertices of $T$.}\label{cond-lem-nec-H^2-1-1}
\end{align}
Let $x$ be the unique vertex of $T$ with $d_{T}(x)=3$.
Without loss of generality, we may assume that $x\in \{u^{(1)}_{1},w^{(1)}_{1},z\}$.

We first assume that $x=u^{(1)}_{1}$.
Then we easily verify that $N_{T}(u^{(1)}_{1})\cap \{w^{(1)}_{j}:1\leq j\leq 5\}\neq \emptyset $ and $N_{T}(N_{T}(u^{(1)}_{1}))\cap \{w^{(1)}_{j}:1\leq j\leq 5\}\neq \emptyset $, and hence $z\notin V(T)$, i.e., $V(T)\subseteq V(A_{1})$, and by the symmetry of vertices of $A_{1}$, we may assume that $T$ is a subgraph of $A_{1}$ induced by
\begin{enumerate}
\item[{$\bullet $}]
$\{w^{(1)}_{2},u^{(1)}_{2},u^{(1)}_{1},u^{(1)}_{4},u^{(1)}_{5},w^{(1)}_{1},v^{(1)}_{1},v^{(1)}_{5},v^{(1)}_{4},v^{(1)}_{3}\}$; or
\item[{$\bullet $}]
$\{w^{(1)}_{2},u^{(1)}_{2},u^{(1)}_{1},v^{(1)}_{1},w^{(1)}_{1},u^{(1)}_{5},u^{(1)}_{4},w^{(1)}_{4},v^{(1)}_{4},v^{(1)}_{3}\}$.
\end{enumerate}
In either case, (\ref{cond-lem-nec-H^2-1-1}) holds.

Next we assume that $x=w^{(1)}_{1}$.
Since $d_{T}(x)=3$, we have $z\in V(T)$.
Furthermore, we see that the path $P$ of $T$ with an end $w^{(1)}_{1}$ having length $5$ intersects with $\bigcup _{2\leq i\leq s}V(A_{i})$.
By the symmetry of indices $2,3,\ldots ,s$, we may assume that $V(P)\cap V(A_{2})\neq \emptyset $ and, by the symmetry of vertices of $A_{2}$, we may further assume that $P=w^{(1)}_{1}zw^{(2)}_{1}u^{(2)}_{1}u^{(2)}_{2}u^{(2)}_{3}$.
Suppose that there is an induced subgraph $S$ of $H^{(2)}_{s}$ isomorphic to $S_{8}(0,0,0,1)$ or $T^{*}_{8}$ containing all vertices of $T$.
Note that $N_{H^{(2)}_{s}}(u^{(2)}_{2})-V(T)=\{w^{(2)}_{2}\}$, $N_{H^{(2)}_{s}}(u^{(2)}_{1})-V(T)=\{u^{(2)}_{5}\}$ and $N_{H^{(2)}_{s}}(u^{(2)}_{5})-V(T)=\{u^{(2)}_{4},w^{(2)}_{5}\}$.
Hence
\begin{enumerate}
\item[{$\bullet $}]
if $S\simeq S_{8}(0,0,0,1)$, then $w^{(2)}_{2}\in V(S)$; and
\item[{$\bullet $}]
if $S\simeq T^{*}_{8}$, then we have $\{u^{(2)}_{5},u^{(2)}_{4}\}\subseteq V(S)$ or $\{u^{(2)}_{5},w^{(2)}_{5}\}\subseteq V(S)$.
\end{enumerate}
In either case, $S$ contains a cycle.
This contradicts the fact that $S$ is a tree, which shows that (\ref{cond-lem-nec-H^2-1-1}) holds.

Finally we assume that $x=z$.
Let $P$ be the path of $T$ with an end $z$ having length $5$.
Then by the symmetry of vertices of $H^{(2)}_{s}$, we may assume that $P=zw^{(1)}_{1}u^{(1)}_{1}u^{(1)}_{2}u^{(1)}_{3}u^{(1)}_{4}$.
Since $N_{H^{(2)}_{s}}(u^{(1)}_{j})-V(T)=\{w^{(1)}_{j}\}$ for each $j\in \{2,3\}$,
\begin{enumerate}
\item[{$\bullet $}]
if an induced subgraph $S$ of $H^{(2)}_{s}$ isomorphic to $S_{8}(0,0,0,1)$ contains all vertices of $T$, then $w^{(1)}_{3}\in V(S)$; and
\item[{$\bullet $}]
if an induced subgraph $S$ of $H^{(2)}_{s}$ isomorphic to $T^{*}_{8}$ contains all vertices of $T$, then $w^{(1)}_{2}\in V(S)$.
\end{enumerate}
In either case, $S$ contains a cycle, which is a contradiction.
This leads to (\ref{cond-lem-nec-H^2-1-1}).
\qed

Let $A'$ be the graph depicted in Figure~\ref{fH3}.
For an integer $s\geq 4$, we let $H^{(3)}_{s}$ be the graph obtained from $s$ pairwise vertex-disjoint copies $A'_{1},\ldots ,A'_{s}$ of $A'$ such that $V(A'_{i})=\{u^{(i)}_{j},v^{(i)}_{j,h},w^{(i)}_{j,h'}:1\leq j\leq 2,~1\leq h\leq 2,~1\leq h'\leq 4\}$ where $u^{(i)}_{j}$, $v^{(i)}_{j,h}$ and $w^{(i)}_{j,h'}$ respectively correspond to $u_{j}$, $v_{j,h}$ and $w_{j,h'}$, by adding edges $u^{(i)}_{2}u^{(i+1)}_{1}$ (indices are to be read modulo $s$).

\begin{figure}
\begin{center}
{\unitlength 0.1in%
\begin{picture}(24.5000,15.0000)(2.0000,-16.8000)%
%
\special{sh 1.000}%
\special{ia 400 995 50 50 0.0000000 6.2831853}%
\special{pn 8}%
\special{ar 400 995 50 50 0.0000000 6.2831853}%
%
\special{sh 1.000}%
\special{ia 800 595 50 50 0.0000000 6.2831853}%
\special{pn 8}%
\special{ar 800 595 50 50 0.0000000 6.2831853}%
%
\special{sh 1.000}%
\special{ia 1200 395 50 50 0.0000000 6.2831853}%
\special{pn 8}%
\special{ar 1200 395 50 50 0.0000000 6.2831853}%
%
\special{sh 1.000}%
\special{ia 1200 795 50 50 0.0000000 6.2831853}%
\special{pn 8}%
\special{ar 1200 795 50 50 0.0000000 6.2831853}%
%
\special{sh 1.000}%
\special{ia 1200 1195 50 50 0.0000000 6.2831853}%
\special{pn 8}%
\special{ar 1200 1195 50 50 0.0000000 6.2831853}%
%
\special{sh 1.000}%
\special{ia 1200 1595 50 50 0.0000000 6.2831853}%
\special{pn 8}%
\special{ar 1200 1595 50 50 0.0000000 6.2831853}%
%
\special{sh 1.000}%
\special{ia 800 1395 50 50 0.0000000 6.2831853}%
\special{pn 8}%
\special{ar 800 1395 50 50 0.0000000 6.2831853}%
%
\special{pn 8}%
\special{pa 400 995}%
\special{pa 800 595}%
\special{fp}%
\special{pa 800 595}%
\special{pa 1200 395}%
\special{fp}%
\special{pa 1200 1595}%
\special{pa 800 1395}%
\special{fp}%
\special{pa 800 1395}%
\special{pa 1200 1195}%
\special{fp}%
\special{pa 1200 795}%
\special{pa 800 595}%
\special{fp}%
\special{pa 800 1395}%
\special{pa 400 995}%
\special{fp}%
%
\special{sh 1.000}%
\special{ia 2600 995 50 50 0.0000000 6.2831853}%
\special{pn 8}%
\special{ar 2600 995 50 50 0.0000000 6.2831853}%
%
\special{sh 1.000}%
\special{ia 2200 595 50 50 0.0000000 6.2831853}%
\special{pn 8}%
\special{ar 2200 595 50 50 0.0000000 6.2831853}%
%
\special{sh 1.000}%
\special{ia 1800 395 50 50 0.0000000 6.2831853}%
\special{pn 8}%
\special{ar 1800 395 50 50 0.0000000 6.2831853}%
%
\special{sh 1.000}%
\special{ia 1800 795 50 50 0.0000000 6.2831853}%
\special{pn 8}%
\special{ar 1800 795 50 50 0.0000000 6.2831853}%
%
\special{sh 1.000}%
\special{ia 1800 1195 50 50 0.0000000 6.2831853}%
\special{pn 8}%
\special{ar 1800 1195 50 50 0.0000000 6.2831853}%
%
\special{sh 1.000}%
\special{ia 1800 1595 50 50 0.0000000 6.2831853}%
\special{pn 8}%
\special{ar 1800 1595 50 50 0.0000000 6.2831853}%
%
\special{sh 1.000}%
\special{ia 2200 1395 50 50 0.0000000 6.2831853}%
\special{pn 8}%
\special{ar 2200 1395 50 50 0.0000000 6.2831853}%
%
\special{pn 8}%
\special{pa 2600 995}%
\special{pa 2200 595}%
\special{fp}%
\special{pa 2200 595}%
\special{pa 1800 395}%
\special{fp}%
\special{pa 1800 1595}%
\special{pa 2200 1395}%
\special{fp}%
\special{pa 2200 1395}%
\special{pa 1800 1195}%
\special{fp}%
\special{pa 1800 795}%
\special{pa 2200 595}%
\special{fp}%
\special{pa 2200 1395}%
\special{pa 2600 995}%
\special{fp}%
%
\special{pn 8}%
\special{pa 1800 395}%
\special{pa 1200 395}%
\special{fp}%
\special{pa 1200 795}%
\special{pa 1800 1195}%
\special{fp}%
\special{pa 1200 1195}%
\special{pa 1800 795}%
\special{fp}%
\special{pa 1800 1595}%
\special{pa 1200 1595}%
\special{fp}%
%
\special{pn 8}%
\special{pa 1200 1595}%
\special{pa 1179 1566}%
\special{pa 1157 1538}%
\special{pa 1137 1509}%
\special{pa 1116 1481}%
\special{pa 1097 1452}%
\special{pa 1079 1423}%
\special{pa 1062 1395}%
\special{pa 1047 1366}%
\special{pa 1034 1337}%
\special{pa 1022 1309}%
\special{pa 1013 1280}%
\special{pa 1006 1252}%
\special{pa 1001 1223}%
\special{pa 1000 1194}%
\special{pa 1002 1166}%
\special{pa 1006 1137}%
\special{pa 1013 1108}%
\special{pa 1022 1080}%
\special{pa 1034 1051}%
\special{pa 1048 1023}%
\special{pa 1063 994}%
\special{pa 1080 965}%
\special{pa 1098 937}%
\special{pa 1117 908}%
\special{pa 1138 879}%
\special{pa 1158 851}%
\special{pa 1180 822}%
\special{pa 1200 795}%
\special{fp}%
%
\special{pn 8}%
\special{pa 1200 1195}%
\special{pa 1179 1166}%
\special{pa 1157 1138}%
\special{pa 1137 1109}%
\special{pa 1116 1081}%
\special{pa 1097 1052}%
\special{pa 1079 1023}%
\special{pa 1062 995}%
\special{pa 1047 966}%
\special{pa 1034 937}%
\special{pa 1022 909}%
\special{pa 1013 880}%
\special{pa 1006 852}%
\special{pa 1001 823}%
\special{pa 1000 794}%
\special{pa 1002 766}%
\special{pa 1006 737}%
\special{pa 1013 708}%
\special{pa 1022 680}%
\special{pa 1034 651}%
\special{pa 1048 623}%
\special{pa 1063 594}%
\special{pa 1080 565}%
\special{pa 1098 537}%
\special{pa 1117 508}%
\special{pa 1138 479}%
\special{pa 1158 451}%
\special{pa 1180 422}%
\special{pa 1200 395}%
\special{fp}%
%
\special{pn 8}%
\special{pa 1800 1195}%
\special{pa 1821 1166}%
\special{pa 1843 1138}%
\special{pa 1863 1109}%
\special{pa 1884 1081}%
\special{pa 1903 1052}%
\special{pa 1921 1023}%
\special{pa 1938 995}%
\special{pa 1953 966}%
\special{pa 1966 937}%
\special{pa 1978 909}%
\special{pa 1987 880}%
\special{pa 1994 852}%
\special{pa 1999 823}%
\special{pa 2000 794}%
\special{pa 1998 766}%
\special{pa 1994 737}%
\special{pa 1987 708}%
\special{pa 1978 680}%
\special{pa 1966 651}%
\special{pa 1952 623}%
\special{pa 1937 594}%
\special{pa 1920 565}%
\special{pa 1902 537}%
\special{pa 1883 508}%
\special{pa 1862 479}%
\special{pa 1842 451}%
\special{pa 1820 422}%
\special{pa 1800 395}%
\special{fp}%
%
\special{pn 8}%
\special{pa 1800 1595}%
\special{pa 1821 1566}%
\special{pa 1843 1538}%
\special{pa 1863 1509}%
\special{pa 1884 1481}%
\special{pa 1903 1452}%
\special{pa 1921 1423}%
\special{pa 1938 1395}%
\special{pa 1953 1366}%
\special{pa 1966 1337}%
\special{pa 1978 1309}%
\special{pa 1987 1280}%
\special{pa 1994 1252}%
\special{pa 1999 1223}%
\special{pa 2000 1194}%
\special{pa 1998 1166}%
\special{pa 1994 1137}%
\special{pa 1987 1108}%
\special{pa 1978 1080}%
\special{pa 1966 1051}%
\special{pa 1952 1023}%
\special{pa 1937 994}%
\special{pa 1920 965}%
\special{pa 1902 937}%
\special{pa 1883 908}%
\special{pa 1862 879}%
\special{pa 1842 851}%
\special{pa 1820 822}%
\special{pa 1800 795}%
\special{fp}%
\put(4.0000,-8.5000){\makebox(0,0){$u_{1}$}}%
\put(26.0000,-8.5000){\makebox(0,0){$u_{2}$}}%
\put(8.0000,-4.4500){\makebox(0,0){$v_{1,1}$}}%
\put(8.0000,-12.4500){\makebox(0,0){$v_{1,2}$}}%
\put(22.0000,-12.4500){\makebox(0,0){$v_{2,2}$}}%
\put(22.0000,-4.4500){\makebox(0,0){$v_{2,1}$}}%
\put(12.0000,-2.4500){\makebox(0,0){$w_{1,1}$}}%
\put(12.0000,-6.4500){\makebox(0,0){$w_{1,2}$}}%
\put(12.0000,-13.4500){\makebox(0,0){$w_{1,3}$}}%
\put(12.0000,-17.4500){\makebox(0,0){$w_{1,4}$}}%
\put(18.0000,-17.4500){\makebox(0,0){$w_{2,4}$}}%
\put(18.0000,-13.4500){\makebox(0,0){$w_{2,3}$}}%
\put(18.0000,-6.4500){\makebox(0,0){$w_{2,2}$}}%
\put(18.0000,-2.4500){\makebox(0,0){$w_{2,1}$}}%
\end{picture}}%

\caption{Graph $A'$}
\label{fH3}
\end{center}
\end{figure}

\begin{lem}
\label{lem-nec-H^3-1}
For an integer $s\geq 4$, $H^{(3)}_{s}$ is $S_{7}(1,0,1)$-free.
\end{lem}
\proof
Suppose that $H^{(3)}_{s}$ contains an induced subgraph $T$ which is isomorphic to $S_{7}(1,0,1)$.
Let $S$ be the subtree of $T$ which is isomorphic to $S_{7}(1,0,0)$.
Let $x$ be the unique vertex of $S$ with $d_{S}(x)=3$ such that no leaf of $S$ is adjacent to $x$.
Without loss of generality, we may assume that $x\in \{w^{(1)}_{1,1},v^{(1)}_{1,1},u^{(1)}_{1}\}$.

We show that $x\neq w^{(1)}_{1,1}$.
Let $S'$ be an induced subtree of $H^{(3)}_{s}$ isomorphic to $T_{5}$ such that $w^{(1)}_{1,1}$ is the unique vertex of $S'$ having degree $3$.
Then we see that $S'$ is equal to either
$$
H^{(3)}_{s}[\{u^{(1)}_{1},v^{(1)}_{1,1},w^{(1)}_{1,1},w^{(1)}_{2,2},w^{(1)}_{1,3},w^{(1)}_{2,1},w^{(1)}_{2,3}\}]\mbox{~~or~~}H^{(3)}_{s}[\{w^{(1)}_{1,2},v^{(1)}_{1,1},w^{(1)}_{1,1},v^{(1)}_{1,2},w^{(1)}_{1,3},w^{(1)}_{2,1},v^{(1)}_{2,1}\}].
$$
In either case, we can verify that no induced subgraph of $H^{(3)}_{s}$ isomorphic to $S_{5}(1)$ contains all vertices of $S'$.
Since $S_{5}(1)\prec S_{7}(1,0,0)$, this implies $x\neq w^{(1)}_{1,1}$.
Thus $x\in \{v^{(1)}_{1,1},u^{(1)}_{1}\}$.

We make the following two observations.
\begin{enumerate}
\item[{$\bullet $}]
Assume that $x=v^{(1)}_{1,1}$.
Then by the symmetry of vertices of $A'_{1}$, we may assume that
$$
S=H^{(3)}_{s}[\{u^{(s)}_{2},u^{(1)}_{1},v^{(1)}_{1,1},w^{(1)}_{1,3},w^{(1)}_{1,1},w^{(1)}_{1,2},w^{(1)}_{1,4},w^{(1)}_{2,3},v^{(1)}_{2,2},u^{(1)}_{2}\}].
$$
Note that $N_{H^{(3)}_{s}}(v^{(1)}_{2,2})-V(S)=\{w^{(1)}_{2,4}\}$.
\item[{$\bullet $}]
Assume that $x=u^{(1)}_{1}$.
Then by the symmetry of vertices of $A'_{1}$ and the symmetry of vertices of $A'_{s}$, we may assume that
$$
S=H^{(3)}_{s}[\{w^{(1)}_{1,1},v^{(1)}_{1,1},u^{(1)}_{1},w^{(1)}_{1,4},v^{(1)}_{1,2},u^{(s)}_{2},v^{(s)}_{2,2},v^{(s)}_{2,1},w^{(s)}_{2,1},w^{(s)}_{1,1}\}].
$$
Note that $N_{H^{(3)}_{s}}(w^{(s)}_{2,1})-V(S)=\{w^{(s)}_{2,3}\}$.
\end{enumerate}
If $x=v^{(1)}_{1,1}$, let $y=w^{(1)}_{2,4}$; if $x=u^{(1)}_{1}$, let $y=w^{(s)}_{2,3}$.
Then in either case, it follows from the above observations that $V(T)-V(S)=\{y\}$, which forces $T$ to contain a cycle, a contradiction.
\qed

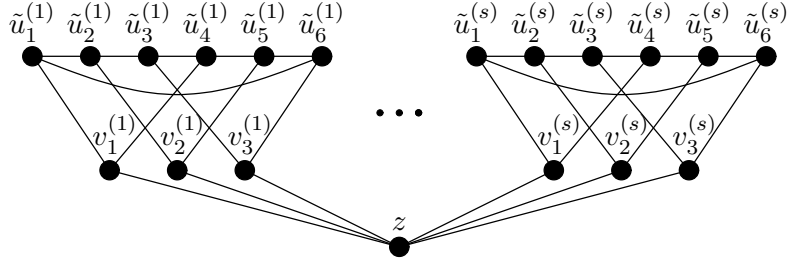
\begin{figure}
\begin{center}
{\unitlength 0.1in%
\begin{picture}(43.9000,12.9000)(0.6000,-14.5000)%
%
\special{sh 1.000}%
\special{ia 600 405 50 50 0.0000000 6.2831853}%
\special{pn 8}%
\special{ar 600 405 50 50 0.0000000 6.2831853}%
%
\special{sh 1.000}%
\special{ia 900 405 50 50 0.0000000 6.2831853}%
\special{pn 8}%
\special{ar 900 405 50 50 0.0000000 6.2831853}%
%
\special{sh 1.000}%
\special{ia 1500 405 50 50 0.0000000 6.2831853}%
\special{pn 8}%
\special{ar 1500 405 50 50 0.0000000 6.2831853}%
%
\special{sh 1.000}%
\special{ia 1200 405 50 50 0.0000000 6.2831853}%
\special{pn 8}%
\special{ar 1200 405 50 50 0.0000000 6.2831853}%
%
\special{sh 1.000}%
\special{ia 2100 405 50 50 0.0000000 6.2831853}%
\special{pn 8}%
\special{ar 2100 405 50 50 0.0000000 6.2831853}%
%
\special{sh 1.000}%
\special{ia 1800 405 50 50 0.0000000 6.2831853}%
\special{pn 8}%
\special{ar 1800 405 50 50 0.0000000 6.2831853}%
\put(6.0000,-2.2500){\makebox(0,0){$\tilde{u}^{(1)}_{1}$}}%
\put(12.0000,-2.2500){\makebox(0,0){$\tilde{u}^{(1)}_{3}$}}%
\put(18.0000,-2.2500){\makebox(0,0){$\tilde{u}^{(1)}_{5}$}}%
\put(21.0000,-2.2500){\makebox(0,0){$\tilde{u}^{(1)}_{6}$}}%
\put(15.0000,-2.2500){\makebox(0,0){$\tilde{u}^{(1)}_{4}$}}%
\put(9.0000,-2.2500){\makebox(0,0){$\tilde{u}^{(1)}_{2}$}}%
%
\special{pn 8}%
\special{pa 600 405}%
\special{pa 631 417}%
\special{pa 662 430}%
\special{pa 786 478}%
\special{pa 816 489}%
\special{pa 909 522}%
\special{pa 940 532}%
\special{pa 1033 559}%
\special{pa 1064 567}%
\special{pa 1126 581}%
\special{pa 1157 587}%
\special{pa 1187 592}%
\special{pa 1249 600}%
\special{pa 1311 604}%
\special{pa 1342 605}%
\special{pa 1373 605}%
\special{pa 1435 601}%
\special{pa 1466 598}%
\special{pa 1528 590}%
\special{pa 1559 584}%
\special{pa 1589 578}%
\special{pa 1620 571}%
\special{pa 1682 555}%
\special{pa 1744 537}%
\special{pa 1806 517}%
\special{pa 1899 484}%
\special{pa 1930 472}%
\special{pa 1960 460}%
\special{pa 2053 424}%
\special{pa 2084 411}%
\special{pa 2100 405}%
\special{fp}%
%
\special{pn 8}%
\special{pa 2100 405}%
\special{pa 600 405}%
\special{fp}%
%
\special{sh 1.000}%
\special{ia 2900 405 50 50 0.0000000 6.2831853}%
\special{pn 8}%
\special{ar 2900 405 50 50 0.0000000 6.2831853}%
%
\special{sh 1.000}%
\special{ia 3200 405 50 50 0.0000000 6.2831853}%
\special{pn 8}%
\special{ar 3200 405 50 50 0.0000000 6.2831853}%
%
\special{sh 1.000}%
\special{ia 3800 405 50 50 0.0000000 6.2831853}%
\special{pn 8}%
\special{ar 3800 405 50 50 0.0000000 6.2831853}%
%
\special{sh 1.000}%
\special{ia 3500 405 50 50 0.0000000 6.2831853}%
\special{pn 8}%
\special{ar 3500 405 50 50 0.0000000 6.2831853}%
%
\special{sh 1.000}%
\special{ia 4400 405 50 50 0.0000000 6.2831853}%
\special{pn 8}%
\special{ar 4400 405 50 50 0.0000000 6.2831853}%
%
\special{sh 1.000}%
\special{ia 4100 405 50 50 0.0000000 6.2831853}%
\special{pn 8}%
\special{ar 4100 405 50 50 0.0000000 6.2831853}%
\put(29.0000,-2.2500){\makebox(0,0){$\tilde{u}^{(s)}_{1}$}}%
\put(35.0000,-2.2500){\makebox(0,0){$\tilde{u}^{(s)}_{3}$}}%
\put(41.0000,-2.2500){\makebox(0,0){$\tilde{u}^{(s)}_{5}$}}%
\put(44.0000,-2.2500){\makebox(0,0){$\tilde{u}^{(s)}_{6}$}}%
\put(38.0000,-2.2500){\makebox(0,0){$\tilde{u}^{(s)}_{4}$}}%
\put(32.0000,-2.2500){\makebox(0,0){$\tilde{u}^{(s)}_{2}$}}%
%
\special{pn 8}%
\special{pa 2900 405}%
\special{pa 2931 417}%
\special{pa 2962 430}%
\special{pa 3086 478}%
\special{pa 3116 489}%
\special{pa 3209 522}%
\special{pa 3240 532}%
\special{pa 3333 559}%
\special{pa 3364 567}%
\special{pa 3426 581}%
\special{pa 3457 587}%
\special{pa 3487 592}%
\special{pa 3549 600}%
\special{pa 3611 604}%
\special{pa 3642 605}%
\special{pa 3673 605}%
\special{pa 3735 601}%
\special{pa 3766 598}%
\special{pa 3828 590}%
\special{pa 3859 584}%
\special{pa 3889 578}%
\special{pa 3920 571}%
\special{pa 3982 555}%
\special{pa 4044 537}%
\special{pa 4106 517}%
\special{pa 4199 484}%
\special{pa 4230 472}%
\special{pa 4260 460}%
\special{pa 4353 424}%
\special{pa 4384 411}%
\special{pa 4400 405}%
\special{fp}%
%
\special{pn 8}%
\special{pa 4400 405}%
\special{pa 2900 405}%
\special{fp}%
%
\special{pn 4}%
\special{sh 1}%
\special{ar 2500 700 16 16 0 6.2831853}%
\special{sh 1}%
\special{ar 2400 700 16 16 0 6.2831853}%
\special{sh 1}%
\special{ar 2600 700 16 16 0 6.2831853}%
\special{sh 1}%
\special{ar 2600 700 16 16 0 6.2831853}%
%
\special{sh 1.000}%
\special{ia 1350 1000 50 50 0.0000000 6.2831853}%
\special{pn 8}%
\special{ar 1350 1000 50 50 0.0000000 6.2831853}%
%
\special{sh 1.000}%
\special{ia 1000 1000 50 50 0.0000000 6.2831853}%
\special{pn 8}%
\special{ar 1000 1000 50 50 0.0000000 6.2831853}%
%
\special{sh 1.000}%
\special{ia 1700 1000 50 50 0.0000000 6.2831853}%
\special{pn 8}%
\special{ar 1700 1000 50 50 0.0000000 6.2831853}%
%
\special{pn 8}%
\special{pa 1000 1000}%
\special{pa 600 400}%
\special{fp}%
\special{pa 1500 400}%
\special{pa 1000 1000}%
\special{fp}%
%
\special{pn 8}%
\special{pa 1700 1000}%
\special{pa 2100 400}%
\special{fp}%
\special{pa 1200 400}%
\special{pa 1700 1000}%
\special{fp}%
%
\special{pn 8}%
\special{pa 1350 1000}%
\special{pa 900 400}%
\special{fp}%
\special{pa 1800 400}%
\special{pa 1350 1000}%
\special{fp}%
\put(10.3000,-8.1000){\makebox(0,0){$v^{(1)}_{1}$}}%
\put(17.3000,-8.1000){\makebox(0,0){$v^{(1)}_{3}$}}%
\put(13.8000,-8.1000){\makebox(0,0){$v^{(1)}_{2}$}}%
%
\special{sh 1.000}%
\special{ia 3650 1000 50 50 0.0000000 6.2831853}%
\special{pn 8}%
\special{ar 3650 1000 50 50 0.0000000 6.2831853}%
%
\special{sh 1.000}%
\special{ia 3300 1000 50 50 0.0000000 6.2831853}%
\special{pn 8}%
\special{ar 3300 1000 50 50 0.0000000 6.2831853}%
%
\special{sh 1.000}%
\special{ia 4000 1000 50 50 0.0000000 6.2831853}%
\special{pn 8}%
\special{ar 4000 1000 50 50 0.0000000 6.2831853}%
%
\special{pn 8}%
\special{pa 3300 1000}%
\special{pa 2900 400}%
\special{fp}%
\special{pa 3800 400}%
\special{pa 3300 1000}%
\special{fp}%
%
\special{pn 8}%
\special{pa 4000 1000}%
\special{pa 4400 400}%
\special{fp}%
\special{pa 3500 400}%
\special{pa 4000 1000}%
\special{fp}%
%
\special{pn 8}%
\special{pa 3650 1000}%
\special{pa 3200 400}%
\special{fp}%
\special{pa 4100 400}%
\special{pa 3650 1000}%
\special{fp}%
\put(33.3000,-8.1000){\makebox(0,0){$v^{(s)}_{1}$}}%
\put(40.3000,-8.1000){\makebox(0,0){$v^{(s)}_{3}$}}%
\put(36.8000,-8.1000){\makebox(0,0){$v^{(s)}_{2}$}}%
%
\special{sh 1.000}%
\special{ia 2500 1400 50 50 0.0000000 6.2831853}%
\special{pn 8}%
\special{ar 2500 1400 50 50 0.0000000 6.2831853}%
%
\special{pn 8}%
\special{pa 2500 1400}%
\special{pa 1000 1000}%
\special{fp}%
\special{pa 1350 1000}%
\special{pa 2500 1400}%
\special{fp}%
\special{pa 2500 1400}%
\special{pa 1700 1000}%
\special{fp}%
%
\special{pn 8}%
\special{pa 2500 1400}%
\special{pa 4000 1000}%
\special{fp}%
\special{pa 3650 1000}%
\special{pa 2500 1400}%
\special{fp}%
\special{pa 2500 1400}%
\special{pa 3300 1000}%
\special{fp}%
\put(25.0000,-12.8000){\makebox(0,0){$z$}}%
\end{picture}}%

\caption{Graph $H^{(4)}_{s}$}
\label{fH4}
\end{center}
\end{figure}

For an integer $s\geq 3$, we let $H^{(4)}_{s}$ be the graph obtained from $s$ pairwise vertex-disjoint copies $B_{1},\ldots ,B_{s}$ of $H^{(1)}_{1}$ such that $V(B_{i})=\{\tilde{u}^{(i)}_{j},v^{(i)}_{h}:1\leq j\leq 6,~1\leq h\leq 3\}$ where $\tilde{u}^{(i)}_{j}$ and $v^{(i)}_{h}$ respectively correspond to $u^{(1)}_{j}$ and $v_{h}$, by adding a new vertex $z$ and edges $zv^{(i)}_{h}~(1\leq i\leq s,~1\leq h\leq 3)$ (see Figure~\ref{fH4}).

\begin{lem}
\label{lem-nec-H^4-1}
For an integer $s\geq 3$, $H^{(4)}_{s}$ is $S_{8}^{(2)}$-free.
\end{lem}
\proof
Note that the subgraph of $H^{(4)}_{s}$ induced by $V(B_{i})\cup \{z\}$ is isomorphic to the Petersen graph.
Let $P$ be an induced path of $H^{(4)}_{s}$ having $8$ vertices.
Since the Petersen graph is $P_{6}$-free, it follows from the symmetry of vertices of $H^{(4)}_{s}$ that we may assume that $P=\tilde{u}^{(1)}_{3}\tilde{u}^{(1)}_{2}\tilde{u}^{(1)}_{1}v^{(1)}_{1}zv^{(2)}_{1}\tilde{u}^{(2)}_{1}\tilde{u}^{(2)}_{2}$.
Then $N_{H^{(4)}_{s}}(v^{(1)}_{1})-V(P)=\{\tilde{u}^{(1)}_{4}\}$.
Since the subgraph of $H^{(4)}_{s}$ induced by $V(P)\cup \{\tilde{u}^{(1)}_{4}\}$ contains a cycle, no induced subgraph of $H_{s}^{(4)}$ isomorphic to $S_{8}^{(2)}$ contains all vertices of $P$.
This implies that $H^{(4)}_{s}$ is $S_{8}^{(2)}$-free.
\qed

\medbreak\noindent\textit{Proof of Theorem~\ref{thm-nec-1}.}\quad
Let $T$ be a tree with ${\rm diam}(T)\geq 7$ which is not a caterpillar, and suppose that $\tilde{\GG}_{3}(\{C_{3},C_{4},T\})$ is a finite family.
By Lemma~\ref{lem-nec-1}, there exist infinitely many connected $3$-regular graphs with girth at least $5$.
Since $\tilde{\GG}_{3}(\{C_{3},C_{4},T\})$ is finite, there exists a connected $3$-regular graph containing an induced subgraph isomorphic to $T$.
This implies that $\Delta (T)\leq 3$.

For each $s\geq 5$, $H^{(1)}_{s}$ is a connected $\{C_{3},C_{4}\}$-free graph with $\delta (H^{(1)}_{s})\geq 3$.
Since $\tilde{\GG}_{3}(\{C_{3},C_{4},T\})$ is finite, there exists $s_{1}\geq 5$ such that $T\prec H^{(1)}_{s_{1}}$.
By Lemma~\ref{lem-nec-H^1-1}(i), we have ${\rm diam}(T)\leq 8$.
Furthermore, it follows from Lemma~\ref{lem-nec-H^1-1}(ii)(iii) that
\begin{align}
&\mbox{if ${\rm diam}(T)=8$, then $T\prec T^{*}_{9}$; and}\label{cond-thm-nec-1-1}\\
&\mbox{if ${\rm diam}(T)=7$, then $T\prec S_{8}(0,0,0,1)$ or $T\prec S_{8}^{(1)}$.}\label{cond-thm-nec-1-2}
\end{align}

For each $s\geq 3$, $H^{(2)}_{s}$ is a connected $\{C_{3},C_{4}\}$-free graph with $\delta (H^{(2)}_{s})\geq 3$.
Since $\tilde{\GG}_{3}(\{C_{3},C_{4},T\})$ is finite, there exists $s_{2}\geq 3$ such that $T\prec H^{(2)}_{s_{2}}$.
Then by Lemma~\ref{lem-nec-H^2-1},
\begin{align}
&S_{8}(0,0,0,1)\not\prec T; \mbox{ and}\label{cond-thm-nec-1-3}\\
&T^{*}_{8}\not\prec T.\label{cond-thm-nec-1-4}
\end{align}
Since $T$ is not a caterpillar, it follows from (\ref{cond-thm-nec-1-1}) and (\ref{cond-thm-nec-1-3}) that if ${\rm diam}(T)=8$, then $T$ is isomorphic to $T_{9}$, as desired.
Thus we may assume that ${\rm diam}(T)=7$.
Note that if a tree $S$ with ${\rm diam}(S)=7$ is not a caterpillar and satisfies $S_{8}(0,0,0,1)\not\prec S\prec S_{8}(0,0,0,1)$, then $S\simeq T_{8}~(\prec S_{8}^{(1)})$.
Hence by (\ref{cond-thm-nec-1-2}) and (\ref{cond-thm-nec-1-3}), we have $T\prec S_{8}^{(1)}$.
This together with (\ref{cond-thm-nec-1-4}) implies that
\begin{align}
T\prec S_{8}(1,1,1,0).\label{cond-thm-nec-1-5}
\end{align}

For each $s\geq 4$, $H^{(3)}_{s}$ is a connected $\{C_{3},C_{4}\}$-free graph with $\delta (H^{(3)}_{s})\geq 3$.
Since $\tilde{\GG}_{3}(\{C_{3},C_{4},T\})$ is finite, there exists $s_{3}\geq 4$ such that $T\prec H^{(3)}_{s_{3}}$.
Then by Lemma~\ref{lem-nec-H^3-1}, $S_{7}(1,0,1)\not\prec T$.
This together with (\ref{cond-thm-nec-1-5}) leads to $T\prec S_{8}(1,1,0,0)$ or $T\prec S_{8}(0,1,1,0)$.
Since $S_{8}(0,1,1,0)\simeq T_{8}^{(1)}$, the latter case leads to the desired conclusion.
Thus we may assume that
\begin{align}
S_{8}(1,0,0,0)\prec T\prec S_{8}(1,1,0,0).\label{cond-thm-nec-1-6}
\end{align}

For each $s\geq 3$, $H^{(4)}_{s}$ is a connected $\{C_{3},C_{4}\}$-free graph with $\delta (H^{(4)}_{s})\geq 3$.
Since $\tilde{\GG}_{3}(\{C_{3},C_{4},T\})$ is finite, there exists $s_{4}\geq 3$ such that $T\prec H^{(4)}_{s_{4}}$.
Then by Lemma~\ref{lem-nec-H^4-1}, $S_{8}^{(2)}\not\prec T$.
This together with (\ref{cond-thm-nec-1-6}) leads to $T\simeq S_{8}(1,0,0,0)~(\simeq T_{8}^{(2)})$, as desired.
\qed

\section{Proof of Theorem~\ref{prop-diambound}}\label{sec-diam}

In the remainder of this paper, we fix a connected $\{C_{3},C_{4}\}$-free graph $G$ with $\delta (G)\geq 3$.
We implicitly use the fact that $N_{G}(u)$ is an independent set of $G$ and $|N_{G}(u)\cap N_{G}(v)|\leq 1$ for any $u,v\in V(G)$ with $u\neq v$.
For $u\in V(G)$ and $U\subseteq V(G)$, we write $N(u)$ and $N(U)$ for $N_{G}(u)$ and $N_{G}(U)$, respectively.
Furthermore, in this section, we fix two vertices $u,v\in V(G)$ with ${\rm dist}_{G}(u,v)={\rm diam}(G)$ and a shortest $u$-$v$ path $P=u_{0}u_{1}\cdots u_{d}$ of $G$ where $d={\rm diam}(G)$, $u_{0}=u$ and $u_{d}=v$.
For each $i~(0\leq i\leq d)$, we fix a vertex $a_{i}\in N(u_{i})-V(P)$.

\begin{lem}
\label{lem-diambound-01}
\begin{enumerate}
\item[{\upshape(i)}]
For $y\in V(G)-V(P)$, $|N(y)\cap V(P)|\leq 1$.
In particular, $N(a_{i})\cap V(P)=\{u_{i}\}$ for $0\leq i\leq d$, and $a_{i}\neq a_{i'}$ for $0\leq i<i'\leq d$.
\item[{\upshape(ii)}]
For $0\leq i<i'\leq d$, if $a_{i}a_{i'}\in E(G)$, then $2\leq i'-i\leq 3$.
\end{enumerate}
\end{lem}
\proof
We first prove (i).
Let $y\in V(G)-V(P)$, and suppose that there exist integers $i$ and $i'$ with $0\leq i<i'\leq d$ such that $yu_{i},yu_{i'}\in E(G)$.
Since $G$ is $\{C_{3},C_{4}\}$-free, $i'-i\geq 3$.
Then $u_{0}u_{1}\cdots u_{i}yu_{i'}u_{i'+1}\cdots u_{d}$ is a $u$-$v$ path of $G$ having length $i+2+(d-i')~(\leq d-1)$, which contradicts the fact that ${\rm dist}_{G}(u,v)=d$.

Next we prove (ii).
Since $G$ is $C_{4}$-free, we have $i'-i\geq 2$.
If $i'-i\geq 4$, then $u_{0}u_{1}\cdots u_{i}a_{i}a_{i'}u_{i'}u_{i'+1}\cdots u_{d}$ is a $u$-$v$ path of $G$ having length $i+3+(d-i')~(\leq d-1)$, which contradicts the fact that ${\rm dist}_{G}(u,v)=d$.
\qed

\medbreak\noindent\textit{Proof of Theorem~\ref{prop-diambound}.}\quad
We first prove (iii).
Suppose that ${\rm diam}(G)\geq 12$.
We show that $T_{9}\prec G$.
Take a vertex $y\in N(a_{6})-\{u_{6}\}$.
By Lemma~\ref{lem-diambound-01}, we have $y\notin V(P)$, $|N(y)\cap V(P)|\leq 1$ and $N(y)\cap V(P)\subseteq \{u_{3},u_{4},u_{8},u_{9}\}$.
We may assume that $yu_{8},yu_{9}\notin E(G)$.
Then $N(y)\cap \{a_{6},u_{i}:5\leq i\leq 12\}=\{a_{6}\}$.
Since $|N(y)\cap N(u_{5})|\leq 1$, we can take a vertex $z\in \{u_{4},a_{5}\}-N(y)$.
Considering Lemma~\ref{lem-diambound-01}, we have $N(z)\cap \{y,a_{6},u_{i}:5\leq i\leq 12\}=\{u_{5}\}$.
Hence $\{z,u_{5},u_{6},y,a_{6},u_{7},\ldots ,u_{12}\}$ induces a copy of $T_{9}$ in $G$, which proves (iii).

Next we prove (i) and (ii).
Let $d_{1}=20$, $d_{2}=16$ and $p\in \{1,2\}$.
We suppose that ${\rm diam}(G)\geq d_{p}$, and show that $T_{8}^{(p)}\prec G$.

\begin{claim}
\label{cl-diambound-1}
Let $i$ be an integer with $5\leq i\leq d-5$, and suppose that there exists a vertex $y\in N(a_{i})-\{u_{i}\}$ such that $yu_{j}\in E(G)$ for some $0\leq j\leq d$ with $|i-j|\geq 3$.
Then $T_{8}^{(1)}\prec G$ and $T_{8}^{(2)}\prec G$.
\end{claim}
\proof
Without loss of generality, we may assume that $i-j\geq 3$ and $y=a_{j}$.
Then by Lemma~\ref{lem-diambound-01}(ii), $j=i-3$, i.e., $a_{i}a_{i-3}\in E(G)$.
If $E_{G}(\{a_{i-3},a_{i}\},\{a_{i+1},a_{i+2},a_{i+3}\})\neq \emptyset $, then it follows from Lemma~\ref{lem-diambound-01}(ii) that $a_{i}a_{h}\in E(G)$ for some $h\in \{i+2,i+3\}$, and hence $u_{0}u_{1}\cdots u_{i-3}a_{i-3}a_{i}a_{h}u_{h}u_{h+1}\cdots u_{d}$ is a $u$-$v$ path of $G$ having length $(i-3)+4+(d-h)~(=d+i-h+1\leq d-1)$, which contradicts the fact that ${\rm dist}_{G}(u,v)=d$.
Thus $E_{G}(\{a_{i-3},a_{i}\},\{a_{i+1},a_{i+2},a_{i+3}\})=\emptyset $.
Since $a_{i+2}a_{i+3}\notin E(G)$ by Lemma~\ref{lem-diambound-01}(ii), $\{u_{i-2},u_{i-1},u_{i},a_{i-3},a_{i},u_{i+1},u_{i+2},a_{i+2},u_{i+3},a_{i+3},u_{i+4},u_{i+5}\}$ induces a copy of $T_{8}^{(1)}$ in $G$.
Furthermore, $\{u_{i-2},u_{i-1},u_{i},a_{i-3},a_{i},u_{i+1},a_{i},u_{i+1},\ldots ,u_{i+5}\}$ induces a copy of $T_{8}^{(2)}$ in $G$.
\qed

\begin{claim}
\label{cl-diambound-2}
Let $i$ be an integer with $5\leq i\leq d-5$, and suppose that $a_{i}a_{i+h}\in E(G)$ for some $h\in \{2,-2\}$.
\begin{enumerate}
\item[{\upshape(a)}]
If $a_{i}a_{i-h}\notin E(G)$, then $T_{8}^{(1)}\prec G$.
\item[{\upshape(b)}]
If $a_{i-1}a_{i+1}\notin E(G)$, then $T_{8}^{(2)}\prec G$.
\end{enumerate}
\end{claim}
\proof
We may assume that $h=-2$, i.e., $a_{i}a_{i-2}\in E(G)$.
In view of Lemma~\ref{lem-diambound-01}(ii) and Claim~\ref{cl-diambound-1}, we may assume that $E(G[\{a_{j}:i-2\leq j\leq i+3\}])\subseteq \{a_{j}a_{j+2}:i-2\leq j\leq i+1\}$.
This implies that if $a_{i}a_{i+2}\notin E(G)$, then $\{a_{i-1},u_{i-1},u_{i},a_{i-2},a_{i},u_{i+1},u_{i+2},a_{i+2},u_{i+3},a_{i+3},u_{i+4},u_{i+5}\}$ induces a copy of $T_{8}^{(1)}$ in $G$; if $a_{i-1}a_{i+1}\notin E(G)$, then $\{a_{i-1},u_{i-1},u_{i},a_{i-2},a_{i},u_{i+1},a_{i+1},u_{i+2},\ldots ,u_{i+5}\}$ induces a copy of $T_{8}^{(2)}$ in $G$.
Thus we obtain the desired conclusions.
\qed

\begin{claim}
\label{cl-diambound-3}
Let $i$ be an integer with $6+p\leq j\leq d-(6+p)$, and let $y\in N(a_{j})-\{u_{j}\}$.
If $N(y)\cap V(P)\neq \emptyset $, then $T_{8}^{(p)}\prec G$.
\end{claim}
\proof
By way of contradiction, suppose that $yu_{j'}\in E(G)$ for some $0\leq j'\leq d$ with $j'\neq j$ and $G$ is $T_{8}^{(p)}$-free.
We may assume that $y=a_{j'}$.
Then by Lemma~\ref{lem-diambound-01}(ii) and Claim~\ref{cl-diambound-1}, we have $|j-j'|=2$.
We may assume that $j-j'=2$, i.e., $a_{j}a_{j-2}\in E(G)$.

Now we claim that
\begin{align}
a_{j}a_{j+2},a_{j+2}a_{j+4}\in E(G).\label{cond-cl-diambound-3-1}
\end{align}
If $p=1$, then $7\leq j<j+2\leq d-5$, and hence applying Claim~\ref{cl-diambound-2}(a) with $(i,h)=(j,-2)$ and then with $(i,h)=(j+2,-2)$, we get $a_{j}a_{j+2},a_{j+2}a_{j+4}\in E(G)$, as desired; if $p=2$, then $8\leq j<j+3\leq d-5$, and hence applying Claim~\ref{cl-diambound-2}(b) with $(i,h)=(j+\alpha ,-2)$ for $0\leq \alpha \leq 3$, we obtain $a_{j-1}a_{j+1},a_{j}a_{j+2},a_{j+1}a_{j+3},a_{j+2}a_{j+4}\in E(G)$, which implies (\ref{cond-cl-diambound-3-1}).

By the fact that $a_{j}a_{j-2}\in E(G)$ and (\ref{cond-cl-diambound-3-1}), $u_{0}u_{1}\cdots u_{j-2}a_{j-2}a_{j}a_{j+2}a_{j+4}u_{j+4}u_{j+5}\cdots u_{d}$ is a $u$-$v$ path of $G$ having length $(j-2)+5+(d-(j+4))~(=d-1)$, which contradicts the fact that ${\rm dist}_{G}(u,v)=d$.
\qed

\begin{claim}
\label{cl-diambound-4}
There exists an integer $i_{0}\in \{7,10,13\}$ such that $N(a_{i_{0}})\cap N(a_{i_{0}+3})=\emptyset $.
\end{claim}
\proof
If $N(a_{i})\cap N(a_{i+3})\neq \emptyset $ for all $i\in \{7,10,13\}$, then there exists an $a_{7}$-$a_{16}$ path $Q$ of $G-V(P)$ having length at most $6$, and hence $u_{0}u_{1}\cdots u_{7}a_{7}Qa_{16}u_{16}u_{17}\cdots u_{d}$ is a $u$-$v$ path of length at most $7+8+(d-16)~(=d-1)$, which contradicts the fact that ${\rm dist}_{G}(u,v)=d$.
\qed

Now we complete the proof of the theorem.
We first consider the case where $p=1$.
Let $i_{0}\in \{7,10,13\}$ be as in Claim~\ref{cl-diambound-4}.
Since $|N(a_{i_{0}})\cap N(a_{i_{0}+2})|\leq 1$, there exists a vertex $y\in N(a_{i_{0}})-(\{u_{i_{0}}\}\cup N(a_{i_{0}+2}))$.
By Claim~\ref{cl-diambound-3}, we may assume that $N(y)\cap V(P)=\emptyset $ and $\{a_{i_{0}},a_{i_{0}+2},a_{i_{0}+3}\}$ is an independent set of $G$.
In particular, $y\notin \{a_{i_{0}+2},a_{i_{0}+3}\}$.
Furthermore, it follows from the definition of $i_{0}$ and $y$ that $\{y,a_{i_{0}+2},a_{i_{0}+3}\}$ is also an independent set of $G$.
Consequently, $\{u_{i_{0}-2},u_{i_{0}-1},u_{i_{0}},y,a_{i_{0}},u_{i_{0}+1},u_{i_{0}+2},a_{i_{0}+2},u_{i_{0}+3},a_{i_{0}+3},u_{i_{0}+4},u_{i_{0}+5}\}$ induces a copy of $T_{8}^{(1)}$ in $G$.

Finally we consider the case where $p=2$.
Recall that $a_{8}\neq a_{9}$ and $a_{8}a_{9}\notin E(G)$.
Since $|N(a_{8})\cap N(a_{9})|\leq 1$, there exists a vertex $y'\in N(a_{8})-(\{u_{8}\}\cup N(a_{9}))$.
By Claim~\ref{cl-diambound-3}, we may assume that $N(y')\cap V(P)=\emptyset $.
Consequently, $\{u_{6},u_{7},u_{8},y',a_{8},u_{9},a_{9},u_{10},\ldots ,u_{13}\}$ induces a copy of $T_{8}^{(2)}$ in $G$.

This completes the proof of Theorem~\ref{prop-diambound}.
\qed

\section{Paths between two vertices and the existence of induced trees}\label{sec-M-cond}

For an integer $k\geq 4$  and two non-adjacent vertices $v$ and $w$ of $G$, an induced $v$-$w$ path $P$ of $G$ with $|V(P)|=k$ is called a {\it $(v,w;k)$-path}.
Let
$$
M_{k}^{w}(v)=\{x\in V(G):\mbox{there exists a $(v,w;k)$-path $P$ such that $N(v)\cap V(P)=\{x\}$}\}.
$$

\begin{lem}
\label{lem-Mk-fund1}
Let $k\geq 4$ be an integer, and let $v$ and $w$ be non-adjacent vertices of $G$.
For $i\in \{1,2\}$, let $Q_{i}=u^{(i)}_{1}u^{(i)}_{2}\cdots u^{(i)}_{k}$ be a $(v,w;k)$-path where $u^{(i)}_{1}=v$ and $u^{(i)}_{k}=w$, and suppose that $u^{(1)}_{2}\neq u^{(2)}_{2}$.
Then the following hold.
\begin{enumerate}
\item[{\upshape(i)}]
We have $\{u^{(1)}_{2},u^{(1)}_{3}\}\cap \{u^{(2)}_{2},u^{(2)}_{3}\}=\emptyset $.
\item[{\upshape(ii)}]
If $k=5$, $\{u^{(1)}_{2},u^{(1)}_{3},u^{(1)}_{4}\}\cap \{u^{(2)}_{2},u^{(2)}_{3}\}=\emptyset $.
\item[{\upshape(iii)}]
If $k=4$, then $E_{G}(\{u^{(1)}_{2},u^{(1)}_{3}\},\{u^{(2)}_{2},u^{(2)}_{3}\})=\emptyset $.
\item[{\upshape(iv)}]
If $k=5$ and $u^{(1)}_{4}\neq u^{(2)}_{4}$, then $E_{G}(\{u^{(1)}_{2},u^{(1)}_{3},u^{(1)}_{4}\},\{u^{(2)}_{2},u^{(2)}_{3},u^{(2)}_{4}\})\subseteq \{u^{(1)}_{2}u^{(2)}_{4},u^{(1)}_{3}u^{(2)}_{3},u^{(1)}_{4}u^{(2)}_{2}\}$.
\item[{\upshape(v)}]
If $k=5$, $u^{(1)}_{4}\neq u^{(2)}_{4}$ and $u^{(1)}_{2}\notin M_{4}^{w}(v)$, then $E_{G}(\{u^{(1)}_{2},u^{(1)}_{3},u^{(1)}_{4}\},\{u^{(2)}_{2},u^{(2)}_{3},u^{(1)}_{4}\})\subseteq \{u^{(1)}_{3}u^{(2)}_{3},u^{(1)}_{4}u^{(2)}_{2}\}$.
\end{enumerate}
\end{lem}
\proof
Since $G$ is $\{C_{3},C_{4}\}$-free, we see that (i) and (iii) hold.
If $k=5$ and $u^{(1)}_{4}\in \{u^{(2)}_{2},u^{(2)}_{3}\}$, then $N(w)\cap \{u^{(2)}_{2},u^{(2)}_{3}\}\neq \emptyset $, and hence $Q_{2}$ is not an induced path, which is a contradiction.
Thus we obtain (ii).
By the $\{C_{3},C_{4}\}$-freeness of $G$ and (ii), we get (iv).
If $k=5$ and $u^{(1)}_{2}u^{(2)}_{4}\in E(G)$, then $vu^{(1)}_{2}u^{(2)}_{4}w$ is a $(v,w;4)$-path.
This together with (iv) leads to (v).
\qed

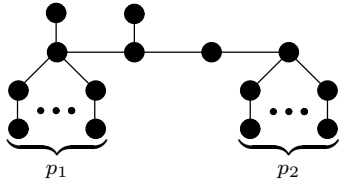
\begin{figure}
\begin{center}
{\unitlength 0.1in%
\begin{picture}(24.7500,7.5000)(1.7000,-10.9000)%
%
\special{sh 1.000}%
\special{ia 2195 995 50 50 0.0000000 6.2831853}%
\special{pn 8}%
\special{ar 2195 995 50 50 0.0000000 6.2831853}%
%
\special{sh 1.000}%
\special{ia 2195 795 50 50 0.0000000 6.2831853}%
\special{pn 8}%
\special{ar 2195 795 50 50 0.0000000 6.2831853}%
%
\special{sh 1.000}%
\special{ia 2395 595 50 50 0.0000000 6.2831853}%
\special{pn 8}%
\special{ar 2395 595 50 50 0.0000000 6.2831853}%
%
\special{sh 1.000}%
\special{ia 2595 795 50 50 0.0000000 6.2831853}%
\special{pn 8}%
\special{ar 2595 795 50 50 0.0000000 6.2831853}%
%
\special{sh 1.000}%
\special{ia 2595 995 50 50 0.0000000 6.2831853}%
\special{pn 8}%
\special{ar 2595 995 50 50 0.0000000 6.2831853}%
%
\special{pn 8}%
\special{pa 2595 995}%
\special{pa 2595 795}%
\special{fp}%
\special{pa 2595 795}%
\special{pa 2395 595}%
\special{fp}%
\special{pa 2395 595}%
\special{pa 2195 795}%
\special{fp}%
\special{pa 2195 795}%
\special{pa 2195 995}%
\special{fp}%
%
\special{sh 1.000}%
\special{ia 1995 595 50 50 0.0000000 6.2831853}%
\special{pn 8}%
\special{ar 1995 595 50 50 0.0000000 6.2831853}%
%
\special{sh 1.000}%
\special{ia 1595 595 50 50 0.0000000 6.2831853}%
\special{pn 8}%
\special{ar 1595 595 50 50 0.0000000 6.2831853}%
%
\special{sh 1.000}%
\special{ia 1195 595 50 50 0.0000000 6.2831853}%
\special{pn 8}%
\special{ar 1195 595 50 50 0.0000000 6.2831853}%
%
\special{sh 1.000}%
\special{ia 1395 795 50 50 0.0000000 6.2831853}%
\special{pn 8}%
\special{ar 1395 795 50 50 0.0000000 6.2831853}%
%
\special{sh 1.000}%
\special{ia 1395 995 50 50 0.0000000 6.2831853}%
\special{pn 8}%
\special{ar 1395 995 50 50 0.0000000 6.2831853}%
%
\special{sh 1.000}%
\special{ia 995 995 50 50 0.0000000 6.2831853}%
\special{pn 8}%
\special{ar 995 995 50 50 0.0000000 6.2831853}%
%
\special{sh 1.000}%
\special{ia 995 795 50 50 0.0000000 6.2831853}%
\special{pn 8}%
\special{ar 995 795 50 50 0.0000000 6.2831853}%
%
\special{sh 1.000}%
\special{ia 1190 390 50 50 0.0000000 6.2831853}%
\special{pn 8}%
\special{ar 1190 390 50 50 0.0000000 6.2831853}%
%
\special{sh 1.000}%
\special{ia 1595 395 50 50 0.0000000 6.2831853}%
\special{pn 8}%
\special{ar 1595 395 50 50 0.0000000 6.2831853}%
%
\special{pn 8}%
\special{pa 1195 595}%
\special{pa 2395 595}%
\special{fp}%
%
\special{pn 8}%
\special{pa 1595 595}%
\special{pa 1595 395}%
\special{fp}%
\put(11.9500,-11.5500){\makebox(0,0){$\underbrace{\hspace{13mm}}_{p_{1}}$}}%
%
\special{pn 4}%
\special{sh 1}%
\special{ar 1190 900 16 16 0 6.2831853}%
\special{sh 1}%
\special{ar 1110 900 16 16 0 6.2831853}%
\special{sh 1}%
\special{ar 1270 900 16 16 0 6.2831853}%
\special{sh 1}%
\special{ar 1270 900 16 16 0 6.2831853}%
%
\special{pn 4}%
\special{sh 1}%
\special{ar 2395 905 16 16 0 6.2831853}%
\special{sh 1}%
\special{ar 2315 905 16 16 0 6.2831853}%
\special{sh 1}%
\special{ar 2475 905 16 16 0 6.2831853}%
\special{sh 1}%
\special{ar 2475 905 16 16 0 6.2831853}%
\put(23.9500,-11.5500){\makebox(0,0){$\underbrace{\hspace{13mm}}_{p_{2}}$}}%
%
\special{pn 8}%
\special{pa 1190 390}%
\special{pa 1190 590}%
\special{fp}%
\special{pa 1190 590}%
\special{pa 990 790}%
\special{fp}%
\special{pa 990 790}%
\special{pa 990 990}%
\special{fp}%
\special{pa 1390 990}%
\special{pa 1390 790}%
\special{fp}%
%
\special{pn 8}%
\special{pa 1190 590}%
\special{pa 1390 790}%
\special{fp}%
\end{picture}}%

\caption{Graph $T_{8}^{*}(p_{1},p_{2})$}
\label{f2}
\end{center}
\end{figure}

For integers $p_{1},p_{2}\geq 1$, let $T_{8}^{*}(p_{1},p_{2})$ be the tree of order $2p_{1}+2p_{2}+6$ depicted in Figure~\ref{f2}.
Note that $T_{8}^{*}(1,2)\simeq T_{8}^{(1)}$ and $T_{8}^{(2)}\prec T_{8}^{*}(2,1)$.

\begin{lem}
\label{lem-Mk-1}
Let $p_{1},p_{2}\geq 1$ be integers, and let $n=2R(3,p_{1}+2)+3R(3,p_{2})+2p_{2}(p_{1}+p_{2}+1)+3$.
Let $v$ and $w$ be non-adjacent vertices of $G$.
If $|M_{4}^{w}(v)|\geq n$, then $T_{8}^{*}(p_{1},p_{2})\prec G$.
\end{lem}
\proof
Let $n_{1}=R(3,p_{1}+2)+p_{2}(p_{1}+p_{2}+1)$ and $n_{2}=R(3,p_{2})+1$.
Note that $n=2n_{1}+3n_{2}$.
Let $a_{1},\ldots ,a_{n}$ be distinct vertices in $M_{4}^{w}(v)$.
Then for each $i~(1\leq i\leq n)$, there exists a vertex $b_{i}\in V(G)$ such that $Q_{i}=va_{i}b_{i}w$ is a $(v,w;4)$-path.
By Lemma~\ref{lem-Mk-fund1}(i)(iii), $\{a_{i},b_{i}\}\cap \{a_{j},b_{j}\}=\emptyset $ and $E_{G}(\{a_{i},b_{i}\},\{a_{j},b_{j}\})=\emptyset $ for $1\leq i<j\leq n$.
For $i~(1\leq i\leq n_{1}+n_{2})$, take $x_{i}\in N(a_{i})-\{v,b_{i}\}$.
Then by Lemma~\ref{lem-Mk-fund1}(iii), $\{x_{i}:1\leq i\leq n_{1}+n_{2}\}\cap \{v,w,a_{i},b_{i}:1\leq i\leq n\}=\emptyset $.
Since $G$ is $\{C_{3},C_{4}\}$-free, $N(x_{i})\cap \{v,w,a_{j}:1\leq j\leq n,~j\neq i\}=\emptyset $ for each $i~(1\leq i\leq n_{1}+n_{2})$ and, in particular, $x_{1},\ldots ,x_{n_{1}+n_{2}}$ are distinct. 
Since $|N(x_{i})\cap \{b_{j}:1\leq j\leq n\}|\leq |N(x_{i})\cap N(w)|\leq 1$, we have $|N(\{x_{i}:1\leq i\leq n_{1}+n_{2}\})\cap \{b_{i}:n_{1}+n_{2}+1\leq i\leq n\}|\leq n_{1}+n_{2}$.
Since $n-(n_{1}+n_{2})=n_{1}+2n_{2}$, we may assume that $N(\{x_{i}:1\leq i\leq n_{1}+n_{2}\})\cap \{b_{i}:n_{1}+n_{2}+1\leq i\leq n_{1}+2n_{2}\}=\emptyset $.
For each $i~(n_{1}+n_{2}+1\leq i\leq n_{1}+2n_{2})$, take $y_{i}\in N(b_{i})-\{w,a_{i}\}$.
Set $Y=\{y_{i}:n_{1}+n_{2}+1\leq i\leq n_{1}+2n_{2}\}$.
Then by Lemma~\ref{lem-Mk-fund1}(iii), $Y\cap \{v,w,a_{i},b_{i}:1\leq i\leq n_{1}+n_{2}\}=\emptyset $.
Since $G$ is $\{C_{3},C_{4}\}$-free, $N(y_{i})\cap \{v,w,b_{j}:1\leq j\leq n_{1}+2n_{2},~j\neq i\}=\emptyset $ for each $i~(n_{1}+n_{2}+1\leq i\leq n_{1}+2n_{2})$ and, in particular, $y_{n_{1}+n_{2}+1},\ldots ,y_{n_{1}+2n_{2}}$ are distinct.
Since $N(\{x_{i}:1\leq i\leq n_{1}+n_{2}\})\cap \{b_{i}:n_{1}+n_{2}+1\leq i\leq n_{1}+2n_{2}\}=\emptyset $, we also have $Y\cap \{x_{i}:1\leq i\leq n_{1}+n_{2}\}=\emptyset $.
Since $|Y|=n_{2}$ and $|N(y_{i})\cap \{a_{j}:1\leq j\leq n_{1}+n_{2}\}|\leq |N(y_{i})\cap N(v)|\leq 1$ for each $i~(n_{1}+n_{2}+1\leq i\leq n_{1}+2n_{2})$, we may assume that $N(Y)\cap \{a_{i}:1\leq i\leq n_{1}\}=\emptyset $.
Set $X=\{x_{i}:1\leq i\leq n_{1}\}$.
We have the following:
\begin{align}
&x_{1},\ldots ,x_{n_{1}},y_{n_{1}+n_{2}+1},\ldots ,y_{n_{1}+2n_{2}}\mbox{ are distinct};\label{cond-lem-Mk-1-0}\\
&(X\cup Y)\cap \left(\bigcup _{1\leq i\leq n_{1}+2n_{2}}V(Q_{i})\right)=\emptyset ;\label{cond-lem-Mk-1-1}\\
&N(x_{i})\cap \left(\bigcup _{n_{1}+n_{2}+1\leq j\leq n_{1}+2n_{2}}V(Q_{j})\right)=\emptyset \mbox{ for }1\leq i\leq n_{1};\label{cond-lem-Mk-1-2}\\
&N(x_{i})\cap \{a_{j}:1\leq j\leq n_{1},~j\neq i\}=\emptyset \mbox{ for }1\leq i\leq n_{1};\label{cond-lem-Mk-1-3}\\
&N(y_{i})\cap \left(\bigcup _{1\leq j\leq n_{1}}V(Q_{j})\right)=\emptyset \mbox{ for }n_{1}+n_{2}+1\leq i\leq n_{1}+2n_{2}; \mbox{ and}\label{cond-lem-Mk-1-4}\\
&N(y_{i})\cap \{b_{j}:n_{1}+n_{2}+1\leq j\leq n_{1}+2n_{2},~j\neq i\}=\emptyset \mbox{ for }n_{1}+n_{2}+1\leq i\leq n_{1}+2n_{2}.\label{cond-lem-Mk-1-5}
\end{align}

\begin{claim}
\label{cl-lem-Mk-1}
There exists $X_{0}\subseteq X$ and there exists $Y_{0}\subseteq Y$ with $|X_{0}|=p_{1}+2$ and $|Y_{0}|=p_{2}$ such that $X_{0}\cup Y_{0}$ is an independent set of $G$.
\end{claim}
\proof
We first assume that there exists $u_{0}\in Y$ such that $|N(u_{0})\cap X|\geq p_{1}+p_{2}+2$.
Note that $N(u_{0})\cap X$ is an independent set of $G$.
Since $G$ is $C_{3}$-free and $|Y-\{u_{0}\}|=n_{2}-1=R(3,p_{2})$, there exists an independent set $Y_{0}\subseteq Y-\{u_{0}\}$ of $G$ with $|Y_{0}|=p_{2}$.
We have $|N(u_{0})\cap X\cap N(y)|\leq |N(u_{0})\cap N(y)|\leq 1$ for each $y\in Y_{0}$.
In particular, $|(N(u_{0})\cap X)-N(Y_{0})|\geq (p_{1}+p_{2}+2)-p_{2}=p_{1}+2$.
Take $X_{0}\subseteq (N(u_{0})\cap X)-N(Y_{0})$ with $|X_{0}|=p_{1}+2$.
Then $X_{0}\cup Y_{0}$ is an independent set of $G$, as desired.
Thus we may assume that
\begin{align}
\mbox{$|N(y)\cap X|\leq p_{1}+p_{2}+1$ for all $y\in Y$.}\label{cond-cl-lem-Mk-1-1}
\end{align}
Since $G$ is $C_{3}$-free and $|Y|>R(3,p_{2})$, there exists an independent set $Y_{0}\subseteq Y$ of $G$ with $|Y_{0}|=p_{2}$.
By (\ref{cond-cl-lem-Mk-1-1}),
\begin{align*}
|X-N(Y_{0})| &\geq n_{1}-\sum _{y\in Y_{0}}|N(y)\cap X|\\
&\geq (R(3,p_{1}+2)+p_{2}(p_{1}+p_{2}+1))-p_{2}(p_{1}+p_{2}+1)\\
&= R(3,p_{1}+2).
\end{align*}
Since $G$ is $C_{3}$-free, there exists an independent set $X_{0}\subseteq X-N(Y_{0})$ of $G$ with $|X_{0}|=p_{1}+2$.
Then $X_{0}\cup Y_{0}$ is an independent set of $G$.
\qed

\begin{figure}
\begin{center}
{\unitlength 0.1in%
\begin{picture}(44.0000,18.7000)(4.5000,-22.7000)%
\put(28.0000,-23.3500){\makebox(0,0){$w$}}%
%
\special{sh 1.000}%
\special{ia 2790 595 50 50 0.0000000 6.2831853}%
\special{pn 8}%
\special{ar 2790 595 50 50 0.0000000 6.2831853}%
\put(27.9000,-4.6500){\makebox(0,0){$v$}}%
%
\special{sh 1.000}%
\special{ia 1000 995 50 50 0.0000000 6.2831853}%
\special{pn 8}%
\special{ar 1000 995 50 50 0.0000000 6.2831853}%
%
\special{sh 1.000}%
\special{ia 1600 995 50 50 0.0000000 6.2831853}%
\special{pn 8}%
\special{ar 1600 995 50 50 0.0000000 6.2831853}%
%
\special{sh 1.000}%
\special{ia 2600 995 50 50 0.0000000 6.2831853}%
\special{pn 8}%
\special{ar 2600 995 50 50 0.0000000 6.2831853}%
%
\special{sh 1.000}%
\special{ia 3000 995 50 50 0.0000000 6.2831853}%
\special{pn 8}%
\special{ar 3000 995 50 50 0.0000000 6.2831853}%
%
\special{sh 1.000}%
\special{ia 3600 995 50 50 0.0000000 6.2831853}%
\special{pn 8}%
\special{ar 3600 995 50 50 0.0000000 6.2831853}%
%
\special{sh 1.000}%
\special{ia 4600 995 50 50 0.0000000 6.2831853}%
\special{pn 8}%
\special{ar 4600 995 50 50 0.0000000 6.2831853}%
%
\special{sh 1.000}%
\special{ia 1000 1800 50 50 0.0000000 6.2831853}%
\special{pn 8}%
\special{ar 1000 1800 50 50 0.0000000 6.2831853}%
%
\special{sh 1.000}%
\special{ia 1600 1800 50 50 0.0000000 6.2831853}%
\special{pn 8}%
\special{ar 1600 1800 50 50 0.0000000 6.2831853}%
%
\special{sh 1.000}%
\special{ia 2600 1800 50 50 0.0000000 6.2831853}%
\special{pn 8}%
\special{ar 2600 1800 50 50 0.0000000 6.2831853}%
%
\special{sh 1.000}%
\special{ia 3000 1800 50 50 0.0000000 6.2831853}%
\special{pn 8}%
\special{ar 3000 1800 50 50 0.0000000 6.2831853}%
%
\special{sh 1.000}%
\special{ia 3600 1800 50 50 0.0000000 6.2831853}%
\special{pn 8}%
\special{ar 3600 1800 50 50 0.0000000 6.2831853}%
%
\special{sh 1.000}%
\special{ia 4600 1800 50 50 0.0000000 6.2831853}%
\special{pn 8}%
\special{ar 4600 1800 50 50 0.0000000 6.2831853}%
%
\special{sh 1.000}%
\special{ia 2800 2205 50 50 0.0000000 6.2831853}%
\special{pn 8}%
\special{ar 2800 2205 50 50 0.0000000 6.2831853}%
%
\special{sh 1.000}%
\special{ia 2400 1405 50 50 0.0000000 6.2831853}%
\special{pn 8}%
\special{ar 2400 1405 50 50 0.0000000 6.2831853}%
%
\special{sh 1.000}%
\special{ia 1400 1405 50 50 0.0000000 6.2831853}%
\special{pn 8}%
\special{ar 1400 1405 50 50 0.0000000 6.2831853}%
%
\special{sh 1.000}%
\special{ia 800 1405 50 50 0.0000000 6.2831853}%
\special{pn 8}%
\special{ar 800 1405 50 50 0.0000000 6.2831853}%
%
\special{sh 1.000}%
\special{ia 3800 1400 50 50 0.0000000 6.2831853}%
\special{pn 8}%
\special{ar 3800 1400 50 50 0.0000000 6.2831853}%
%
\special{sh 1.000}%
\special{ia 4800 1400 50 50 0.0000000 6.2831853}%
\special{pn 8}%
\special{ar 4800 1400 50 50 0.0000000 6.2831853}%
\put(8.5000,-10.0000){\makebox(0,0){$a_{1}$}}%
\put(14.5000,-10.0000){\makebox(0,0){$a_{2}$}}%
\put(23.5000,-10.0000){\makebox(0,0){$a_{p_{1}+1}$}}%
\put(32.5000,-10.0000){\makebox(0,0){$a_{p_{1}+2}$}}%
\put(39.7000,-10.0000){\makebox(0,0){$a_{n_{1}+n_{2}+1}$}}%
\put(50.0000,-10.0000){\makebox(0,0){$a_{n_{1}+n_{2}+p_{2}}$}}%
\put(6.5000,-14.0500){\makebox(0,0){$x_{1}$}}%
\put(21.5000,-14.0500){\makebox(0,0){$x_{p_{1}+1}$}}%
\put(41.7000,-14.0500){\makebox(0,0){$y_{n_{1}+n_{2}+1}$}}%
\put(52.0000,-14.0500){\makebox(0,0){$y_{n_{1}+n_{2}+p_{2}}$}}%
\put(8.5000,-18.0500){\makebox(0,0){$b_{1}$}}%
\put(14.5000,-18.0500){\makebox(0,0){$b_{2}$}}%
\put(23.5000,-18.0500){\makebox(0,0){$b_{p_{1}+1}$}}%
\put(32.5000,-18.0500){\makebox(0,0){$b_{p_{1}+2}$}}%
\put(39.7000,-18.0500){\makebox(0,0){$b_{n_{1}+n_{2}+1}$}}%
\put(50.0000,-18.0500){\makebox(0,0){$b_{n_{1}+n_{2}+p_{2}}$}}%
\put(12.5000,-14.0500){\makebox(0,0){$x_{2}$}}%
%
\special{pn 4}%
\special{sh 1}%
\special{ar 2000 1200 16 16 0 6.2831853}%
\special{sh 1}%
\special{ar 2200 1200 16 16 0 6.2831853}%
\special{sh 1}%
\special{ar 2100 1200 16 16 0 6.2831853}%
\special{sh 1}%
\special{ar 2100 1200 16 16 0 6.2831853}%
%
\special{pn 4}%
\special{sh 1}%
\special{ar 4200 1600 16 16 0 6.2831853}%
\special{sh 1}%
\special{ar 4400 1600 16 16 0 6.2831853}%
\special{sh 1}%
\special{ar 4300 1600 16 16 0 6.2831853}%
\special{sh 1}%
\special{ar 4300 1600 16 16 0 6.2831853}%
%
\special{pn 20}%
\special{pa 2800 600}%
\special{pa 1000 1000}%
\special{fp}%
\special{pa 1000 1000}%
\special{pa 800 1400}%
\special{fp}%
\special{pa 1000 1000}%
\special{pa 1000 1800}%
\special{fp}%
\special{pa 1000 1800}%
\special{pa 2800 2200}%
\special{fp}%
\special{pa 2800 2200}%
\special{pa 3600 1800}%
\special{fp}%
\special{pa 3600 1800}%
\special{pa 3800 1400}%
\special{fp}%
\special{pa 4800 1400}%
\special{pa 4600 1800}%
\special{fp}%
\special{pa 4600 1800}%
\special{pa 2800 2200}%
\special{fp}%
\special{pa 2800 600}%
\special{pa 1600 1000}%
\special{fp}%
\special{pa 1600 1000}%
\special{pa 1400 1400}%
\special{fp}%
\special{pa 2400 1400}%
\special{pa 2600 1000}%
\special{fp}%
\special{pa 2600 1000}%
\special{pa 2800 600}%
\special{fp}%
\special{pa 2800 600}%
\special{pa 3000 1000}%
\special{fp}%
%
\special{pn 4}%
\special{pa 2800 600}%
\special{pa 3600 1000}%
\special{fp}%
\special{pa 3600 1000}%
\special{pa 3600 1800}%
\special{fp}%
\special{pa 4600 1800}%
\special{pa 4600 1000}%
\special{fp}%
\special{pa 4600 1000}%
\special{pa 2800 600}%
\special{fp}%
\special{pa 3000 1000}%
\special{pa 3000 1800}%
\special{fp}%
\special{pa 3000 1800}%
\special{pa 2800 2200}%
\special{fp}%
\special{pa 2800 2200}%
\special{pa 2600 1800}%
\special{fp}%
\special{pa 2600 1800}%
\special{pa 2600 1000}%
\special{fp}%
\special{pa 1600 1000}%
\special{pa 1600 1800}%
\special{fp}%
\special{pa 1600 1800}%
\special{pa 2800 2200}%
\special{fp}%
\end{picture}}%

\caption{An induced copy of $T_{8}^{*}(p_{1},p_{2})$ in Lemma~\ref{lem-Mk-1}}
\label{flem4.2}
\end{center}
\end{figure}
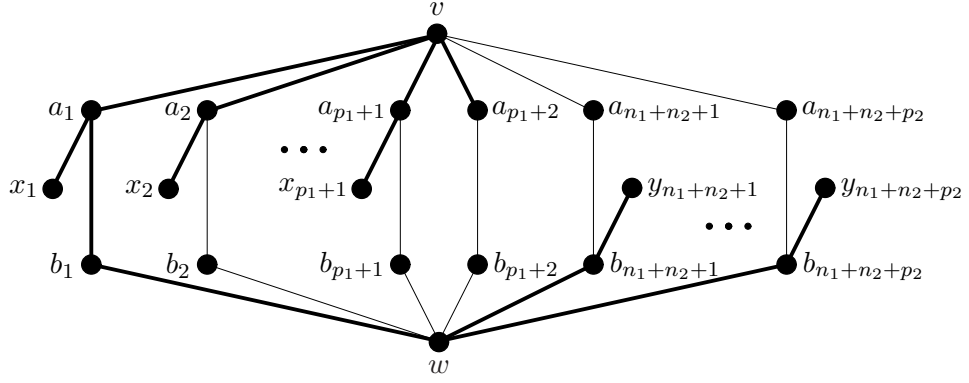

By Claim~\ref{cl-lem-Mk-1}, we may assume that
\begin{align}
\{x_{1},\ldots ,x_{p_{1}+2},y_{n_{1}+n_{2}+1},\ldots ,y_{n_{1}+n_{2}+p_{2}}\}\mbox{ is an independent set of }G.\label{cond-lem-Mk-1-5+}
\end{align}
For each $i~(1\leq i\leq p_{1}+2)$, let $I_{i}=\{j:1\leq j\leq p_{1}+2,~x_{j}b_{i}\in E(G)\}$.
For $1\leq j\leq p_{1}+2$, since $|N(x_{j})\cap N(w)|\leq 1$, $j$ belongs to at most one of $I_{1},\ldots ,I_{p_{1}+2}$.
In particular, $I_{1},\ldots ,I_{p_{1}+2}$ are pairwise disjoint.
Consequently, $\sum _{1\leq i\leq p_{1}+2}|I_{i}|=|\bigcup _{1\leq i\leq p_{1}+2}I_{i}|\leq p_{1}+2$.
This implies that there exists $i_{0}$ with $1\leq i_{0}\leq p_{1}+2$ such that $|I_{i_{0}}|\leq 1$.
We may assume that $i_{0}=1$ and $2,\ldots ,p_{1}+1\notin I_{1}$, that is,
\begin{align}
N(b_{1})\cap \{x_{j}:2\leq j\leq p_{1}+1\}=\emptyset .\label{cond-lem-Mk-1-6}
\end{align}
Since $b_{1}x_{1}\notin E(G)$ by the $C_{3}$-freeness of $G$, it follows from (\ref{cond-lem-Mk-1-0})--(\ref{cond-lem-Mk-1-5}), (\ref{cond-lem-Mk-1-5+}) and (\ref{cond-lem-Mk-1-6}) that
$$
\{a_{i},x_{i}:1\leq i\leq p_{1}+1\}\cup \{v,a_{p_{1}+2},b_{1},w\}\cup \{b_{i},y_{i}:n_{1}+n_{2}+1\leq i\leq n_{1}+n_{2}+p_{2}\}
$$
induces a copy of $T_{8}^{*}(p_{1},p_{2})$ in $G$ (see Figure~\ref{flem4.2}).
\qed

Recall that $T_{8}^{*}(1,2)\simeq T_{8}^{(1)}$ and $T_{8}^{(2)}\prec T_{8}^{*}(2,1)$.
Since $R(3,1)=1$, $R(3,2)=3$, $R(3,3)=6$ and $R(3,4)=9$ (see \cite{R}), we obtain the following as a corollary of Lemma~\ref{lem-Mk-1}.

\begin{lem}
\label{lem-Mk-1-cor}
Let $v$ and $w$ be non-adjacent vertices of $G$.
\begin{enumerate}
\item[{\upshape(i)}]
If $|M_{4}^{w}(v)|\geq 40$, then $T_{8}^{(1)}\prec G$.
\item[{\upshape(ii)}]
If $|M_{4}^{w}(v)|\geq 32$, then $T_{8}^{(2)}\prec G$.
\end{enumerate}
\end{lem}

\begin{lem}
\label{lem-Mk-2}
Let $v$ and $w$ be non-adjacent vertices of $G$.
\begin{enumerate}
\item[{\upshape(i)}]
If $|M_{4}^{w}(v)\cup M_{5}^{w}(v)|\geq 4837$, then $T_{8}^{(1)}\prec G$.
\item[{\upshape(ii)}]
If $|M_{4}^{w}(v)\cup M_{5}^{w}(v)|\geq 2047$, then $T_{8}^{(2)}\prec G$.
\end{enumerate}
\end{lem}
\proof
Let $h_{1}=40$, $m_{1}=124$, $h_{2}=32$ and $m_{2}=66$.
Note that $h_{1}m_{1}-m_{1}+1=4837$ and $h_{2}m_{2}-m_{2}+1=2047$.
Let $p\in \{1,2\}$, and suppose that $|M_{4}^{w}(v)\cup M_{5}^{w}(v)|\geq h_{p}m_{p}-m_{p}+1$.
It suffices to show that $T_{8}^{(p)}\prec G$.

By Lemma~\ref{lem-Mk-1-cor}, we may assume that $|M_{4}^{w}(v)|\leq h_{p}-1$.
Then $|M_{5}^{w}(v)-M_{4}^{w}(v)|\geq h_{p}m_{p}-m_{p}+1-(h_{p}-1)=(h_{p}-1)(m_{p}-1)+1$.
Take $(h_{p}-1)(m_{p}-1)+1$ distinct vertices $a_{1},\ldots ,a_{(h_{p}-1)(m_{p}-1)+1}$ in $M_{5}^{w}(v)-M_{4}^{w}(v)$.
Then for $i~(1\leq i\leq (h_{p}-1)(m_{p}-1)+1)$, there exists a $(v,w;5)$-path $Q_{i}=va_{i}b_{i}c_{i}w$ of $G$.
If $|\{c_{i}:1\leq i\leq (h_{p}-1)(m_{p}-1)+1\}|\leq m_{p}-1$, then there exists a vertex $w'$ of $G$ such that
$$
|M_{4}^{w'}(v)|\geq |\{a_{i}:1\leq i\leq (h_{p}-1)(m_{p}-1)+1,~c_{i}=w'\}|\geq \left\lceil \frac{(h_{p}-1)(m_{p}-1)+1}{m_{p}-1}\right\rceil ~(=h_{p}),
$$
and hence by Lemma~\ref{lem-Mk-1-cor}, we have $T_{8}^{(p)}\prec G$.
Thus we may assume that $|\{c_{i}:1\leq i\leq (h_{p}-1)(m_{p}-1)+1\}|\geq m_{p}$.
Relabeling the paths $Q_{i}$, we may also assume that $c_{i}\neq c_{j}$ for $1\leq i<j\leq m_{p}$.
By Lemma~\ref{lem-Mk-fund1}(ii), $Q_{1},\ldots ,Q_{m_{p}}$ are pairwise internally disjoint.
Furthermore, by Lemma~\ref{lem-Mk-fund1}(v), we have
\begin{align}
E\left(G\left[\bigcup _{1\leq i\leq m_{p}}V(Q_{i})\right]\right)-\left(\bigcup _{1\leq i\leq m_{p}}E(Q_{i})\right)=E(G[\{b_{i}:1\leq i\leq m_{p}\}]).\label{cond-lem-Mk-2-1}
\end{align}

\medskip
\noindent
\textbf{Case 1:} $p=1$, i.e., $m_{p}=124$.

Take $x_{1}\in N(a_{1})-\{v,b_{1}\}$ and $y_{1}\in N(b_{1})-\{a_{1},c_{1}\}$.
Since $G$ is $\{C_{3},C_{4}\}$-free, we have $x_{1}\neq y_{1}$.
By (\ref{cond-lem-Mk-2-1}), $x_{1}\notin \{w,c_{1},a_{i},b_{i},c_{i}:2\leq i\leq 124\}$ and $y_{1}\notin \{v,w,a_{i},c_{i}:2\leq i\leq 124\}$.
Since $b_{1},\ldots ,b_{124}$ are distinct vertices, we may assume that $y_{1}\notin \{b_{i}:1\leq i\leq 123\}$.
Consequently,
\begin{align}
x_{1}\neq y_{1}\mbox{ and }\{x_{1},y_{1}\}\cap \left(\bigcup _{1\leq i\leq 123}V(Q_{i})\right)=\emptyset .\label{cond-lem-Mk-2-case1-1}
\end{align}

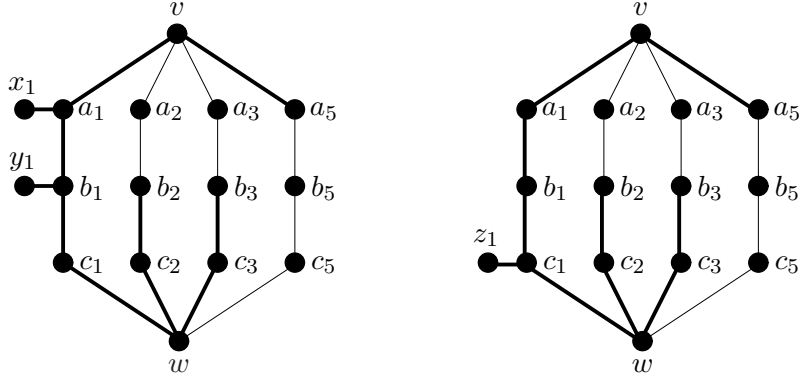
\begin{figure}
\begin{center}
{\unitlength 0.1in%
\begin{picture}(40.6000,18.6500)(15.9000,-22.6500)%
\put(26.0000,-23.3000){\makebox(0,0){$w$}}%
%
\special{sh 1.000}%
\special{ia 2590 595 50 50 0.0000000 6.2831853}%
\special{pn 8}%
\special{ar 2590 595 50 50 0.0000000 6.2831853}%
\put(25.9000,-4.6500){\makebox(0,0){$v$}}%
%
\special{sh 1.000}%
\special{ia 2600 2200 50 50 0.0000000 6.2831853}%
\special{pn 8}%
\special{ar 2600 2200 50 50 0.0000000 6.2831853}%
%
\special{sh 1.000}%
\special{ia 2000 990 50 50 0.0000000 6.2831853}%
\special{pn 8}%
\special{ar 2000 990 50 50 0.0000000 6.2831853}%
%
\special{sh 1.000}%
\special{ia 2000 1390 50 50 0.0000000 6.2831853}%
\special{pn 8}%
\special{ar 2000 1390 50 50 0.0000000 6.2831853}%
%
\special{sh 1.000}%
\special{ia 2000 1790 50 50 0.0000000 6.2831853}%
\special{pn 8}%
\special{ar 2000 1790 50 50 0.0000000 6.2831853}%
%
\special{sh 1.000}%
\special{ia 2400 990 50 50 0.0000000 6.2831853}%
\special{pn 8}%
\special{ar 2400 990 50 50 0.0000000 6.2831853}%
%
\special{sh 1.000}%
\special{ia 2400 1390 50 50 0.0000000 6.2831853}%
\special{pn 8}%
\special{ar 2400 1390 50 50 0.0000000 6.2831853}%
%
\special{sh 1.000}%
\special{ia 2400 1790 50 50 0.0000000 6.2831853}%
\special{pn 8}%
\special{ar 2400 1790 50 50 0.0000000 6.2831853}%
%
\special{sh 1.000}%
\special{ia 2800 1790 50 50 0.0000000 6.2831853}%
\special{pn 8}%
\special{ar 2800 1790 50 50 0.0000000 6.2831853}%
%
\special{sh 1.000}%
\special{ia 2800 1390 50 50 0.0000000 6.2831853}%
\special{pn 8}%
\special{ar 2800 1390 50 50 0.0000000 6.2831853}%
%
\special{sh 1.000}%
\special{ia 2800 990 50 50 0.0000000 6.2831853}%
\special{pn 8}%
\special{ar 2800 990 50 50 0.0000000 6.2831853}%
%
\special{sh 1.000}%
\special{ia 3200 990 50 50 0.0000000 6.2831853}%
\special{pn 8}%
\special{ar 3200 990 50 50 0.0000000 6.2831853}%
%
\special{sh 1.000}%
\special{ia 3200 1390 50 50 0.0000000 6.2831853}%
\special{pn 8}%
\special{ar 3200 1390 50 50 0.0000000 6.2831853}%
%
\special{sh 1.000}%
\special{ia 3200 1790 50 50 0.0000000 6.2831853}%
\special{pn 8}%
\special{ar 3200 1790 50 50 0.0000000 6.2831853}%
%
\special{sh 1.000}%
\special{ia 1800 990 50 50 0.0000000 6.2831853}%
\special{pn 8}%
\special{ar 1800 990 50 50 0.0000000 6.2831853}%
%
\special{sh 1.000}%
\special{ia 1800 1390 50 50 0.0000000 6.2831853}%
\special{pn 8}%
\special{ar 1800 1390 50 50 0.0000000 6.2831853}%
%
\special{pn 20}%
\special{pa 1800 1390}%
\special{pa 2000 1390}%
\special{fp}%
\special{pa 2400 1390}%
\special{pa 2400 1790}%
\special{fp}%
\special{pa 2400 1790}%
\special{pa 2600 2190}%
\special{fp}%
\special{pa 2600 2190}%
\special{pa 2800 1790}%
\special{fp}%
\special{pa 2800 1790}%
\special{pa 2800 1390}%
\special{fp}%
\special{pa 2600 2190}%
\special{pa 2000 1790}%
\special{fp}%
\special{pa 2000 1790}%
\special{pa 2000 990}%
\special{fp}%
\special{pa 2000 990}%
\special{pa 2600 590}%
\special{fp}%
\special{pa 2600 590}%
\special{pa 3200 990}%
\special{fp}%
\special{pa 2000 990}%
\special{pa 1800 990}%
\special{fp}%
%
\special{pn 4}%
\special{pa 2600 600}%
\special{pa 2400 1000}%
\special{fp}%
\special{pa 2400 1000}%
\special{pa 2400 1400}%
\special{fp}%
\special{pa 2800 1400}%
\special{pa 2800 1000}%
\special{fp}%
\special{pa 2800 1000}%
\special{pa 2600 600}%
\special{fp}%
\special{pa 3200 1000}%
\special{pa 3200 1800}%
\special{fp}%
\special{pa 3200 1800}%
\special{pa 2600 2200}%
\special{fp}%
\put(21.5000,-14.0000){\makebox(0,0){$b_{1}$}}%
\put(25.5000,-14.0000){\makebox(0,0){$b_{2}$}}%
\put(29.5000,-14.0000){\makebox(0,0){$b_{3}$}}%
\put(33.5000,-14.0000){\makebox(0,0){$b_{5}$}}%
\put(21.5000,-10.0000){\makebox(0,0){$a_{1}$}}%
\put(25.5000,-10.0000){\makebox(0,0){$a_{2}$}}%
\put(29.5000,-10.0000){\makebox(0,0){$a_{3}$}}%
\put(33.5000,-10.0000){\makebox(0,0){$a_{5}$}}%
\put(33.5000,-18.0000){\makebox(0,0){$c_{5}$}}%
\put(29.5000,-18.0000){\makebox(0,0){$c_{3}$}}%
\put(25.5000,-18.0000){\makebox(0,0){$c_{2}$}}%
\put(21.5000,-18.0000){\makebox(0,0){$c_{1}$}}%
\put(17.9000,-8.6500){\makebox(0,0){$x_{1}$}}%
\put(17.9000,-12.6000){\makebox(0,0){$y_{1}$}}%
\put(50.0000,-23.3000){\makebox(0,0){$w$}}%
%
\special{sh 1.000}%
\special{ia 4990 595 50 50 0.0000000 6.2831853}%
\special{pn 8}%
\special{ar 4990 595 50 50 0.0000000 6.2831853}%
\put(49.9000,-4.6500){\makebox(0,0){$v$}}%
%
\special{sh 1.000}%
\special{ia 5000 2200 50 50 0.0000000 6.2831853}%
\special{pn 8}%
\special{ar 5000 2200 50 50 0.0000000 6.2831853}%
%
\special{sh 1.000}%
\special{ia 4400 990 50 50 0.0000000 6.2831853}%
\special{pn 8}%
\special{ar 4400 990 50 50 0.0000000 6.2831853}%
%
\special{sh 1.000}%
\special{ia 4400 1390 50 50 0.0000000 6.2831853}%
\special{pn 8}%
\special{ar 4400 1390 50 50 0.0000000 6.2831853}%
%
\special{sh 1.000}%
\special{ia 4400 1790 50 50 0.0000000 6.2831853}%
\special{pn 8}%
\special{ar 4400 1790 50 50 0.0000000 6.2831853}%
%
\special{sh 1.000}%
\special{ia 4800 990 50 50 0.0000000 6.2831853}%
\special{pn 8}%
\special{ar 4800 990 50 50 0.0000000 6.2831853}%
%
\special{sh 1.000}%
\special{ia 4800 1390 50 50 0.0000000 6.2831853}%
\special{pn 8}%
\special{ar 4800 1390 50 50 0.0000000 6.2831853}%
%
\special{sh 1.000}%
\special{ia 4800 1790 50 50 0.0000000 6.2831853}%
\special{pn 8}%
\special{ar 4800 1790 50 50 0.0000000 6.2831853}%
%
\special{sh 1.000}%
\special{ia 5200 1790 50 50 0.0000000 6.2831853}%
\special{pn 8}%
\special{ar 5200 1790 50 50 0.0000000 6.2831853}%
%
\special{sh 1.000}%
\special{ia 5200 1390 50 50 0.0000000 6.2831853}%
\special{pn 8}%
\special{ar 5200 1390 50 50 0.0000000 6.2831853}%
%
\special{sh 1.000}%
\special{ia 5200 990 50 50 0.0000000 6.2831853}%
\special{pn 8}%
\special{ar 5200 990 50 50 0.0000000 6.2831853}%
%
\special{sh 1.000}%
\special{ia 5600 990 50 50 0.0000000 6.2831853}%
\special{pn 8}%
\special{ar 5600 990 50 50 0.0000000 6.2831853}%
%
\special{sh 1.000}%
\special{ia 5600 1390 50 50 0.0000000 6.2831853}%
\special{pn 8}%
\special{ar 5600 1390 50 50 0.0000000 6.2831853}%
%
\special{sh 1.000}%
\special{ia 5600 1790 50 50 0.0000000 6.2831853}%
\special{pn 8}%
\special{ar 5600 1790 50 50 0.0000000 6.2831853}%
%
\special{sh 1.000}%
\special{ia 4200 1790 50 50 0.0000000 6.2831853}%
\special{pn 8}%
\special{ar 4200 1790 50 50 0.0000000 6.2831853}%
%
\special{pn 4}%
\special{pa 5000 600}%
\special{pa 4800 1000}%
\special{fp}%
\special{pa 4800 1000}%
\special{pa 4800 1400}%
\special{fp}%
\special{pa 5200 1400}%
\special{pa 5200 1000}%
\special{fp}%
\special{pa 5200 1000}%
\special{pa 5000 600}%
\special{fp}%
\special{pa 5600 1000}%
\special{pa 5600 1800}%
\special{fp}%
\special{pa 5600 1800}%
\special{pa 5000 2200}%
\special{fp}%
\put(45.5000,-14.0000){\makebox(0,0){$b_{1}$}}%
\put(49.5000,-14.0000){\makebox(0,0){$b_{2}$}}%
\put(53.5000,-14.0000){\makebox(0,0){$b_{3}$}}%
\put(57.5000,-14.0000){\makebox(0,0){$b_{5}$}}%
\put(45.5000,-10.0000){\makebox(0,0){$a_{1}$}}%
\put(49.5000,-10.0000){\makebox(0,0){$a_{2}$}}%
\put(53.5000,-10.0000){\makebox(0,0){$a_{3}$}}%
\put(57.5000,-10.0000){\makebox(0,0){$a_{5}$}}%
\put(57.5000,-18.0000){\makebox(0,0){$c_{5}$}}%
\put(53.5000,-18.0000){\makebox(0,0){$c_{3}$}}%
\put(49.5000,-18.0000){\makebox(0,0){$c_{2}$}}%
\put(45.5000,-18.0000){\makebox(0,0){$c_{1}$}}%
\put(41.9000,-16.6000){\makebox(0,0){$z_{1}$}}%
%
\special{pn 20}%
\special{pa 4990 2200}%
\special{pa 5190 1800}%
\special{fp}%
\special{pa 5190 1800}%
\special{pa 5190 1400}%
\special{fp}%
\special{pa 4790 1400}%
\special{pa 4790 1800}%
\special{fp}%
\special{pa 4790 1800}%
\special{pa 4990 2200}%
\special{fp}%
\special{pa 4990 2200}%
\special{pa 4390 1800}%
\special{fp}%
\special{pa 4390 1800}%
\special{pa 4390 1000}%
\special{fp}%
\special{pa 4390 1000}%
\special{pa 4990 600}%
\special{fp}%
\special{pa 4990 600}%
\special{pa 5590 1000}%
\special{fp}%
\special{pa 4390 1800}%
\special{pa 4190 1800}%
\special{fp}%
\end{picture}}%

\caption{Induced copies of $T_{8}^{(1)}$ and $T_{8}^{(2)}$ in Lemma~\ref{lem-Mk-2}}
\label{flem4.4}
\end{center}
\end{figure}

If $x_{1}w\in E(G)$, then $va_{1}x_{1}w$ is a $(v,w;4)$-path, and hence $a_{1}\in M_{4}^{w}(v)$, which is a contradiction.
This together with the fact that $G$ is $\{C_{3},C_{4}\}$-free implies that
\begin{align}
N(x_{1})\cap \{v,w,y_{1},b_{1},c_{1},a_{i}:2\leq i\leq 123\}=\emptyset \mbox{ and }N(y_{1})\cap \{v,w,a_{1},c_{1}\}=\emptyset .\label{cond-lem-Mk-2-case1-2}
\end{align}
By (\ref{cond-lem-Mk-2-case1-2}) and the fact that $b_{1}w\notin E(G)$, the set $M_{4}^{w}(u)$ is well-defined for $u\in \{b_{1},x_{1},y_{1}\}$.
Considering Lemma~\ref{lem-Mk-1-cor}(i), we may assume that
\begin{align}
|N(u)\cap \{b_{i}:2\leq i\leq 123\}|\leq |M_{4}^{w}(u)|\leq 39\mbox{ for }u\in \{b_{1},x_{1},y_{1}\}.\label{cond-lem-Mk-2-case1-4}
\end{align}
By (\ref{cond-lem-Mk-2-1}), we have
\begin{align}
N(b_{1})\cap \{c_{i}:2\leq i\leq 123\}=\emptyset .\label{cond-lem-Mk-2-case1-5}
\end{align}
Furthermore, we have
\begin{align}
|N(u)\cap \{c_{i}:2\leq i\leq 123\}|\leq |N(u)\cap N(w)|\leq 1\mbox{ for }u\in \{x_{1},y_{1}\}.\label{cond-lem-Mk-2-case1-6}
\end{align}
By (\ref{cond-lem-Mk-2-case1-4})--(\ref{cond-lem-Mk-2-case1-6}),
$$
|N(\{b_{1},x_{1},y_{1}\})\cap \{b_{i},c_{i}:2\leq i\leq 123\}|\leq \sum _{u\in \{b_{1},x_{1},y_{1}\}}|N(u)\cap \{b_{i},c_{i}:2\leq i\leq 123\}|\leq 3\cdot 39+2.
$$
We may assume that $E_{G}(\{b_{1},x_{1},y_{1}\},\{b_{i},c_{i}:2\leq i\leq 4\})=\emptyset $.
By (\ref{cond-lem-Mk-2-case1-2}), we have
\begin{align}
N(x_{1})\cap \{v,w,y_{1},b_{1},c_{1},a_{i},b_{i},c_{i}:2\leq i\leq 4\}=\emptyset \mbox{ and }N(y_{1})\cap \{v,w,a_{1},c_{1},b_{i},c_{i}:2\leq i\leq 4\}=\emptyset .\label{cond-lem-Mk-2-case1-7}
\end{align}
Since $G[\{b_{2},b_{3},b_{4}\}]$ is not complete, we may assume that $b_{2}b_{3}\notin E(G)$, which implies that
\begin{align}
\{b_{1},b_{2},b_{3}\}\mbox{ is an independent set of }G.\label{cond-lem-Mk-2-case1-8}
\end{align}
Since $|N(y_{1})\cap \{a_{5},a_{6}\}|\leq |N(y_{1})\cap N(v)|\leq 1$, we may assume that
\begin{align}
y_{1}a_{5}\notin E(G).\label{cond-lem-Mk-2-case1-9}
\end{align}
Consequently, it follows from (\ref{cond-lem-Mk-2-1}), (\ref{cond-lem-Mk-2-case1-1}) and (\ref{cond-lem-Mk-2-case1-7})--(\ref{cond-lem-Mk-2-case1-9}) that $\{b_{2},c_{2},w,b_{3},c_{3},c_{1},b_{1},y_{1},a_{1},x_{1},v,a_{5}\}$ induces a copy of $T_{8}^{(1)}$ in $G$ (see the left graph in Figure~\ref{flem4.4}).

\medskip
\noindent
\textbf{Case 2:} $p=2$, i.e., $m_{p}=66$.

For $i~(1\leq i\leq 66)$, take $z_{i}\in N(c_{i})-\{w,b_{i}\}$.
Since $G$ is $\{C_{3},C_{4}\}$-free, $z_{1},\ldots ,z_{66}$ are distinct.
By (\ref{cond-lem-Mk-2-1}), we have
\begin{align}
\{z_{i}:1\leq i\leq 66\}\cap \left(\bigcup _{1\leq i\leq 66}V(Q_{i})\right)=\emptyset .\label{cond-lem-Mk-2-case2-1}
\end{align}

If $\{z_{i}:1\leq i\leq 66\}\subseteq N(v)$, then $|M_{4}^{w}(v)|\geq 66>32$ because $vz_{i}c_{i}w~(1\leq i\leq 66)$ are distinct $(v,w;4)$-paths of $G$, which contradicts the assumption in the second paragraph of the proof.
Thus, without loss of generality, we may assume that $vz_{1}\notin E(G)$.
This together with the $\{C_{3},C_{4}\}$-freeness of $G$ implies that
\begin{align}
N(z_{1})\cap \{v,w,a_{1},b_{1},c_{i}:2\leq i\leq 66\}=\emptyset .\label{cond-lem-Mk-2-case2-2}
\end{align}
By (\ref{cond-lem-Mk-2-case2-2}) and the fact that $b_{1}w\notin E(G)$, the set $M_{4}^{w}(u)$ is well-defined for $u\in \{b_{1},z_{1}\}$.
Considering Lemma~\ref{lem-Mk-1-cor}(ii), we may assume that
\begin{align}
|N(u)\cap \{b_{i}:2\leq i\leq 66\}|\leq |M_{4}^{w}(u)|\leq 31\mbox{ for }u\in \{b_{1},z_{1}\}.\label{cond-lem-Mk-2-case2-4}
\end{align}
By (\ref{cond-lem-Mk-2-case2-4}), 
$$
|N(\{b_{1},z_{1}\})\cap \{b_{i}:2\leq i\leq 66\}|\leq 2\cdot 31.
$$
We may assume that $N(\{b_{1},z_{1}\})\cap \{b_{i}:2\leq i\leq 4\}=\emptyset $.
This together with (\ref{cond-lem-Mk-2-1}) and (\ref{cond-lem-Mk-2-case2-2}) leads to
\begin{align}
N(b_{1})\cap \{v,w,a_{i},b_{i},c_{i}:2\leq i\leq 4\}=\emptyset \mbox{ and }N(z_{1})\cap \{v,w,a_{1},b_{1},b_{i},c_{i}:2\leq i\leq 4\}=\emptyset .\label{cond-lem-Mk-2-case2-5}
\end{align}
Since $G[\{b_{2},b_{3},b_{4}\}]\not\simeq C_{3}$, we may assume that
\begin{align}
b_{2}b_{3}\notin E(G).\label{cond-lem-Mk-2-case2-6}
\end{align}
Since $|N(z_{1})\cap \{a_{5},a_{6}\}|\leq |N(z_{1})\cap N(v)|\leq 1$, we may assume that
\begin{align}
z_{1}a_{5}\notin E(G).\label{cond-lem-Mk-2-case2-7}
\end{align}
Consequently, it follows from (\ref{cond-lem-Mk-2-1}), (\ref{cond-lem-Mk-2-case2-1}) and (\ref{cond-lem-Mk-2-case2-5})--(\ref{cond-lem-Mk-2-case2-7}) that $\{b_{2},c_{2},w,b_{3},c_{3},c_{1},z_{1},b_{1},a_{1},v,a_{4}\}$ induces a copy of $T_{8}^{(2)}$ in $G$ (see the right graph in Figure~\ref{flem4.4}).
\qed

\begin{figure}
\begin{center}
{\unitlength 0.1in%
\begin{picture}(16.0000,18.6500)(16.5000,-22.6500)%
\put(26.0000,-23.3000){\makebox(0,0){$w$}}%
%
\special{sh 1.000}%
\special{ia 2590 595 50 50 0.0000000 6.2831853}%
\special{pn 8}%
\special{ar 2590 595 50 50 0.0000000 6.2831853}%
\put(25.9000,-4.6500){\makebox(0,0){$v$}}%
%
\special{sh 1.000}%
\special{ia 2600 2200 50 50 0.0000000 6.2831853}%
\special{pn 8}%
\special{ar 2600 2200 50 50 0.0000000 6.2831853}%
%
\special{sh 1.000}%
\special{ia 2000 990 50 50 0.0000000 6.2831853}%
\special{pn 8}%
\special{ar 2000 990 50 50 0.0000000 6.2831853}%
%
\special{sh 1.000}%
\special{ia 2000 1390 50 50 0.0000000 6.2831853}%
\special{pn 8}%
\special{ar 2000 1390 50 50 0.0000000 6.2831853}%
%
\special{sh 1.000}%
\special{ia 2000 1790 50 50 0.0000000 6.2831853}%
\special{pn 8}%
\special{ar 2000 1790 50 50 0.0000000 6.2831853}%
%
\special{sh 1.000}%
\special{ia 2400 990 50 50 0.0000000 6.2831853}%
\special{pn 8}%
\special{ar 2400 990 50 50 0.0000000 6.2831853}%
%
\special{sh 1.000}%
\special{ia 2400 1390 50 50 0.0000000 6.2831853}%
\special{pn 8}%
\special{ar 2400 1390 50 50 0.0000000 6.2831853}%
%
\special{sh 1.000}%
\special{ia 2400 1790 50 50 0.0000000 6.2831853}%
\special{pn 8}%
\special{ar 2400 1790 50 50 0.0000000 6.2831853}%
%
\special{sh 1.000}%
\special{ia 2800 1790 50 50 0.0000000 6.2831853}%
\special{pn 8}%
\special{ar 2800 1790 50 50 0.0000000 6.2831853}%
%
\special{sh 1.000}%
\special{ia 2800 1390 50 50 0.0000000 6.2831853}%
\special{pn 8}%
\special{ar 2800 1390 50 50 0.0000000 6.2831853}%
%
\special{sh 1.000}%
\special{ia 2800 990 50 50 0.0000000 6.2831853}%
\special{pn 8}%
\special{ar 2800 990 50 50 0.0000000 6.2831853}%
%
\special{sh 1.000}%
\special{ia 3200 990 50 50 0.0000000 6.2831853}%
\special{pn 8}%
\special{ar 3200 990 50 50 0.0000000 6.2831853}%
%
\special{sh 1.000}%
\special{ia 3200 1390 50 50 0.0000000 6.2831853}%
\special{pn 8}%
\special{ar 3200 1390 50 50 0.0000000 6.2831853}%
%
\special{sh 1.000}%
\special{ia 3200 1790 50 50 0.0000000 6.2831853}%
\special{pn 8}%
\special{ar 3200 1790 50 50 0.0000000 6.2831853}%
\put(18.5000,-10.0000){\makebox(0,0){$a_{1}$}}%
%
\special{pn 20}%
\special{pa 2600 590}%
\special{pa 2000 990}%
\special{fp}%
\special{pa 2000 990}%
\special{pa 2000 1390}%
\special{fp}%
\special{pa 2400 1390}%
\special{pa 2400 990}%
\special{fp}%
\special{pa 2400 990}%
\special{pa 2600 590}%
\special{fp}%
\special{pa 2600 590}%
\special{pa 2800 990}%
\special{fp}%
\special{pa 2800 990}%
\special{pa 2800 1790}%
\special{fp}%
\special{pa 2800 1790}%
\special{pa 2600 2190}%
\special{fp}%
\special{pa 2600 2190}%
\special{pa 3200 1790}%
\special{fp}%
\special{pa 3200 1790}%
\special{pa 3200 1390}%
\special{fp}%
%
\special{pn 4}%
\special{pa 3200 1390}%
\special{pa 3200 990}%
\special{fp}%
\special{pa 3200 990}%
\special{pa 2600 590}%
\special{fp}%
\special{pa 2400 1390}%
\special{pa 2400 1390}%
\special{fp}%
\special{pa 2400 1790}%
\special{pa 2400 1390}%
\special{fp}%
\special{pa 2400 1790}%
\special{pa 2600 2190}%
\special{fp}%
\special{pa 2600 2190}%
\special{pa 2000 1790}%
\special{fp}%
\special{pa 2000 1790}%
\special{pa 2000 1390}%
\special{fp}%
\put(18.5000,-14.0000){\makebox(0,0){$b_{1}$}}%
\put(18.5000,-18.0000){\makebox(0,0){$c_{1}$}}%
\put(22.5000,-10.0000){\makebox(0,0){$a_{2}$}}%
\put(22.5000,-14.0000){\makebox(0,0){$b_{2}$}}%
\put(22.5000,-18.0000){\makebox(0,0){$c_{2}$}}%
\put(26.5000,-18.0000){\makebox(0,0){$c_{3}$}}%
\put(26.5000,-14.0000){\makebox(0,0){$b_{3}$}}%
\put(26.5000,-10.0000){\makebox(0,0){$a_{3}$}}%
\put(30.5000,-10.0000){\makebox(0,0){$a_{4}$}}%
\put(30.5000,-14.0000){\makebox(0,0){$b_{4}$}}%
\put(30.5000,-18.0000){\makebox(0,0){$c_{4}$}}%
\end{picture}}%

\caption{An induced copy of $T_{9}$ in Lemma~\ref{lem-Mk-3} and Lemma~\ref{lem-Mk-5}}
\label{flem4.5}
\end{center}
\end{figure}
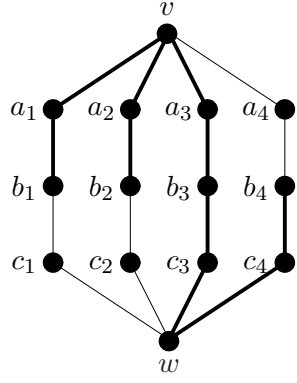

\begin{lem}
\label{lem-Mk-3}
Let $v$ and $w$ be non-adjacent vertices of $G$.
If there are $25$ pairwise internally disjoint $(v,w;5)$-paths of $G$, then $T_{9}\prec G$.
\end{lem}
\proof
Let $Q_{i}=va_{i}b_{i}c_{i}w~(1\leq i\leq 25)$ be $25$ pairwise internally disjoint $(v,w;5)$-paths of $G$.
Then by Lemma~\ref{lem-Mk-fund1}(iv), we have
\begin{align}
E\left(G\left[\bigcup _{1\leq i\leq 25}V(Q_{i})\right]\right)&-\left(\bigcup _{1\leq i\leq 25}E(Q_{i})\right)\nonumber \\
&= E(G[\{b_{i}:1\leq i\leq 25\}])\cup E_{G}(\{a_{i}:1\leq i\leq 25\},\{c_{i}:1\leq i\leq 25\}).\label{cond-lem-Mk-3-1}
\end{align}
Now we define a graph $H$ on $\{1,\ldots ,25\}$ by joining $i$ and $j$ with $i\neq j$ if and only if $a_{i}c_{j}\in E(G)$ or $a_{j}c_{i}\in E(G)$.
Since $|N(a_{i})\cap \{c_{j}:1\leq j\leq 25\}|\leq |N(a_{i})\cap N(w)|\leq 1$ and $|N(c_{i})\cap \{a_{j}:1\leq j\leq 25\}|\leq |N(c_{i})\cap N(v)|\leq 1$ for all $1\leq i\leq 25$, we have $\Delta (H)\leq 2$.
Hence there exists an independent set $I\subseteq \{1,\ldots ,25\}$ of $H$ with $|I|=\lceil \frac{25}{3}\rceil =9$.
Then $a_{i}c_{j}\notin E(G)$ for all $i,j\in I$.
Since $G$ is $C_{3}$-free and $R(3,4)=9$, there exists $I'\subseteq I$ with $|I'|=4$ such that $\{b_{i}:i\in I'\}$ is an independent set of $G$.
We may assume that $I'=\{1,2,3,4\}$.
Then it follows from (\ref{cond-lem-Mk-3-1}) that $E(G[\bigcup _{1\leq i\leq 4}V(Q_{i})])=\bigcup _{1\leq i\leq 4}E(Q_{i})$, and hence $\{b_{1},a_{1},v,b_{2},a_{2},a_{3},b_{3},c_{3},w,c_{4},b_{4}\}$ induces a copy of $T_{9}$ in $G$ (see Figure~\ref{flem4.5}).
\qed

\begin{lem}
\label{lem-Mk-4}
Let $v$ and $w$ be non-adjacent vertices of $G$.
If $|M_{4}^{w}(v)|\geq 75$, then $T_{9}\prec G$.
\end{lem}
\proof
Let $a_{1},\ldots ,a_{75}$ be distinct vertices in $M_{4}^{w}(v)$.
Then for each $i~(1\leq i\leq 75)$, there exists a vertex $b_{i}\in V(G)$ such that $Q_{i}=ua_{i}b_{i}v$ is a $(v,w;4)$-path.
By Lemma~\ref{lem-Mk-fund1}(i)(iii), $\{a_{i},b_{i}\}\cap \{a_{j},b_{j}\}=\emptyset $ and $E_{G}(\{a_{i},b_{i}\},\{a_{j},b_{j}\})=\emptyset $ for $1\leq i<j\leq 75$.
For $i~(1\leq i\leq 13)$, take $x_{i}\in N(a_{i})-\{v,b_{i}\}$.
Then
$$
|N(\{x_{i}:1\leq i\leq 13\})\cap \{b_{i}:14\leq i\leq 75\}|\leq \sum _{1\leq i\leq 13}|N(x_{i})\cap N(w)|\leq 13.
$$
We may assume that
\begin{align}
N(\{x_{i}:1\leq i\leq 13\})\cap \{b_{i}:14\leq i\leq 62\}=\emptyset .\label{cond-lem-Mk-4-0}
\end{align}
For each $i~(14\leq i\leq 62)$, take $y_{i}\in N(b_{i})-\{w,a_{i}\}$.
By  Lemma~\ref{lem-Mk-fund1}(i)(iii), (\ref{cond-lem-Mk-4-0}) and the $\{C_{3},C_{4}\}$-freeness of $G$,
\begin{align}
&x_{1},\ldots ,x_{13},y_{14},\ldots ,y_{62}\mbox{ are distinct};\label{cond-lem-Mk-4-1}\\
&(\{x_{i}:1\leq i\leq 13\}\cup \{y_{i}:14\leq i\leq 62\})\cap \left(\bigcup _{1\leq i\leq 62}V(Q_{i})\right)=\emptyset ;\label{cond-lem-Mk-4-2}\\
&N(x_{i})\cap \{v,w,a_{j}:1\leq j\leq 62,~j\neq i\}=\emptyset \mbox{ for }1\leq i\leq 13;\mbox{ and }\label{cond-lem-Mk-4-3}
\end{align}
\begin{align}
N(y_{i})\cap \{v,w,b_{j}:1\leq j\leq 62,~j\neq i\}=\emptyset \mbox{ for }14\leq i\leq 62.\label{cond-lem-Mk-4-4}
\end{align}

For $i~(14\leq i\leq 62)$, since $|N(y_{i})\cap N(v)|\leq 1$, there exists a vertex $z_{i}\in N(y_{i})-(\{b_{i}\}\cup N(v))$.
Then $\{z_{i}:14\leq i\leq 62\}\cap \{a_{i}:1\leq i\leq 62\}=\emptyset $.
Since $N(y_{i})\cap N(w)=\{b_{i}\}$, we have $z_{i}\notin N(w)$, and hence $\{z_{i}:14\leq i\leq 62\}\cap \{b_{i}:1\leq i\leq 62\}=\emptyset $.
By (\ref{cond-lem-Mk-4-4}), we have $z_{i}\notin \{v,w\}$ for $14\leq i\leq 62$.
Consequently,
\begin{align}
\{z_{i}:14\leq i\leq 62\}\cap \left(\bigcup _{1\leq i\leq 62}V(Q_{i})\right)=\emptyset .\label{cond-lem-Mk-4-5}
\end{align}
By the choice of $z_{i}$ and the $\{C_{3},C_{4}\}$-freeness of $G$, we have
\begin{align}
N(z_{i})\cap \{v,w,b_{i}\}=\emptyset \mbox{ for }14\leq i\leq 62.\label{cond-lem-Mk-4-6}
\end{align}

Assume for the moment that $|\{z_{i}:14\leq i\leq 62\}|\leq 2$.
Then there exists a vertex $z$ of $G$ such that $|\{i:14\leq i\leq 62,~z_{i}=z\}|\geq \lceil \frac{49}{2} \rceil =25$.
We may assume that $z_{i}=z$ for all $i~(14\leq i\leq 38)$.
Then the paths $va_{i}b_{i}y_{i}z~(14\leq i\leq 38)$ are $25$ pairwise internally disjoint $(v,z;5)$-paths of $G$, and hence it follows from Lemma~\ref{lem-Mk-4} that $T_{9}\prec G$.
Thus we may assume that $|\{z_{i}:14\leq i\leq 62\}|\geq 3$.
We may assume that $z_{14},z_{15},z_{16}$ are distinct.

\begin{claim}
\label{cl-lem-Mk-4-1}
There exist $i,j\in \{14,15,16\}$ with $i\neq j$ such that $\{y_{i},z_{i}\}\cap \{y_{j},z_{j}\}=\emptyset $.
\end{claim}
\proof
By (\ref{cond-lem-Mk-4-1}), $y_{14},y_{15},y_{16}$ are distinct.
Since $G[\{y_{14},y_{15},y_{16}\}]\not\simeq C_{3}$, there exist $14\leq i<j\leq 16$ such that $y_{i}y_{j}\notin E(G)$.
This implies that $z_{i}\neq y_{j}$ and $z_{j}\neq y_{i}$.
From $y_{i}\neq y_{j}$ and $z_{i}\neq z_{j}$, it follows that $\{y_{i},z_{i}\}\cap \{y_{j},z_{j}\}=\emptyset $.
\qed

By Claim~\ref{cl-lem-Mk-4-1}, we may assume that
\begin{align}
\{y_{14},z_{14}\}\cap \{y_{15},z_{15}\}=\emptyset .\label{cond-lem-Mk-4-6++}
\end{align}
Since $|N(u)\cap \{a_{i}:1\leq i\leq 13\}|\leq |N(u)\cap N(v)|\leq 1$ for all $u\in \{y_{14},y_{15},z_{14},z_{15}\}$, we may assume that 
\begin{align}
N(\{y_{14},y_{15},z_{14},z_{15}\})\cap \{a_{i}:1\leq i\leq 9\}=\emptyset .\label{cond-lem-Mk-4-6+}
\end{align}
This together with (\ref{cond-lem-Mk-4-1}) and (\ref{cond-lem-Mk-4-6++}) implies that
\begin{align}
x_{1},\ldots ,x_{9},y_{14},y_{15},z_{14},z_{15}\mbox{ are distinct}.\label{cond-lem-Mk-4-7}
\end{align}

\begin{claim}
\label{cl-lem-Mk-4-2}
There exists $i_{0}\in \{14,15\}$ and there exists $J\subseteq \{1,\ldots ,9\}$ with $|J|=3$ such that $N(\{y_{i_{0}},z_{i_{0}}\})\cap \{x_{j}:j\in J\}=\emptyset $.
\end{claim}
\proof
Suppose that $|N(\{y_{i},z_{i}\})\cap \{x_{j}:1\leq j\leq 9\}|\geq 7$ for each $i\in \{14,15\}$.
Then $|N(\{y_{14},z_{14}\})\cap N(\{y_{15},z_{15}\})\cap \{x_{j}:1\leq j\leq 9\}|\geq 5$.
Hence for some $u\in \{y_{14},z_{14}\}$ and some $u'\in \{y_{15},z_{15}\}$, $|N(u)\cap N(u')|\geq |N(u)\cap N(u')\cap \{x_{j}:1\leq j\leq 9\}|\geq 2$, which is a contradiction.
Thus for some $i_{0}\in \{14,15\}$, we have $|N(\{y_{i_{0}},z_{i_{0}}\})\cap \{x_{j}:1\leq j\leq 9\}|\leq 6$.
This implies that $|\{j:1\leq j\leq 9,~x_{j}\notin N(\{y_{i_{0}},z_{i_{0}}\})\}|\geq 3$.
Hence $i_{0}$ and a set $J\subseteq \{j:1\leq j\leq 9,~x_{j}\notin N(\{y_{i_{0}},z_{i_{0}}\})\}$ with $|J|=3$ satisfy the desired property.
\qed

\begin{figure}
\begin{center}
{\unitlength 0.1in%
\begin{picture}(26.0000,18.6500)(4.5000,-22.6500)%
\put(18.0000,-23.3000){\makebox(0,0){$w$}}%
%
\special{sh 1.000}%
\special{ia 1790 595 50 50 0.0000000 6.2831853}%
\special{pn 8}%
\special{ar 1790 595 50 50 0.0000000 6.2831853}%
\put(17.9000,-4.6500){\makebox(0,0){$v$}}%
%
\special{sh 1.000}%
\special{ia 1000 995 50 50 0.0000000 6.2831853}%
\special{pn 8}%
\special{ar 1000 995 50 50 0.0000000 6.2831853}%
%
\special{sh 1.000}%
\special{ia 1600 995 50 50 0.0000000 6.2831853}%
\special{pn 8}%
\special{ar 1600 995 50 50 0.0000000 6.2831853}%
%
\special{sh 1.000}%
\special{ia 1000 1795 50 50 0.0000000 6.2831853}%
\special{pn 8}%
\special{ar 1000 1795 50 50 0.0000000 6.2831853}%
%
\special{sh 1.000}%
\special{ia 1600 1795 50 50 0.0000000 6.2831853}%
\special{pn 8}%
\special{ar 1600 1795 50 50 0.0000000 6.2831853}%
%
\special{sh 1.000}%
\special{ia 1800 2200 50 50 0.0000000 6.2831853}%
\special{pn 8}%
\special{ar 1800 2200 50 50 0.0000000 6.2831853}%
%
\special{sh 1.000}%
\special{ia 1400 1400 50 50 0.0000000 6.2831853}%
\special{pn 8}%
\special{ar 1400 1400 50 50 0.0000000 6.2831853}%
%
\special{sh 1.000}%
\special{ia 800 1400 50 50 0.0000000 6.2831853}%
\special{pn 8}%
\special{ar 800 1400 50 50 0.0000000 6.2831853}%
\put(8.5000,-10.0000){\makebox(0,0){$a_{1}$}}%
\put(14.5000,-10.0000){\makebox(0,0){$a_{2}$}}%
\put(6.5000,-14.0000){\makebox(0,0){$x_{1}$}}%
\put(12.5000,-14.0000){\makebox(0,0){$x_{2}$}}%
\put(8.5000,-18.0000){\makebox(0,0){$b_{1}$}}%
\put(14.5000,-18.0000){\makebox(0,0){$b_{2}$}}%
%
\special{sh 1.000}%
\special{ia 2000 990 50 50 0.0000000 6.2831853}%
\special{pn 8}%
\special{ar 2000 990 50 50 0.0000000 6.2831853}%
%
\special{sh 1.000}%
\special{ia 2000 1790 50 50 0.0000000 6.2831853}%
\special{pn 8}%
\special{ar 2000 1790 50 50 0.0000000 6.2831853}%
\put(18.5000,-9.9500){\makebox(0,0){$a_{4}$}}%
\put(18.5000,-17.9500){\makebox(0,0){$b_{4}$}}%
%
\special{sh 1.000}%
\special{ia 2600 990 50 50 0.0000000 6.2831853}%
\special{pn 8}%
\special{ar 2600 990 50 50 0.0000000 6.2831853}%
%
\special{sh 1.000}%
\special{ia 2600 1790 50 50 0.0000000 6.2831853}%
\special{pn 8}%
\special{ar 2600 1790 50 50 0.0000000 6.2831853}%
\put(24.2000,-10.0000){\makebox(0,0){$a_{14}$}}%
\put(24.3000,-17.5500){\makebox(0,0){$b_{14}$}}%
%
\special{sh 1.000}%
\special{ia 2800 1390 50 50 0.0000000 6.2831853}%
\special{pn 8}%
\special{ar 2800 1390 50 50 0.0000000 6.2831853}%
%
\special{sh 1.000}%
\special{ia 3000 1390 50 50 0.0000000 6.2831853}%
\special{pn 8}%
\special{ar 3000 1390 50 50 0.0000000 6.2831853}%
\put(27.9000,-12.6500){\makebox(0,0){$y_{14}$}}%
\put(29.9000,-12.6000){\makebox(0,0){$z_{14}$}}%
%
\special{pn 20}%
\special{pa 1800 600}%
\special{pa 1000 1000}%
\special{fp}%
\special{pa 1000 1000}%
\special{pa 800 1400}%
\special{fp}%
\special{pa 1400 1400}%
\special{pa 1600 1000}%
\special{fp}%
\special{pa 1600 1000}%
\special{pa 1800 600}%
\special{fp}%
\special{pa 1800 600}%
\special{pa 2000 1000}%
\special{fp}%
\special{pa 2000 1000}%
\special{pa 2000 1800}%
\special{fp}%
\special{pa 2000 1800}%
\special{pa 1800 2200}%
\special{fp}%
\special{pa 1800 2200}%
\special{pa 2600 1800}%
\special{fp}%
\special{pa 2600 1800}%
\special{pa 2800 1400}%
\special{fp}%
\special{pa 2800 1400}%
\special{pa 3000 1400}%
\special{fp}%
%
\special{pn 4}%
\special{pa 2600 1800}%
\special{pa 2600 1000}%
\special{fp}%
\special{pa 2600 1000}%
\special{pa 1800 600}%
\special{fp}%
\special{pa 1600 1000}%
\special{pa 1600 1800}%
\special{fp}%
\special{pa 1600 1800}%
\special{pa 1800 2200}%
\special{fp}%
\special{pa 1800 2200}%
\special{pa 1000 1800}%
\special{fp}%
\special{pa 1000 1800}%
\special{pa 1000 1000}%
\special{fp}%
\end{picture}}%

\caption{An induced copy of $T_{9}$ in Lemma~\ref{lem-Mk-4}}
\label{flem4.6}
\end{center}
\end{figure}
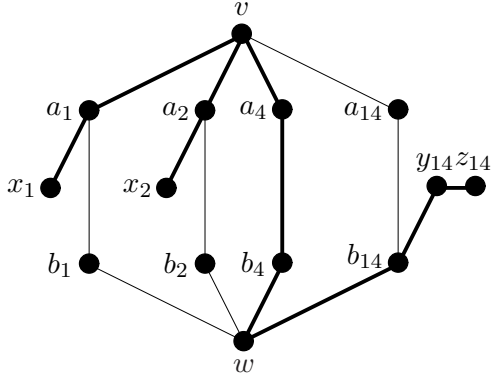

Let $i_{0}$ and $J$ be as in Claim~\ref{cl-lem-Mk-4-2}.
We may assume that $i_{0}=14$ and $J=\{1,2,3\}$.
Then it follows from (\ref{cond-lem-Mk-4-4}) and (\ref{cond-lem-Mk-4-6+}) that
\begin{align}
N(y_{14})\cap \left(\{x_{1},x_{2},x_{3}\}\cup \left(\bigcup _{1\leq i\leq 9}V(Q_{i})\right)\right)=\emptyset \mbox{ and }N(z_{14})\cap \{x_{1},x_{2},x_{3}\}=\emptyset .\label{cond-lem-Mk-4-8}
\end{align}
Since $G[\{x_{1},x_{2},x_{3}\}]\not\simeq C_{3}$, we may assume that
\begin{align}
x_{1}x_{2}\notin E(G).\label{cond-lem-Mk-4-8+}
\end{align}
Since $|N(\{x_{1},x_{2},z_{14}\})\cap \{b_{i}:4\leq i\leq 7\}|\leq \sum _{u\in \{x_{1},x_{2},z_{14}\}}|N(u)\cap N(w)|\leq 3$, we may assume that
\begin{align}
ub_{4}\notin E(G)\mbox{ for }u\in \{x_{1},x_{2},z_{14}\}.\label{cond-lem-Mk-4-9}
\end{align}
By (\ref{cond-lem-Mk-4-0}), (\ref{cond-lem-Mk-4-3}) and (\ref{cond-lem-Mk-4-8})--(\ref{cond-lem-Mk-4-9}),
\begin{align}
N(x_{i})\cap \{v,w,a_{3-i},a_{4},b_{4},b_{14},x_{3-i},y_{14},z_{14}\}=\emptyset \mbox{ for }i\in \{1,2\}.\label{cond-lem-Mk-4-10}
\end{align}
By (\ref{cond-lem-Mk-4-6}), (\ref{cond-lem-Mk-4-6+}), (\ref{cond-lem-Mk-4-8}) and (\ref{cond-lem-Mk-4-9}),
\begin{align}
N(z_{14})\cap \{v,w,a_{1},a_{2},a_{4},b_{4},b_{14},x_{1},x_{2}\}=\emptyset .\label{cond-lem-Mk-4-11}
\end{align}
Consequently, it follows from (\ref{cond-lem-Mk-4-2}), (\ref{cond-lem-Mk-4-5}), (\ref{cond-lem-Mk-4-7}), (\ref{cond-lem-Mk-4-8}), (\ref{cond-lem-Mk-4-10}) and (\ref{cond-lem-Mk-4-11}) that
$$
\{x_{1},a_{1},v,x_{2},a_{2},a_{4},b_{4},w,b_{14},y_{14},z_{14}\}
$$
induces a copy of $T_{9}$ in $G$ (see Figure~\ref{flem4.6}).
\qed

\begin{lem}
\label{lem-Mk-5}
Let $v$ and $w$ be non-adjacent vertices of $G$.
If $|M_{4}^{w}(v)\cup M_{5}^{w}(v)|\geq 667$, then $T_{9}\prec G$.
\end{lem}
\proof
Suppose that $|M_{4}^{w}(v)\cup M_{5}^{w}(v)|\geq 667$.
By Lemma~\ref{lem-Mk-4}, we may assume that $|M_{4}^{w}(v)|\leq 74$.
Then $|M_{5}^{w}(v)-M_{4}^{w}(v)|\geq 667-74=593$.
Take $593$ distinct vertices $a_{1},\ldots ,a_{593}$ in $M_{5}^{w}(v)-M_{4}^{w}(v)$.
Then for each $i~(1\leq i\leq 593)$, there exists a $(v,w;5)$-path $Q_{i}=va_{i}b_{i}c_{i}w$ of $G$.
If $|\{c_{i}:1\leq i\leq 593\}|\leq 8$, then there exists a vertex $w'$ of $G$ such that
$$
|M_{4}^{w'}(v)|\geq |\{a_{i}:1\leq i\leq 593,~c_{i}=w'\}|\geq \left\lceil \frac{593}{8} \right\rceil =75,
$$
and hence by Lemma~\ref{lem-Mk-4}, we have $T_{9}\prec G$.
Thus we may assume that $|\{c_{i}:1\leq i\leq 593\}|\geq 9$.
We may assume that $c_{i}\neq c_{j}$ for $1\leq i<j\leq 9$.
By Lemma~\ref{lem-Mk-fund1}(ii), $Q_{1},\ldots ,Q_{9}$ are pairwise internally disjoint.
Furthermore, by Lemma~\ref{lem-Mk-fund1}(v), we have 
\begin{align}
E\left(G\left[\bigcup _{1\leq i\leq 9}V(Q_{i})\right]\right)-\left(\bigcup _{1\leq i\leq 9}E(Q_{i})\right)=E(G[\{b_{i}:1\leq i\leq 9\}]).\label{cond-lem-Mk-5-1}
\end{align}
Since $G$ is $C_{3}$-free and $R(3,4)=9$, there exists an independent set $X\subseteq \{b_{i}:1\leq i\leq 9\}$ of $G$ with $|X|=4$.
We may assume that $X=\{b_{i}:1\leq i\leq 4\}$.
Then it follows from (\ref{cond-lem-Mk-5-1}) that $\{b_{1},a_{1},v,b_{2},a_{2},a_{3},b_{3},c_{3},w,c_{4},b_{4}\}$ induces a copy of $T_{9}$ in $G$ (see Figure~\ref{flem4.5}).
\qed

\section{Proof of Theorem~\ref{prop-maxdegbound}}\label{sec-maxdeg}

Recall that we have fixed a connected $\{C_{3},C_{4}\}$-free graph $G$ with $\delta (G)\geq 3$.
For a vertex $u$ of $G$ and a non-negative integer $d$, let $N_{d}(u)$ be the set of vertices $v$ of $G$ such that ${\rm dist}_{G}(u,v)=d$, and let $N_{\leq d}(u)=\bigcup _{0\leq i\leq d}N_{i}(u)$ and $N_{\geq d}(u)=\bigcup _{i\geq d}N_{i}(u)$.
Note that $N_{0}(u)=\{u\}$ and $N_{1}(u)=N(u)$.
We further fix a vertex $w$ of $G$ with $d_{G}(w)=\Delta (G)$.
For a subset $U$ of $V(G)$, let $L(U)$ be the set of vertices $v\in N_{2}(w)\cup N_{3}(w)$ such that there exists a $v$-$w$ path $Q$ of $G$ with $V(Q)\cap U=\emptyset $ having order $4$.

\begin{lem}
\label{lem-Xunion}
Let $X$ be a subset of $N_{\geq 2}(w)$.
Let $Y_{1}=(X\cup N(X))\cap N_{2}(w)$, $Y_{2}=N(Y_{1})\cap N(w)$, $Z_{1}=N(X)\cap L(X)$, $Z_{2}=(X\cup N(Z_{1}\cup X))\cap N_{2}(w)$ and $Z_{3}=N(Z_{2})\cap N(w)$.
Then the following hold.
\begin{enumerate}
\item[{\upshape(i)}]
For $a\in N(w)-Y_{2}$, $E_{G}(X,N_{\leq 1}(a))=\emptyset $.
In particular, $E_{G}(X,N(w)-Y_{2})=\emptyset $.
\item[{\upshape(ii)}]
For $a\in N(w)-Z_{3}$, $E_{G}(X,N_{\leq 2}(a)-Z_{3})=\emptyset $.
In particular, $E_{G}(X,N(w)-Z_{3})=\emptyset $.
\end{enumerate}
\end{lem}
\proof
Since $X\subseteq N_{\geq 2}(w)$, a vertex $x\in X$ satisfies $N(x)\cap N(w)\neq \emptyset $ if and only if $x\in N_{2}(w)$.
This implies that $N(X)\cap N(w)=N(X\cap N_{2}(w))\cap N(w)$.
Since $X\cap N_{2}(w) \subseteq Y_{1}$,
\begin{align}
N(X)\cap N(w) = N(X\cap N_{2}(w))\cap N(w)\subseteq N(Y_{1})\cap N(w)=Y_{2}.\label{cond-lem-Xunion-Y2}
\end{align}
Since $Y_{1}=(X\cup N(X))\cap N_{2}(w)\subseteq (X\cup N(Z_{1}\cup X))\cap N_{2}(w)=Z_{2}$, we have
\begin{align}
Y_{2}=N(Y_{1})\cap N(w)\subseteq N(Z_{2})\cap N(w)=Z_{3}.\label{cond-lem-Xunion-YZ}
\end{align}

We first prove (i).
Suppose that there exists $a\in N(w)-Y_{2}$ such that $E_{G}(X,N_{\leq 1}(a))\neq \emptyset $.
Then there exists a vertex $y\in N(X)\cap N_{\leq 1}(a)$.
Clearly $y\neq w$.
If $y=a$, then by (\ref{cond-lem-Xunion-Y2}), $a\in N(X)\cap N(w)\subseteq Y_{2}$; if $y\in N(a)-\{w\}$, then $y\in N(X)\cap N_{2}(w)\subseteq Y_{1}$, and hence $a\in N(y)\cap N(w)\subseteq N(Y_{1})\cap N(w)=Y_{2}$.
In either case, we obtain $a\in Y_{2}$, which is a contradiction.
Therefore (i) holds.

Next we prove (ii).
Suppose that there exists $a'\in N(w)-Z_{3}$ such that $E_{G}(X,N_{\leq 2}(a')-Z_{3})\neq \emptyset $.
Then there exists a vertex $y'\in N(X)\cap (N_{\leq 2}(a')-Z_{3})$.
By (\ref{cond-lem-Xunion-YZ}), we have $a'\in N(w)-Z_{3}\subseteq N(w)-Y_{2}$.
This together with (i) implies that
\begin{align}
E_{G}(X,N_{\leq 1}(a')\cup (N(w)-Z_{3}))=\emptyset .\label{cond-lem-Xunion-adist}
\end{align}
Hence $y'\in N_{2}(a')-N(w)$.
Then there exists a vertex $a''\in N(y')\cap N(a')$.
Note that $a''\in N_{2}(w)$ and $y'\in N_{2}(w)\cup N_{3}(w)$.
We clearly have $a'\notin X$ and, since $a',a''\in N_{\leq 1}(a')$, it follows from (\ref{cond-lem-Xunion-adist}) that $a'',y'\notin X$.
Considering a $y'$-$w$ path $y'a''a'w$ of $G$, this leads to $y'\in N(X)\cap L(X)=Z_{1}$, and hence
$$
a''\in N(y')\cap N_{2}(w)\subseteq N(Z_{1}\cup X)\cap N_{2}(w)\subseteq Z_{2}.
$$
Consequently,
$$
a'\in N(a'')\cap N(w)\subseteq N(Z_{2})\cap N(w)=Z_{3},
$$
which is a contradiction.
Therefore (ii) holds.
\qed

\begin{lem}
\label{lem-Mk-nei}
Let $X$ be a subset of $N_{\geq 2}(w)$.
Then $N(X)\cap N_{2}(w)\subseteq \bigcup _{x\in X}M_{4}^{w}(x)$ and $N(X)\cap L(X)\subseteq \bigcup _{x\in X}M_{5}^{w}(x)$.
\end{lem}
\proof
Let $a\in N(X)\cap N_{2}(w)$, and take $y\in X\cap N(a)$ and $a'\in N(a)\cap N(w)$.
Then $yaa'w$ is a $(y,w;4)$-path of $G$, and hence $a\in M_{4}^{w}(y)\subseteq \bigcup _{x\in X}M_{4}^{w}(x)$.
Since $a$ is arbitrary, we get $N(X)\cap N_{2}(w)\subseteq \bigcup _{x\in X}M_{4}^{w}(x)$.

Let $b\in N(X)\cap L(X)$, and take $z\in X\cap N(b)$.
Since $b\in L(X)$, there exists a $b$-$w$ path $Q=bb'b''w$ with $V(Q)\cap X=\emptyset $.
In particular, $z\notin V(Q)$.
By the definition of $L(X)$, $b\in L(X)\subseteq N_{2}(w)\cup N_{3}(w)$.
Since $bz\in E(G)$ and $z\in N_{\geq 2}(w)$, the $\{C_{3},C_{4}\}$-freeness of $G$ assures us $N(z)\cap V(Q)=\{b\}$.
Hence $zbb'b''w$ is a $(z,w;5)$-path of $G$.
Consequently, $b\in M_{5}^{w}(z)\subseteq \bigcup _{x\in X}M_{5}^{w}(x)$.
Since $b$ is arbitrary, we get $N(X)\cap L(X)\subseteq \bigcup _{x\in X}M_{5}^{w}(x)$.
\qed

\begin{lem}
\label{lem-Mk-nei3}
Let $X$ be a subset of $N_{\geq 2}(w)$ with $|X|\geq 2$ such that $G[X]$ is connected.
Let $Y_{1}$, $Y_{2}$, $Z_{1}$, $Z_{2}$ and $Z_{3}$ be the sets obtained from $X$ as in Lemma~\ref{lem-Xunion}.
Then
\begin{enumerate}
\item[{\upshape(i)}]
$|Y_{2}|\leq |Y_{1}|\leq \sum _{x\in X}|M_{4}^{w}(x)|$,
\item[{\upshape(ii)}]
$|Z_{1}|\leq \sum _{x\in X}|M_{5}^{w}(x)|$, and
\item[{\upshape(iii)}]
$|Z_{3}|\leq |Z_{2}|\leq \sum _{x\in Z_{1}\cup X}|M_{4}^{w}(x)|$.
\end{enumerate}
\end{lem}
\proof
Since $|X|\geq 2$ and $G[X]$ is connected, we have $X\subseteq N(X)$, and hence $Y_{1}=N(X)\cap N_{2}(w)$ and $Z_{2}=N(Z_{1}\cup X)\cap N_{2}(w)$.
By Lemma~\ref{lem-Mk-nei}, $|Y_{1}|=|N(X)\cap N_{2}(w)|\leq |\bigcup _{x\in X}M_{4}^{w}(x)|\leq \sum _{x\in X}|M_{4}^{w}(x)|$.
Since $|N(y)\cap N(w)|=1$ for each $y\in Y_{1}~(\subseteq N_{2}(w))$, we also have $|Y_{2}|=|N(Y_{1})\cap N(w)|\leq \sum _{y\in Y_{1}}|N(y)\cap N(w)|=|Y_{1}|$, which proves (i).

By Lemma~\ref{lem-Mk-nei}, $|Z_{1}|=|N(X)\cap L(X)|\leq |\bigcup _{x\in X}M_{5}^{w}(x)|\leq \sum _{x\in X}|M_{5}^{w}(x)|$, which proves (ii).

Since $Z_{1}=N(X)\cap L(X)\subseteq N(X)\cap (N_{2}(w)\cup N_{3}(w))$, we have $Z_{1}\cup X\subseteq N_{\geq 2}(w)$.
Applying Lemma~\ref{lem-Mk-nei} to $Z_{1}\cup X$, we get $|Z_{2}|=|N(Z_{1}\cup X)\cap N_{2}(w)|\leq |\bigcup _{x\in Z_{1}\cup X}M_{4}^{w}(x)|\leq \sum _{x\in Z_{1}\cup X}|M_{4}^{w}(x)|$.
Since $|N(z)\cap N(w)|=1$ for each $z\in Z_{2}~(\subseteq N_{2}(w))$, we also have $|Z_{3}|=|N(Z_{2})\cap N(w)|\leq \sum _{z\in Z_{2}}|N(z)\cap N(w)|=|Z_{2}|$, which proves (iii).
\qed

Now we prove Theorem~\ref{prop-maxdegbound}.
First we prove Theorem~\ref{prop-maxdegbound}(i).

\begin{prop}
\label{prop-M81}
If $\Delta (G)\geq 943218$, then $T_{8}^{(1)}\prec G$.
\end{prop}
\proof
Note that $d_{G}(w)\geq 943218$.
We start with three claims.

\begin{figure}
\begin{center}
{\unitlength 0.1in%
\begin{picture}(24.5000,11.6500)(2.0000,-13.6500)%
\put(14.0000,-14.3000){\makebox(0,0){$w$}}%
%
\special{sh 1.000}%
\special{ia 1400 1300 50 50 0.0000000 6.2831853}%
\special{pn 8}%
\special{ar 1400 1300 50 50 0.0000000 6.2831853}%
%
\special{pn 8}%
\special{pa 200 900}%
\special{pa 2600 900}%
\special{pa 2600 1100}%
\special{pa 200 1100}%
\special{pa 200 900}%
\special{pa 2600 900}%
\special{fp}%
%
\special{pn 8}%
\special{pa 200 200}%
\special{pa 2600 200}%
\special{pa 2600 800}%
\special{pa 200 800}%
\special{pa 200 200}%
\special{pa 2600 200}%
\special{fp}%
%
\special{sh 1.000}%
\special{ia 800 1000 50 50 0.0000000 6.2831853}%
\special{pn 8}%
\special{ar 800 1000 50 50 0.0000000 6.2831853}%
%
\special{sh 1.000}%
\special{ia 800 700 50 50 0.0000000 6.2831853}%
\special{pn 8}%
\special{ar 800 700 50 50 0.0000000 6.2831853}%
%
\special{sh 1.000}%
\special{ia 600 600 50 50 0.0000000 6.2831853}%
\special{pn 8}%
\special{ar 600 600 50 50 0.0000000 6.2831853}%
%
\special{sh 1.000}%
\special{ia 400 400 50 50 0.0000000 6.2831853}%
\special{pn 8}%
\special{ar 400 400 50 50 0.0000000 6.2831853}%
%
\special{sh 1.000}%
\special{ia 1000 600 50 50 0.0000000 6.2831853}%
\special{pn 8}%
\special{ar 1000 600 50 50 0.0000000 6.2831853}%
%
\special{sh 1.000}%
\special{ia 1200 400 50 50 0.0000000 6.2831853}%
\special{pn 8}%
\special{ar 1200 400 50 50 0.0000000 6.2831853}%
%
\special{sh 1.000}%
\special{ia 1800 1000 50 50 0.0000000 6.2831853}%
\special{pn 8}%
\special{ar 1800 1000 50 50 0.0000000 6.2831853}%
%
\special{sh 1.000}%
\special{ia 1600 700 50 50 0.0000000 6.2831853}%
\special{pn 8}%
\special{ar 1600 700 50 50 0.0000000 6.2831853}%
%
\special{sh 1.000}%
\special{ia 2000 700 50 50 0.0000000 6.2831853}%
\special{pn 8}%
\special{ar 2000 700 50 50 0.0000000 6.2831853}%
%
\special{sh 1.000}%
\special{ia 1600 400 50 50 0.0000000 6.2831853}%
\special{pn 8}%
\special{ar 1600 400 50 50 0.0000000 6.2831853}%
%
\special{sh 1.000}%
\special{ia 2200 1000 50 50 0.0000000 6.2831853}%
\special{pn 8}%
\special{ar 2200 1000 50 50 0.0000000 6.2831853}%
\put(6.5000,-10.0000){\makebox(0,0){$a_{1}$}}%
\put(8.0000,-5.8000){\makebox(0,0){$a_{2}$}}%
\put(4.8000,-6.7000){\makebox(0,0){$a_{3}$}}%
\put(5.5000,-3.4000){\makebox(0,0){$a_{4}$}}%
\put(10.5000,-3.4000){\makebox(0,0){$a'_{4}$}}%
\put(11.5000,-6.7000){\makebox(0,0){$a'_{3}$}}%
\put(17.5000,-4.0000){\makebox(0,0){$b_{3}$}}%
\put(17.5000,-7.0000){\makebox(0,0){$b_{2}$}}%
\put(21.5000,-7.0000){\makebox(0,0){$b'_{2}$}}%
\put(16.3000,-9.9000){\makebox(0,0){$b_{1}$}}%
\put(23.6000,-10.0000){\makebox(0,0){$c_{1}$}}%
\put(26.5000,-11.0000){\makebox(0,0)[lb]{$N(w)$}}%
\put(26.5000,-8.0000){\makebox(0,0)[lb]{$N_{\geq 2}(w)$}}%
%
\special{pn 20}%
\special{pa 1400 1300}%
\special{pa 800 1000}%
\special{fp}%
\special{pa 800 1000}%
\special{pa 800 700}%
\special{fp}%
\special{pa 800 700}%
\special{pa 600 600}%
\special{fp}%
\special{pa 600 600}%
\special{pa 400 400}%
\special{fp}%
\special{pa 1200 400}%
\special{pa 1000 600}%
\special{fp}%
\special{pa 1000 600}%
\special{pa 800 700}%
\special{fp}%
\special{pa 1400 1300}%
\special{pa 1800 1000}%
\special{fp}%
\special{pa 1800 1000}%
\special{pa 1600 700}%
\special{fp}%
\special{pa 1600 700}%
\special{pa 1600 400}%
\special{fp}%
\special{pa 1800 1000}%
\special{pa 2000 700}%
\special{fp}%
\special{pa 1400 1300}%
\special{pa 2200 1000}%
\special{fp}%
\end{picture}}%

\caption{An induced copy of $T_{8}^{(1)}$ in Claim~\ref{cl-prop-M81-1}}
\label{fcl5.1}
\end{center}
\end{figure}
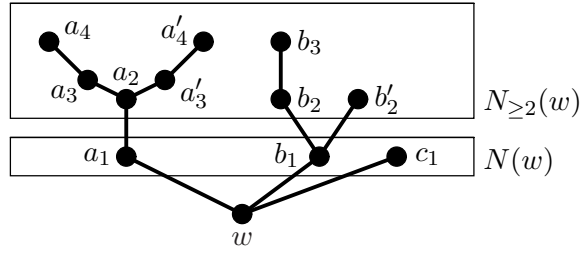

\begin{claim}
\label{cl-prop-M81-1}
If $G[N_{\geq 2}(w)]$ has an induced path $a_{4}a_{3}a_{2}a'_{3}a'_{4}$ with $a_{2}\in N_{2}(w)$, then $T_{8}^{(1)}\prec G$.
\end{claim}
\proof
Take $a_{1}\in N(a_{2})\cap N(w)$.
By Lemmas~\ref{lem-Mk-1-cor}(i) and \ref{lem-Mk-2}(i), we may assume that
\begin{align}
|M_{4}^{w}(x)|\leq 39\mbox{ and }|M_{5}^{w}(x)|\leq 4836 \mbox{ for all }x\in N_{\geq 2}(w).\label{cond-cl-prop-M81-1-1}
\end{align}
Set $X=\{a_{2},a_{3},a_{4},a'_{3},a'_{4}\}$, and let $Z_{1}$, $Z_{2}$ and $Z_{3}$ be the sets obtained from $X$ as in Lemma~\ref{lem-Xunion}.
Then by Lemma~\ref{lem-Mk-nei3}(ii)(iii) and (\ref{cond-cl-prop-M81-1-1}),
\begin{align*}
|Z_{1}| & \leq \sum _{x\in X}|M_{5}^{w}(x)|\leq 5\cdot 4836=24180,\mbox{ and}\\
|Z_{3}| &\leq |Z_{2}| \leq \sum _{x\in Z_{1}\cup X}|M_{4}^{w}(x)|\leq (24180+5)\cdot 39=943215.
\end{align*}
Since $|N(w)-Z_{3}|=d_{G}(w)-|Z_{3}|\geq 943218-943215>0$, there exists a vertex $b_{1}\in N(w)-Z_{3}$.
Take $b_{2},b'_{2}\in N(b_{1})-\{w\}$ with $b_{2}\neq b'_{2}$.
Since $|N(a_{1})\cap N(b_{2})|\leq 1$, there exists a vertex $b_{3}\in N(b_{2})-(\{b_{1}\}\cup N(a_{1}))$.
Note that $b_{3}\in N_{\geq 2}(w)$.
Since $|N(w)\cap N(b_{3})|\leq 1$, $|N(w)-(Z_{3}\cup \{b_{1}\}\cup N(b_{3}))|\geq d_{G}(w)-|Z_{3}|-2>0$.
Let $c_{1}\in N(w)-(Z_{3}\cup \{b_{1}\}\cup N(b_{3}))$.
Since $a_{1}b_{3}\notin E(G)$ and $G$ is $\{C_{3},C_{4}\}$-free, 
$$
E(G[\{w,a_{1},b_{1},b_{2},b_{3},b'_{2},c_{1}\}])=\{wa_{1},wb_{1},b_{1}b_{2},b_{2}b_{3},b_{1}b'_{2},wc_{1}\}.
$$
Since $\{b_{1},b_{2},b_{3},b'_{2}\}\subseteq N_{\leq 2}(b_{1})-Z_{3}$ and $c_{1}\in N(w)-Z_{3}$, it follows from Lemma~\ref{lem-Xunion}(ii) that $E_{G}(X,\{b_{1},b_{2},b_{3},b'_{2},c_{1}\})=\emptyset $.
Consequently, $\{a_{4},a_{3},a_{2},a'_{4},a'_{3},a_{1},w,c_{1},b_{1},b'_{2},b_{2},b_{3}\}$ induces a copy of $T_{8}^{(1)}$ in $G$ (see Figure~\ref{fcl5.1}).
\qed

\begin{claim}
\label{cl-prop-M81-2}
If $N_{\geq 3}(w)\neq \emptyset $ or $G[N_{2}(w)]$ is not $2$-regular, then $T_{8}^{(1)}\prec G$.
\end{claim}
\proof
Suppose that $N_{\geq 3}(w)\neq \emptyset $ or $G[N_{2}(w)]$ is not $2$-regular.
Since $|N(u)\cap N(w)|=1$ for $u\in N_{2}(w)$ and $N(u')\cap N(w)=\emptyset $ for $u'\in N_{3}(w)$, there exists $a_{3}\in N_{\geq 2}(w)$ such that $N(a_{3})\cap N_{2}(w)\neq \emptyset $ and $|N(a_{3})\cap N_{\geq 2}(w)|\geq 3$.
Take $a_{2}\in N(a_{3})\cap N_{2}(w)$.
Since $|N(a_{2})\cap N(w)|=1$, there exists a vertex $a'_{3}\in N(a_{2})-(\{a_{3}\}\cup N(w))$.
Since $|N(a'_{3})\cap N(w)|\leq 1$, there exists a vertex $a'_{4}\in N(a'_{3})-(\{a_{2}\}\cup N(w))$.
Since $G$ is $\{C_{3},C_{4}\}$-free, $a'_{4}\neq a_{3}$ and $a_{3}a_{2}a'_{3}a'_{4}$ is an induced path of $G$.
Since $|(N(a_{3})\cap N_{\geq 2}(w))-\{a_{2}\}|\geq 2$ and $|N(a_{3})\cap N(a_{4})|\leq 1$, we can choose $a_{4}\in (N(a_{3})\cap N_{\geq 2}(w))-\{a_{2}\}$ so that $a_{4}a'_{4}\notin E(G)$.
Then $a_{4}\neq a'_{4}$ and the path $a_{4}a_{3}a_{2}a'_{3}a'_{4}$ satisfies the assumption of Claim~\ref{cl-prop-M81-1}, and hence $T_{8}^{(1)}\prec G$.
\qed

\begin{claim}
\label{cl-prop-M81-3}
If there exists a component $C$ of $G[N_{2}(w)]$ with $C\not\simeq C_{5}$, then $T_{8}^{(1)}\prec G$.
\end{claim}
\proof
Suppose that there exists a component $C$ of $G[N_{2}(w)]$ with $C\not\simeq C_{5}$.
By Claim~\ref{cl-prop-M81-2}, we may assume that $C$ is a cycle.
Since $G$ is $\{C_{3},C_{4}\}$-free, $C$ is a cycle of order at least $6$.
Then $C$ contains an induced path of order $5$, and hence $T_{8}^{(1)}\prec G$ by Claim~\ref{cl-prop-M81-1}.
\qed

By Claims~\ref{cl-prop-M81-2} and \ref{cl-prop-M81-3}, we may assume that $N_{\geq 3}(w)=\emptyset $ and every component of $G[N_{2}(w)]$ is a cycle of order $5$.
Take $a_{1}\in N(w)$ and $a_{2},a'_{2}\in N(a_{1})-\{w\}$ with $a_{2}\neq a'_{2}$.
Note that $a_{2},a'_{2}\in N_{2}(w)$.
Let $C$ and $C'$ be the components of $G[N_{2}(w)]$ containing $a_{2}$ and $a'_{2}$, respectively.
Since $G$ is $\{C_{3},C_{4}\}$-free, we see that $C\neq C'$, $N(a_{1})\cap V(C)=\{a_{2}\}$ and $N(a_{1})\cap V(C')=\{a'_{2}\}$.
Since $C$ is a cycle of order $5$, $C$ contains an induced path $a_{4}a_{3}a_{2}a'_{3}$.

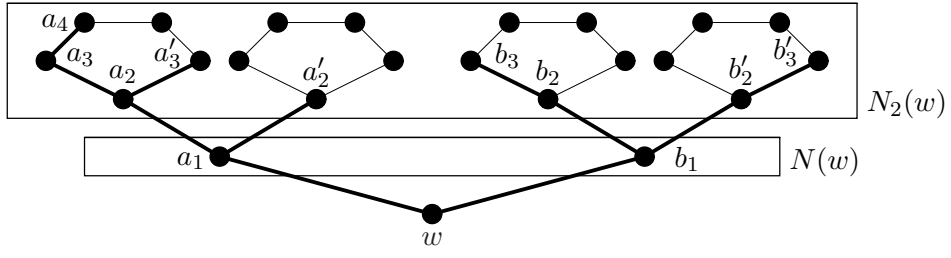
\begin{figure}
\begin{center}
{\unitlength 0.1in%
\begin{picture}(44.5000,11.6500)(10.0000,-22.6500)%
\put(32.0000,-23.3000){\makebox(0,0){$w$}}%
%
\special{sh 1.000}%
\special{ia 3200 2200 50 50 0.0000000 6.2831853}%
\special{pn 8}%
\special{ar 3200 2200 50 50 0.0000000 6.2831853}%
%
\special{pn 8}%
\special{pa 1400 1800}%
\special{pa 5000 1800}%
\special{pa 5000 2000}%
\special{pa 1400 2000}%
\special{pa 1400 1800}%
\special{pa 5000 1800}%
\special{fp}%
%
\special{pn 8}%
\special{pa 1000 1100}%
\special{pa 5400 1100}%
\special{pa 5400 1700}%
\special{pa 1000 1700}%
\special{pa 1000 1100}%
\special{pa 5400 1100}%
\special{fp}%
\put(50.5000,-20.0000){\makebox(0,0)[lb]{$N(w)$}}%
\put(54.5000,-17.0000){\makebox(0,0)[lb]{$N_{2}(w)$}}%
%
\special{sh 1.000}%
\special{ia 1400 1200 50 50 0.0000000 6.2831853}%
\special{pn 8}%
\special{ar 1400 1200 50 50 0.0000000 6.2831853}%
%
\special{sh 1.000}%
\special{ia 1800 1200 50 50 0.0000000 6.2831853}%
\special{pn 8}%
\special{ar 1800 1200 50 50 0.0000000 6.2831853}%
%
\special{sh 1.000}%
\special{ia 1200 1400 50 50 0.0000000 6.2831853}%
\special{pn 8}%
\special{ar 1200 1400 50 50 0.0000000 6.2831853}%
%
\special{sh 1.000}%
\special{ia 1600 1600 50 50 0.0000000 6.2831853}%
\special{pn 8}%
\special{ar 1600 1600 50 50 0.0000000 6.2831853}%
%
\special{sh 1.000}%
\special{ia 2000 1400 50 50 0.0000000 6.2831853}%
\special{pn 8}%
\special{ar 2000 1400 50 50 0.0000000 6.2831853}%
%
\special{pn 4}%
\special{pa 2000 1400}%
\special{pa 1600 1600}%
\special{fp}%
\special{pa 1600 1600}%
\special{pa 1200 1400}%
\special{fp}%
\special{pa 1200 1400}%
\special{pa 1400 1200}%
\special{fp}%
\special{pa 1400 1200}%
\special{pa 1800 1200}%
\special{fp}%
\special{pa 1800 1200}%
\special{pa 2000 1400}%
\special{fp}%
%
\special{sh 1.000}%
\special{ia 2400 1200 50 50 0.0000000 6.2831853}%
\special{pn 8}%
\special{ar 2400 1200 50 50 0.0000000 6.2831853}%
%
\special{sh 1.000}%
\special{ia 2800 1200 50 50 0.0000000 6.2831853}%
\special{pn 8}%
\special{ar 2800 1200 50 50 0.0000000 6.2831853}%
%
\special{sh 1.000}%
\special{ia 2200 1400 50 50 0.0000000 6.2831853}%
\special{pn 8}%
\special{ar 2200 1400 50 50 0.0000000 6.2831853}%
%
\special{sh 1.000}%
\special{ia 2600 1600 50 50 0.0000000 6.2831853}%
\special{pn 8}%
\special{ar 2600 1600 50 50 0.0000000 6.2831853}%
%
\special{sh 1.000}%
\special{ia 3000 1400 50 50 0.0000000 6.2831853}%
\special{pn 8}%
\special{ar 3000 1400 50 50 0.0000000 6.2831853}%
%
\special{pn 4}%
\special{pa 3000 1400}%
\special{pa 2600 1600}%
\special{fp}%
\special{pa 2600 1600}%
\special{pa 2200 1400}%
\special{fp}%
\special{pa 2200 1400}%
\special{pa 2400 1200}%
\special{fp}%
\special{pa 2400 1200}%
\special{pa 2800 1200}%
\special{fp}%
\special{pa 2800 1200}%
\special{pa 3000 1400}%
\special{fp}%
%
\special{sh 1.000}%
\special{ia 3600 1200 50 50 0.0000000 6.2831853}%
\special{pn 8}%
\special{ar 3600 1200 50 50 0.0000000 6.2831853}%
%
\special{sh 1.000}%
\special{ia 4000 1200 50 50 0.0000000 6.2831853}%
\special{pn 8}%
\special{ar 4000 1200 50 50 0.0000000 6.2831853}%
%
\special{sh 1.000}%
\special{ia 3400 1400 50 50 0.0000000 6.2831853}%
\special{pn 8}%
\special{ar 3400 1400 50 50 0.0000000 6.2831853}%
%
\special{sh 1.000}%
\special{ia 3800 1600 50 50 0.0000000 6.2831853}%
\special{pn 8}%
\special{ar 3800 1600 50 50 0.0000000 6.2831853}%
%
\special{sh 1.000}%
\special{ia 4200 1400 50 50 0.0000000 6.2831853}%
\special{pn 8}%
\special{ar 4200 1400 50 50 0.0000000 6.2831853}%
%
\special{pn 4}%
\special{pa 4200 1400}%
\special{pa 3800 1600}%
\special{fp}%
\special{pa 3800 1600}%
\special{pa 3400 1400}%
\special{fp}%
\special{pa 3400 1400}%
\special{pa 3600 1200}%
\special{fp}%
\special{pa 3600 1200}%
\special{pa 4000 1200}%
\special{fp}%
\special{pa 4000 1200}%
\special{pa 4200 1400}%
\special{fp}%
%
\special{sh 1.000}%
\special{ia 4600 1200 50 50 0.0000000 6.2831853}%
\special{pn 8}%
\special{ar 4600 1200 50 50 0.0000000 6.2831853}%
%
\special{sh 1.000}%
\special{ia 5000 1200 50 50 0.0000000 6.2831853}%
\special{pn 8}%
\special{ar 5000 1200 50 50 0.0000000 6.2831853}%
%
\special{sh 1.000}%
\special{ia 4400 1400 50 50 0.0000000 6.2831853}%
\special{pn 8}%
\special{ar 4400 1400 50 50 0.0000000 6.2831853}%
%
\special{sh 1.000}%
\special{ia 4800 1600 50 50 0.0000000 6.2831853}%
\special{pn 8}%
\special{ar 4800 1600 50 50 0.0000000 6.2831853}%
%
\special{sh 1.000}%
\special{ia 5200 1400 50 50 0.0000000 6.2831853}%
\special{pn 8}%
\special{ar 5200 1400 50 50 0.0000000 6.2831853}%
%
\special{pn 4}%
\special{pa 5200 1400}%
\special{pa 4800 1600}%
\special{fp}%
\special{pa 4800 1600}%
\special{pa 4400 1400}%
\special{fp}%
\special{pa 4400 1400}%
\special{pa 4600 1200}%
\special{fp}%
\special{pa 4600 1200}%
\special{pa 5000 1200}%
\special{fp}%
\special{pa 5000 1200}%
\special{pa 5200 1400}%
\special{fp}%
%
\special{sh 1.000}%
\special{ia 2100 1900 50 50 0.0000000 6.2831853}%
\special{pn 8}%
\special{ar 2100 1900 50 50 0.0000000 6.2831853}%
\put(19.5000,-19.2000){\makebox(0,0){$a_{1}$}}%
%
\special{sh 1.000}%
\special{ia 4300 1900 50 50 0.0000000 6.2831853}%
\special{pn 8}%
\special{ar 4300 1900 50 50 0.0000000 6.2831853}%
\put(45.2000,-19.1000){\makebox(0,0){$b_{1}$}}%
%
\special{pn 20}%
\special{pa 2100 1900}%
\special{pa 1600 1600}%
\special{fp}%
\special{pa 2600 1600}%
\special{pa 2100 1900}%
\special{fp}%
%
\special{pn 20}%
\special{pa 4300 1900}%
\special{pa 3800 1600}%
\special{fp}%
\special{pa 4800 1600}%
\special{pa 4300 1900}%
\special{fp}%
%
\special{pn 20}%
\special{pa 3200 2200}%
\special{pa 2100 1900}%
\special{fp}%
\special{pa 4300 1900}%
\special{pa 3200 2200}%
\special{fp}%
\put(16.0000,-14.7000){\makebox(0,0){$a_{2}$}}%
\put(13.8000,-13.8000){\makebox(0,0){$a_{3}$}}%
\put(12.5000,-12.0000){\makebox(0,0){$a_{4}$}}%
\put(18.3000,-13.6000){\makebox(0,0){$a'_{3}$}}%
\put(26.0000,-14.6000){\makebox(0,0){$a'_{2}$}}%
\put(38.0000,-14.6000){\makebox(0,0){$b_{2}$}}%
\put(48.0000,-14.5000){\makebox(0,0){$b'_{2}$}}%
\put(50.3000,-13.5000){\makebox(0,0){$b'_{3}$}}%
\put(35.9000,-13.7000){\makebox(0,0){$b_{3}$}}%
%
\special{pn 20}%
\special{pa 1600 1600}%
\special{pa 1200 1400}%
\special{fp}%
\special{pa 1200 1400}%
\special{pa 1400 1200}%
\special{fp}%
\special{pa 1600 1600}%
\special{pa 2000 1400}%
\special{fp}%
\special{pa 3400 1400}%
\special{pa 3800 1600}%
\special{fp}%
\special{pa 4800 1600}%
\special{pa 5200 1400}%
\special{fp}%
\end{picture}}%

\caption{An induced copy of $T_{8}^{(1)}$ in Proposition~\ref{prop-M81}}
\label{fprop5.4}
\end{center}
\end{figure}

Since $|N(V(C)\cup V(C'))\cap N(w)|\leq \sum _{x\in V(C)\cup V(C')}|N(x)\cap N(w)|=10$, $|N(w)-N(V(C)\cup V(C'))|\geq 943218-10>0$.
Take $b_{1}\in N(w)-N(V(C)\cup V(C'))$ and $b_{2},b'_{2}\in N(b_{1})-\{w\}$ with $b_{2}\neq b'_{2}$.
Note that $b_{2},b'_{2}\in N_{2}(w)$.
Let $D$ and $D'$ be the components of $G[N_{2}(w)]$ containing $b_{2}$ and $b'_{2}$, respectively.
Then $C,C',D,D'$ are pairwise distinct and there is no edge between any two of $C,C',D,D'$.
Since $G$ is $C_{4}$-free, $a_{1}b_{2},a_{1}b'_{2}\notin E(G)$.
Furthermore, since $|N(b_{2})\cap N(a_{1})|\leq 1$ and $|N(b'_{2})\cap N(a_{1})|\leq 1$, there exist vertices $b_{3}\in N_{D}(b_{2})-N(a_{1})$ and $b'_{3}\in N_{D'}(b'_{2})-N(a_{1})$.
Then $\{b_{3},b_{2},b_{1},b'_{3},b'_{2},w,a_{1},a'_{2},a_{2},a'_{3},a_{3},a_{4}\}$ induces a copy of $T_{8}^{(1)}$ in $G$ (see Figure~\ref{fprop5.4}).
\qed

Next we prove Theorem~\ref{prop-maxdegbound}(ii).

\begin{prop}
\label{prop-M82}
If $\Delta (G)\geq 190375$, then $T_{8}^{(2)}\prec G$.
\end{prop}
\proof
Note that $d_{G}(w)\geq 190375$.
By Lemmas~\ref{lem-Mk-1-cor}(ii) and \ref{lem-Mk-2}(ii), we may assume that
\begin{align}
|M_{4}^{w}(x)|\leq 31\mbox{ and }|M_{5}^{w}(x)|\leq 2046 \mbox{ for all }x\in N_{\geq 2}(w).\label{cond-prop-M82-1}
\end{align}
Take a vertex $a_{1}\in N(w)$ which attains the minimum of $|N(x)-\{w\}|~(=d_{G}(x)-1)$ as $x$ ranges over $N(w)$, and take $a_{2}\in N(a_{1})-\{w\}$ and $a_{3}\in N(a_{2})-\{a_{1}\}$.
Note that $a_{2}\in N_{2}(w)$.
Since $N(a_{2})\cap N(w)=\{a_{1}\}$, we have $a_{3}\in N_{\geq 2}(w)$.
Since $|N(a_{3})\cap N(w)|\leq 1$, there exists a vertex $a_{4}\in N(a_{3})-(\{a_{2}\}\cup N(w))$.
Since $G$ is $\{C_{3},C_{4}\}$-free, $a_{1}a_{2}a_{3}a_{4}$ is an induced path of $G[N_{\geq 2}(w)]$.

Set $X=\{a_{2},a_{3},a_{4}\}$, and let $Z_{1}$, $Z_{2}$ and $Z_{3}$ be the sets obtained from $X$ as in Lemma~\ref{lem-Xunion}.
Then by Lemma~\ref{lem-Mk-nei3}(ii)(iii) and (\ref{cond-prop-M82-1}),
\begin{align*}
|Z_{1}| & \leq \sum _{x\in X}|M_{5}^{w}(x)|\leq 3\cdot 2046=6138,\mbox{ and}\\
|Z_{3}| &\leq |Z_{2}| \leq \sum _{x\in Z_{1}\cup X}|M_{4}^{w}(x)|\leq (6138+3)\cdot 31=190371.
\end{align*}
Since $|N(w)-Z_{3}|=d_{G}(w)-|Z_{3}|\geq 190375-190371>0$, there exists a vertex $b_{1}\in N(w)-Z_{3}$.

Since $G$ is $\{C_{3},C_{4}\}$-free, $N(x)\cap N(x')\cap N_{2}(w)=\emptyset $ for any $x,x'\in N(b_{1})-\{w\}$ with $x\neq x'$, and hence
\begin{align}
\sum _{x\in N(b_{1})-\{w\}}|N(x)\cap N(a_{1})\cap N_{2}(w)| = |N(N(b_{1})-\{w\})\cap N(a_{1})\cap N_{2}(w)|.\label{cond-propM82-new1}
\end{align}
Since $N(b_{1})-\{w\}\subseteq N_{\leq 2}(b_{1})-Z_{3}$, it follows from Lemma~\ref{lem-Xunion}(ii) that $N(X)\cap (N(b_{1})-\{w\})=\emptyset $.
In particular, $N(a_{2})\cap (N(b_{1})-\{w\})=\emptyset $, that is,
\begin{align}
N(N(b_{1})-\{w\})\subseteq V(G)-\{a_{2}\}.\label{cond-propM82-new2}
\end{align}
By (\ref{cond-propM82-new1}), (\ref{cond-propM82-new2}) and the minimality of $|N(a_{1})-\{w\}|$, we have
\begin{align*}
\sum _{x\in N(b_{1})-\{w\}}|N(x)\cap N(a_{1})\cap N_{2}(w)| &= |N(N(b_{1})-\{w\})\cap N(a_{1})\cap N_{2}(w)|\\
&\leq |(N(a_{1})\cap N_{2}(w))-\{a_{2}\}|\\
&= |N(a_{1})-\{w\}|-1\\
&\leq |N(b_{1})-\{w\}|-1.
\end{align*}
This inequality assures us the existence of $b_{2}\in N(b_{1})-\{w\}$ such that $N(b_{2})\cap N(a_{1})\cap N_{2}(w)=\emptyset $.
Take $b'_{2}\in N(b_{1})-\{w,b_{2}\}$.
Note that $b'_{2}\in N_{2}(w)$.
Since $|N(b'_{2})\cap N(a_{1})|\leq 1$, there exists a vertex $b'_{3}\in N(b'_{2})-(\{b_{1}\}\cup N(a_{1}))$.
Since $|N(b_{2})\cap N(b'_{3})|\leq 1$, there exists a vertex $b_{3}\in N(b_{2})-(\{b_{1}\}\cup N(b'_{3}))$.
Note that $b_{3}\in N_{2}(w)\cup N_{3}(w)$.
Since $N(b_{2})\cap N(a_{1})\cap N_{2}(w)=\emptyset $ and $a_{1}\in N(w)$, it follows that $b_{3}a_{1}\notin E(G)$.
Consequently, $N(a_{1})\cap \{b_{1},b_{2},b_{3},b'_{2},b'_{3}\}=\emptyset $.
Since $|N(b_{3})\cap N(w)|\leq 1$ and $|N(b'_{3})\cap N(w)|\leq 1$, we have $|N(w)-(Z_{3}\cup \{b_{1}\}\cup N(\{b_{3},b'_{3}\}))|\geq 190375-(190371+1+2)>0$.
Take $c_{1}\in N(w)-(Z_{3}\cup \{b_{1}\}\cup N(\{b_{3},b'_{3}\}))$.
Then
$$
E(G[\{w,a_{1},b_{1},b_{2},b_{3},b'_{2},b'_{3},c_{1}\}])=\{wa_{1},wb_{1},b_{1}b_{2},b_{2}b_{3},b_{1}b'_{2},b'_{2}b'_{3},wc_{1}\}.
$$
Recall that $a_{1}a_{2}a_{3}a_{4}$ is an induced path of $G$.

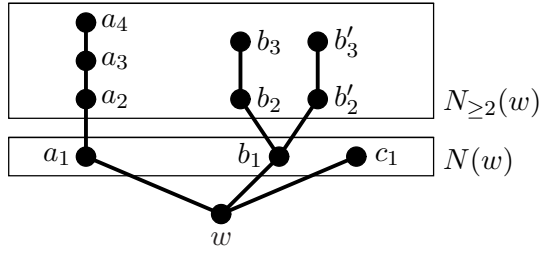
\begin{figure}
\begin{center}
{\unitlength 0.1in%
\begin{picture}(22.5000,11.6500)(4.0000,-13.6500)%
\put(15.0000,-14.3000){\makebox(0,0){$w$}}%
%
\special{sh 1.000}%
\special{ia 1500 1300 50 50 0.0000000 6.2831853}%
\special{pn 8}%
\special{ar 1500 1300 50 50 0.0000000 6.2831853}%
%
\special{pn 8}%
\special{pa 400 900}%
\special{pa 2600 900}%
\special{pa 2600 1100}%
\special{pa 400 1100}%
\special{pa 400 900}%
\special{pa 2600 900}%
\special{fp}%
%
\special{pn 8}%
\special{pa 400 200}%
\special{pa 2600 200}%
\special{pa 2600 800}%
\special{pa 400 800}%
\special{pa 400 200}%
\special{pa 2600 200}%
\special{fp}%
%
\special{sh 1.000}%
\special{ia 800 1000 50 50 0.0000000 6.2831853}%
\special{pn 8}%
\special{ar 800 1000 50 50 0.0000000 6.2831853}%
%
\special{sh 1.000}%
\special{ia 800 700 50 50 0.0000000 6.2831853}%
\special{pn 8}%
\special{ar 800 700 50 50 0.0000000 6.2831853}%
%
\special{sh 1.000}%
\special{ia 1800 1000 50 50 0.0000000 6.2831853}%
\special{pn 8}%
\special{ar 1800 1000 50 50 0.0000000 6.2831853}%
%
\special{sh 1.000}%
\special{ia 1600 700 50 50 0.0000000 6.2831853}%
\special{pn 8}%
\special{ar 1600 700 50 50 0.0000000 6.2831853}%
%
\special{sh 1.000}%
\special{ia 2000 700 50 50 0.0000000 6.2831853}%
\special{pn 8}%
\special{ar 2000 700 50 50 0.0000000 6.2831853}%
%
\special{sh 1.000}%
\special{ia 1600 400 50 50 0.0000000 6.2831853}%
\special{pn 8}%
\special{ar 1600 400 50 50 0.0000000 6.2831853}%
%
\special{sh 1.000}%
\special{ia 2200 1000 50 50 0.0000000 6.2831853}%
\special{pn 8}%
\special{ar 2200 1000 50 50 0.0000000 6.2831853}%
\put(6.5000,-10.0000){\makebox(0,0){$a_{1}$}}%
\put(17.5000,-4.0000){\makebox(0,0){$b_{3}$}}%
\put(17.5000,-7.0000){\makebox(0,0){$b_{2}$}}%
\put(21.5000,-7.0000){\makebox(0,0){$b'_{2}$}}%
\put(16.5000,-9.9000){\makebox(0,0){$b_{1}$}}%
\put(23.6000,-10.0000){\makebox(0,0){$c_{1}$}}%
\put(26.5000,-11.0000){\makebox(0,0)[lb]{$N(w)$}}%
\put(26.5000,-8.0000){\makebox(0,0)[lb]{$N_{\geq 2}(w)$}}%
%
\special{sh 1.000}%
\special{ia 800 500 50 50 0.0000000 6.2831853}%
\special{pn 8}%
\special{ar 800 500 50 50 0.0000000 6.2831853}%
%
\special{sh 1.000}%
\special{ia 800 300 50 50 0.0000000 6.2831853}%
\special{pn 8}%
\special{ar 800 300 50 50 0.0000000 6.2831853}%
\put(21.5000,-4.0000){\makebox(0,0){$b'_{3}$}}%
%
\special{sh 1.000}%
\special{ia 2000 400 50 50 0.0000000 6.2831853}%
\special{pn 8}%
\special{ar 2000 400 50 50 0.0000000 6.2831853}%
\put(9.5000,-7.0000){\makebox(0,0){$a_{2}$}}%
\put(9.5000,-5.0000){\makebox(0,0){$a_{3}$}}%
\put(9.5000,-3.0000){\makebox(0,0){$a_{4}$}}%
%
\special{pn 20}%
\special{pa 1500 1300}%
\special{pa 800 1000}%
\special{fp}%
\special{pa 800 1000}%
\special{pa 800 300}%
\special{fp}%
\special{pa 1500 1300}%
\special{pa 1800 1000}%
\special{fp}%
\special{pa 1800 1000}%
\special{pa 1600 700}%
\special{fp}%
\special{pa 1600 700}%
\special{pa 1600 400}%
\special{fp}%
\special{pa 2000 400}%
\special{pa 2000 700}%
\special{fp}%
\special{pa 2000 700}%
\special{pa 1800 1000}%
\special{fp}%
\special{pa 1500 1300}%
\special{pa 2200 1000}%
\special{fp}%
\end{picture}}%

\caption{An induced copy of $T_{8}^{(2)}$ in Proposition~\ref{prop-M82}}
\label{fprop5.5}
\end{center}
\end{figure}

Since $\{b_{1},b_{2},b_{3},b'_{2},b'_{3}\}\subseteq N_{\leq 2}(b_{1})-Z_{3}$ and $c_{1}\in N(w)-Z_{3}$, it follows from Lemma~\ref{lem-Xunion}(ii) that $E_{G}(X,\{b_{1},b_{2},b_{3},b'_{2},b'_{3},c_{1}\})=\emptyset $.
Consequently, $\{b_{3},b_{2},b_{1},b'_{3},b'_{2},w,c_{1},a_{1},a_{2},a_{3},a_{4}\}$ induces a copy of $T_{8}^{(2)}$ in $G$ (see Figure~\ref{fprop5.5}).
\qed

Finally we prove Theorem~\ref{prop-maxdegbound}(iii).

\begin{prop}
\label{prop-M91}
If $\Delta (G)\geq 197433$, then $T_{9}\prec G$. 
\end{prop}
\proof
Note that $d_{G}(w)\geq 197433$.
By Lemmas~\ref{lem-Mk-4} and \ref{lem-Mk-5}, we may assume that
\begin{align}
|M_{4}^{w}(x)|\leq 74\mbox{ and }|M_{5}^{w}(x)|\leq 666 \mbox{ for all }x\in N_{\geq 2}(w).\label{cond-prop-M91-1}
\end{align}
Set
$$
k_{0}=\max\{k:\mbox{there is an induced path }a_{1}a_{2}\cdots a_{k}\mbox{ with }a_{1}\in N(w)\mbox{ and }\{a_{2},\ldots ,a_{k}\}\subseteq N_{\geq 2}(w)\}.
$$

\medskip
\noindent
\textbf{Case 1:} $k_{0}\geq 6$.

\begin{figure}
\begin{center}
{\unitlength 0.1in%
\begin{picture}(22.5000,11.6500)(4.0000,-25.6500)%
\put(15.0000,-26.3000){\makebox(0,0){$w$}}%
%
\special{sh 1.000}%
\special{ia 1500 2500 50 50 0.0000000 6.2831853}%
\special{pn 8}%
\special{ar 1500 2500 50 50 0.0000000 6.2831853}%
%
\special{pn 8}%
\special{pa 400 2100}%
\special{pa 2600 2100}%
\special{pa 2600 2300}%
\special{pa 400 2300}%
\special{pa 400 2100}%
\special{pa 2600 2100}%
\special{fp}%
%
\special{pn 8}%
\special{pa 400 1400}%
\special{pa 2600 1400}%
\special{pa 2600 2000}%
\special{pa 400 2000}%
\special{pa 400 1400}%
\special{pa 2600 1400}%
\special{fp}%
%
\special{sh 1.000}%
\special{ia 800 2200 50 50 0.0000000 6.2831853}%
\special{pn 8}%
\special{ar 800 2200 50 50 0.0000000 6.2831853}%
%
\special{sh 1.000}%
\special{ia 800 1900 50 50 0.0000000 6.2831853}%
\special{pn 8}%
\special{ar 800 1900 50 50 0.0000000 6.2831853}%
%
\special{sh 1.000}%
\special{ia 1800 2200 50 50 0.0000000 6.2831853}%
\special{pn 8}%
\special{ar 1800 2200 50 50 0.0000000 6.2831853}%
%
\special{sh 1.000}%
\special{ia 1800 1700 50 50 0.0000000 6.2831853}%
\special{pn 8}%
\special{ar 1800 1700 50 50 0.0000000 6.2831853}%
%
\special{sh 1.000}%
\special{ia 2200 2200 50 50 0.0000000 6.2831853}%
\special{pn 8}%
\special{ar 2200 2200 50 50 0.0000000 6.2831853}%
\put(6.5000,-22.0000){\makebox(0,0){$a_{1}$}}%
\put(19.5000,-17.0500){\makebox(0,0){$b_{2}$}}%
\put(16.5000,-21.9000){\makebox(0,0){$b_{1}$}}%
\put(23.6000,-22.0000){\makebox(0,0){$c_{1}$}}%
\put(26.5000,-23.0000){\makebox(0,0)[lb]{$N(w)$}}%
\put(26.5000,-20.0000){\makebox(0,0)[lb]{$N_{\geq 2}(w)$}}%
%
\special{sh 1.000}%
\special{ia 800 1500 50 50 0.0000000 6.2831853}%
\special{pn 8}%
\special{ar 800 1500 50 50 0.0000000 6.2831853}%
\put(9.5000,-19.0000){\makebox(0,0){$a_{2}$}}%
\put(9.5000,-15.0000){\makebox(0,0){$a_{6}$}}%
%
\special{sh 1.000}%
\special{ia 2200 1700 50 50 0.0000000 6.2831853}%
\special{pn 8}%
\special{ar 2200 1700 50 50 0.0000000 6.2831853}%
\put(23.5000,-17.0500){\makebox(0,0){$c_{2}$}}%
%
\special{pn 4}%
\special{sh 1}%
\special{ar 800 1655 16 16 0 6.2831853}%
\special{sh 1}%
\special{ar 800 1755 16 16 0 6.2831853}%
\special{sh 1}%
\special{ar 800 1705 16 16 0 6.2831853}%
\special{sh 1}%
\special{ar 800 1705 16 16 0 6.2831853}%
%
\special{pn 20}%
\special{pa 800 1900}%
\special{pa 800 1800}%
\special{fp}%
\special{pa 800 1600}%
\special{pa 800 1500}%
\special{fp}%
%
\special{pn 20}%
\special{pa 1500 2500}%
\special{pa 800 2200}%
\special{fp}%
\special{pa 800 2200}%
\special{pa 800 1900}%
\special{fp}%
\special{pa 1500 2500}%
\special{pa 1800 2200}%
\special{fp}%
\special{pa 1800 2200}%
\special{pa 1800 1700}%
\special{fp}%
\special{pa 2200 1700}%
\special{pa 2200 2200}%
\special{fp}%
\special{pa 2200 2200}%
\special{pa 1500 2500}%
\special{fp}%
\end{picture}}%

\caption{An induced copy of $T_{9}$ in Case~1 of Proposition~\ref{prop-M91}}
\label{fprop5.6-1}
\end{center}
\end{figure}
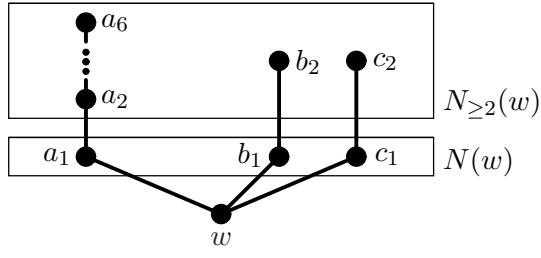

Let $a_{1}a_{2}\cdots a_{6}$ be an induced path with $a_{1}\in N(w)$ and $\{a_{2},\ldots ,a_{6}\}\subseteq N_{\geq 2}(w)$.
Set $X=\{a_{2},\ldots ,a_{6}\}$, and let $Y_{1}$ and $Y_{2}$ be the sets obtained from $X$ as in Lemma~\ref{lem-Xunion}.
Then by Lemma~\ref{lem-Mk-nei3}(i) and (\ref{cond-prop-M91-1}),
$$
|Y_{2}|\leq |Y_{1}| \leq \sum _{x\in X}|M_{4}^{w}(x)|\leq 5\cdot 74=370.
$$
Since $|N(w)-Y_{2}|=d_{G}(w)-|Y_{2}|\geq 197433-370>1$, there exist vertices $b_{1},c_{1}\in N(w)-Y_{2}$ with $b_{1}\neq c_{1}$.
Take $b_{2}\in N(b_{1})-\{w\}$.
Since $|N(c_{1})\cap N(b_{2})|\leq 1$, there exists a vertex $c_{2}\in N(c_{1})-(\{w\}\cup N(b_{2}))$.
Then $b_{2},c_{2}\in N_{2}(w)$, $b_{2}\neq c_{2}$ and $b_{2}a_{1},b_{2}c_{1},c_{2}a_{1},c_{2}b_{1}\notin E(G)$, and hence
$$
E(G[\{w,a_{1},b_{1},b_{2},c_{1},c_{2}\}])=\{wa_{1},wb_{1},b_{1}b_{2},wc_{1},c_{1}c_{2}\}.
$$
Since $\{b_{1},b_{2}\}\subseteq N_{\leq 1}(b_{1})$ and $\{c_{1},c_{2}\}\subseteq N_{\leq 1}(c_{1})$, it follows from Lemma~\ref{lem-Xunion}(i) that $E_{G}(X,\{w,b_{1},b_{2},c_{1},c_{2}\})=\emptyset$.
Consequently, $\{b_{2},b_{1},w,c_{2},c_{1},a_{1},\ldots ,a_{6}\}$ induces a copy of $T_{9}$ in $G$ (see Figure~\ref{fprop5.6-1}).

\medskip
\noindent
\textbf{Case 2:} $k_{0}=5$.

Let $a_{1}a_{2}\cdots a_{5}$ be an induced path with $a_{1}\in N(w)$ and $\{a_{2},\ldots ,a_{5}\}\subseteq N_{\geq 2}(w)$.
Set $X=\{a_{2},\ldots ,a_{5}\}$, and let $Z_{1}$, $Z_{2}$ and $Z_{3}$ be the sets obtained from $X$ as in Lemma~\ref{lem-Xunion}.
Then by Lemma~\ref{lem-Mk-nei3}(ii)(iii) and (\ref{cond-prop-M91-1}),
\begin{align*}
|Z_{1}| & \leq \sum _{x\in X}|M_{5}^{w}(x)|\leq 4\cdot 666=2664,\mbox{ and}\\
|Z_{3}| &\leq |Z_{2}| \leq \sum _{x\in Z_{1}\cup X}|M_{4}^{w}(x)|\leq (2664+4)\cdot 74=197432.
\end{align*}
Since $|N(w)-Z_{3}|=d_{G}(w)-|Z_{3}|\geq 197433-197432>0$, there exists a vertex $b_{1}\in N(w)-Z_{3}$.
Take $b_{2},b'_{2}\in N(b_{1})-\{w\}$ with $b_{2}\neq b'_{2}$.
Then $b_{2},b'_{2}\in N_{2}(w)$.

\begin{figure}
\begin{center}
{\unitlength 0.1in%
\begin{picture}(22.5000,11.6500)(4.0000,-25.6500)%
\put(15.0000,-26.3000){\makebox(0,0){$w$}}%
%
\special{sh 1.000}%
\special{ia 1500 2500 50 50 0.0000000 6.2831853}%
\special{pn 8}%
\special{ar 1500 2500 50 50 0.0000000 6.2831853}%
%
\special{pn 8}%
\special{pa 400 2100}%
\special{pa 2600 2100}%
\special{pa 2600 2300}%
\special{pa 400 2300}%
\special{pa 400 2100}%
\special{pa 2600 2100}%
\special{fp}%
%
\special{pn 8}%
\special{pa 400 1400}%
\special{pa 2600 1400}%
\special{pa 2600 2000}%
\special{pa 400 2000}%
\special{pa 400 1400}%
\special{pa 2600 1400}%
\special{fp}%
%
\special{sh 1.000}%
\special{ia 800 2200 50 50 0.0000000 6.2831853}%
\special{pn 8}%
\special{ar 800 2200 50 50 0.0000000 6.2831853}%
%
\special{sh 1.000}%
\special{ia 800 1900 50 50 0.0000000 6.2831853}%
\special{pn 13}%
\special{ar 800 1900 50 50 0.0000000 6.2831853}%
%
\special{sh 1.000}%
\special{ia 2000 2200 50 50 0.0000000 6.2831853}%
\special{pn 8}%
\special{ar 2000 2200 50 50 0.0000000 6.2831853}%
%
\special{sh 1.000}%
\special{ia 1800 1900 50 50 0.0000000 6.2831853}%
\special{pn 8}%
\special{ar 1800 1900 50 50 0.0000000 6.2831853}%
\put(6.5000,-22.0000){\makebox(0,0){$a_{1}$}}%
\put(17.1000,-17.8000){\makebox(0,0){$b_{2}$}}%
\put(21.5000,-22.0000){\makebox(0,0){$b_{1}$}}%
\put(26.5000,-23.0000){\makebox(0,0)[lb]{$N(w)$}}%
\put(26.5000,-20.0000){\makebox(0,0)[lb]{$N_{\geq 2}(w)$}}%
%
\special{sh 1.000}%
\special{ia 800 1500 50 50 0.0000000 6.2831853}%
\special{pn 8}%
\special{ar 800 1500 50 50 0.0000000 6.2831853}%
\put(9.5000,-19.0000){\makebox(0,0){$a_{2}$}}%
\put(9.5000,-15.0000){\makebox(0,0){$a_{5}$}}%
%
\special{sh 1.000}%
\special{ia 2200 1900 50 50 0.0000000 6.2831853}%
\special{pn 8}%
\special{ar 2200 1900 50 50 0.0000000 6.2831853}%
%
\special{pn 4}%
\special{sh 1}%
\special{ar 800 1655 16 16 0 6.2831853}%
\special{sh 1}%
\special{ar 800 1755 16 16 0 6.2831853}%
\special{sh 1}%
\special{ar 800 1705 16 16 0 6.2831853}%
\special{sh 1}%
\special{ar 800 1705 16 16 0 6.2831853}%
%
\special{pn 4}%
\special{pa 800 1900}%
\special{pa 800 1800}%
\special{fp}%
\special{pa 800 1600}%
\special{pa 800 1500}%
\special{fp}%
%
\special{sh 1.000}%
\special{ia 1800 1610 50 50 0.0000000 6.2831853}%
\special{pn 8}%
\special{ar 1800 1610 50 50 0.0000000 6.2831853}%
%
\special{sh 1.000}%
\special{ia 2200 1610 50 50 0.0000000 6.2831853}%
\special{pn 8}%
\special{ar 2200 1610 50 50 0.0000000 6.2831853}%
%
\special{sh 1.000}%
\special{ia 1500 1900 50 50 0.0000000 6.2831853}%
\special{pn 8}%
\special{ar 1500 1900 50 50 0.0000000 6.2831853}%
%
\special{sh 1.000}%
\special{ia 2500 1900 50 50 0.0000000 6.2831853}%
\special{pn 8}%
\special{ar 2500 1900 50 50 0.0000000 6.2831853}%
%
\special{pn 4}%
\special{pa 2500 1900}%
\special{pa 2438 1918}%
\special{pa 2407 1926}%
\special{pa 2377 1935}%
\special{pa 2346 1944}%
\special{pa 2315 1952}%
\special{pa 2284 1961}%
\special{pa 2253 1969}%
\special{pa 2222 1978}%
\special{pa 2160 1994}%
\special{pa 2129 2003}%
\special{pa 2067 2019}%
\special{pa 2036 2026}%
\special{pa 1974 2042}%
\special{pa 1881 2063}%
\special{pa 1849 2070}%
\special{pa 1818 2077}%
\special{pa 1756 2089}%
\special{pa 1724 2096}%
\special{pa 1693 2101}%
\special{pa 1661 2107}%
\special{pa 1630 2112}%
\special{pa 1598 2117}%
\special{pa 1567 2122}%
\special{pa 1535 2127}%
\special{pa 1504 2132}%
\special{pa 1408 2144}%
\special{pa 1377 2148}%
\special{pa 1345 2152}%
\special{pa 1313 2155}%
\special{pa 1281 2159}%
\special{pa 1250 2162}%
\special{pa 1058 2180}%
\special{pa 1026 2182}%
\special{pa 962 2188}%
\special{pa 930 2190}%
\special{pa 898 2193}%
\special{pa 834 2197}%
\special{pa 802 2200}%
\special{pa 800 2200}%
\special{fp}%
%
\special{pn 20}%
\special{pa 800 2200}%
\special{pa 1500 1900}%
\special{fp}%
\special{pa 1500 1900}%
\special{pa 1800 1900}%
\special{fp}%
\special{pa 1800 1900}%
\special{pa 1800 1600}%
\special{fp}%
\special{pa 1800 1600}%
\special{pa 2200 1600}%
\special{fp}%
\special{pa 2200 1600}%
\special{pa 2200 1900}%
\special{fp}%
%
\special{pn 4}%
\special{pa 2000 2200}%
\special{pa 1800 1900}%
\special{fp}%
\special{pa 2200 1900}%
\special{pa 2000 2200}%
\special{fp}%
\put(14.3000,-17.9000){\makebox(0,0){$u_{2}$}}%
\put(25.2000,-17.7000){\makebox(0,0){$u'_{2}$}}%
\put(18.0000,-14.9000){\makebox(0,0){$b_{3}$}}%
\put(22.5000,-14.9000){\makebox(0,0){$b'_{3}$}}%
\put(21.1000,-17.7000){\makebox(0,0){$b'_{2}$}}%
%
\special{pn 4}%
\special{pa 800 2200}%
\special{pa 800 1900}%
\special{fp}%
\special{pa 800 2200}%
\special{pa 1500 2500}%
\special{fp}%
%
\special{pn 4}%
\special{pa 1500 2500}%
\special{pa 2000 2200}%
\special{fp}%
%
\special{pn 4}%
\special{pa 2200 1900}%
\special{pa 2500 1900}%
\special{fp}%
\end{picture}}%

\caption{An induced path of order $6$ in Claim~\ref{cl-prop-M91-case2-1}}
\label{fcl5.4}
\end{center}
\end{figure}
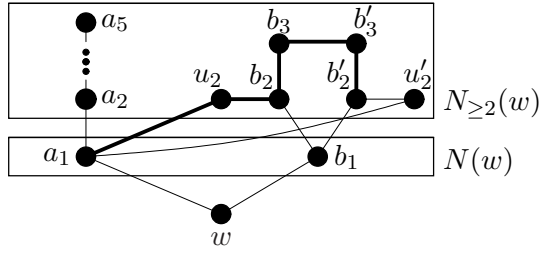

\begin{claim}
\label{cl-prop-M91-case2-1}
There exist $b_{3}\in N(b_{2})-\{b_{1}\}$ and $b'_{3}\in N(b'_{2})-\{b_{1}\}$ such that $b_{3},b'_{3}\in N_{\geq 2}(w)$, $b_{3}b'_{3}\notin E(G)$ and $N(a_{1})\cap \{b_{3},b'_{3}\}=\emptyset $.
\end{claim}
\proof
We first consider the case where $b_{2}$ or $b'_{2}$, say $b'_{2}$, satisfies $N(a_{1})\cap N(b'_{2})=\emptyset $.
Since $|N(b_{2})\cap N(a_{1})|\leq 1$, there exists a vertex $b_{3}\in N(b_{2})-(\{b_{1}\}\cap N(a_{1}))$.
Since $|N(b'_{2})\cap N(b_{3})|\leq 1$, there exists a vertex $b'_{3}\in N(b'_{2})-(\{b_{1}\}\cap N(b_{3}))$.
Then $b_{3},b'_{3}\in N_{\geq 2}(w)$ and $b_{3}b'_{3}\notin E(G)$.
By the assumption that $N(a_{1})\cap N(b'_{2})=\emptyset $, we have $b'_{3}a_{1}\notin E(G)$.
Consequently, $b_{3}$ and $b'_{3}$ satisfy the desired properties.
Thus we may assume that $N(a_{1})\cap N(b)\neq \emptyset $ for each $b\in \{b_{2},b'_{2}\}$.
Write $N(a_{1})\cap N(b_{2})=\{u_{2}\}$ and $N(a_{1})\cap N(b'_{2})=\{u'_{2}\}$.
Note that $u_{2},u'_{2}\in N_{2}(w)$.
Since $N(b_{2})\cap N(b'_{2})=\{b_{1}\}$, $u_{2}b'_{2}\notin E(G)$.
Take $b_{3}\in N(b_{2})-\{b_{1},u_{2}\}$ and $b'_{3}\in N(b'_{2})-\{b_{1},u'_{2}\}$.
Then $b_{3}\neq b'_{3}$ and $b_{3},b'_{3}\in N_{\geq 2}(w)$.
Since $N(a_{1})\cap N(b_{2})=\{u_{2}\}$ and $N(a_{1})\cap N(b'_{2})=\{u'_{2}\}$, we have $N(a_{1})\cap \{b_{3},b'_{3}\}=\emptyset $.
If $b_{3}b'_{3}\in E(G)$, then $a_{1}u_{2}b_{2}b_{3}b'_{3}b'_{2}$ is an induced path, which contradicts the assumption of Case~2 (see Figure~\ref{fcl5.4}).
Thus $b_{3}b'_{3}\notin E(G)$.
Hence $b_{3}$ and $b'_{3}$ satisfy the desired properties.
\qed

\begin{figure}
\begin{center}
{\unitlength 0.1in%
\begin{picture}(22.5000,11.6500)(4.0000,-25.6500)%
\put(15.0000,-26.3000){\makebox(0,0){$w$}}%
%
\special{sh 1.000}%
\special{ia 1500 2500 50 50 0.0000000 6.2831853}%
\special{pn 8}%
\special{ar 1500 2500 50 50 0.0000000 6.2831853}%
%
\special{pn 8}%
\special{pa 400 2100}%
\special{pa 2600 2100}%
\special{pa 2600 2300}%
\special{pa 400 2300}%
\special{pa 400 2100}%
\special{pa 2600 2100}%
\special{fp}%
%
\special{pn 8}%
\special{pa 400 1400}%
\special{pa 2600 1400}%
\special{pa 2600 2000}%
\special{pa 400 2000}%
\special{pa 400 1400}%
\special{pa 2600 1400}%
\special{fp}%
%
\special{sh 1.000}%
\special{ia 800 2200 50 50 0.0000000 6.2831853}%
\special{pn 8}%
\special{ar 800 2200 50 50 0.0000000 6.2831853}%
%
\special{sh 1.000}%
\special{ia 800 1900 50 50 0.0000000 6.2831853}%
\special{pn 13}%
\special{ar 800 1900 50 50 0.0000000 6.2831853}%
%
\special{sh 1.000}%
\special{ia 2000 2200 50 50 0.0000000 6.2831853}%
\special{pn 8}%
\special{ar 2000 2200 50 50 0.0000000 6.2831853}%
%
\special{sh 1.000}%
\special{ia 1800 1900 50 50 0.0000000 6.2831853}%
\special{pn 8}%
\special{ar 1800 1900 50 50 0.0000000 6.2831853}%
\put(6.5000,-22.0000){\makebox(0,0){$a_{1}$}}%
\put(21.5000,-22.0000){\makebox(0,0){$b_{1}$}}%
\put(26.5000,-23.0000){\makebox(0,0)[lb]{$N(w)$}}%
\put(26.5000,-20.0000){\makebox(0,0)[lb]{$N_{\geq 2}(w)$}}%
%
\special{sh 1.000}%
\special{ia 800 1500 50 50 0.0000000 6.2831853}%
\special{pn 8}%
\special{ar 800 1500 50 50 0.0000000 6.2831853}%
\put(9.5000,-19.0000){\makebox(0,0){$a_{2}$}}%
\put(9.5000,-15.0000){\makebox(0,0){$a_{5}$}}%
%
\special{sh 1.000}%
\special{ia 2200 1900 50 50 0.0000000 6.2831853}%
\special{pn 8}%
\special{ar 2200 1900 50 50 0.0000000 6.2831853}%
%
\special{pn 4}%
\special{sh 1}%
\special{ar 800 1655 16 16 0 6.2831853}%
\special{sh 1}%
\special{ar 800 1755 16 16 0 6.2831853}%
\special{sh 1}%
\special{ar 800 1705 16 16 0 6.2831853}%
\special{sh 1}%
\special{ar 800 1705 16 16 0 6.2831853}%
%
\special{pn 20}%
\special{pa 800 1900}%
\special{pa 800 1800}%
\special{fp}%
\special{pa 800 1600}%
\special{pa 800 1500}%
\special{fp}%
%
\special{sh 1.000}%
\special{ia 1800 1610 50 50 0.0000000 6.2831853}%
\special{pn 8}%
\special{ar 1800 1610 50 50 0.0000000 6.2831853}%
%
\special{sh 1.000}%
\special{ia 2200 1610 50 50 0.0000000 6.2831853}%
\special{pn 8}%
\special{ar 2200 1610 50 50 0.0000000 6.2831853}%
%
\special{pn 20}%
\special{pa 2000 2200}%
\special{pa 1800 1900}%
\special{fp}%
\special{pa 2200 1900}%
\special{pa 2000 2200}%
\special{fp}%
\put(23.4000,-19.0000){\makebox(0,0){$b'_{2}$}}%
%
\special{pn 20}%
\special{pa 800 2200}%
\special{pa 800 1900}%
\special{fp}%
\special{pa 800 2200}%
\special{pa 1500 2500}%
\special{fp}%
%
\special{pn 20}%
\special{pa 1500 2500}%
\special{pa 2000 2200}%
\special{fp}%
%
\special{pn 20}%
\special{pa 1800 1600}%
\special{pa 1800 1900}%
\special{fp}%
\special{pa 2200 1900}%
\special{pa 2200 1600}%
\special{fp}%
\put(19.4000,-19.0000){\makebox(0,0){$b_{2}$}}%
\put(19.4000,-16.0000){\makebox(0,0){$b_{3}$}}%
\put(23.4000,-16.0000){\makebox(0,0){$b'_{3}$}}%
\end{picture}}%

\caption{An induced copy of $T_{9}$ in Case~2 of Proposition~\ref{prop-M91}}
\label{fprop5.6-2}
\end{center}
\end{figure}
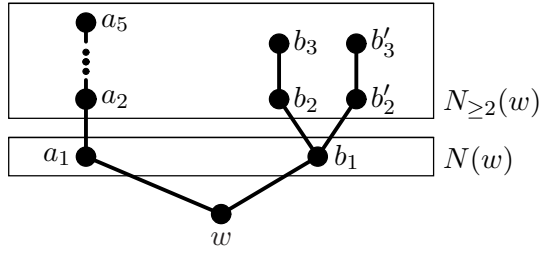

Let $b_{3}$ and $b'_{3}$ be as in Claim~\ref{cl-prop-M91-case2-1}.
Since $N(b_{2})\cap N(b'_{2})=\{b_{1}\}$, $b_{3}\neq b'_{3}$.
Hence
$$
E(G[\{w,a_{1},b_{1},b_{2},b_{3},b'_{2},b'_{3}\}])=\{wa_{1},wb_{1},b_{1}b_{2},b_{2}b_{3},b_{1}b'_{2},b'_{2}b'_{3}\}.
$$
Since $\{b_{1},b_{2},b_{3},b'_{2},b'_{3}\}\subseteq N_{\leq 2}(b_{1})-Z_{3}$, it follows from Lemma~\ref{lem-Xunion}(ii) that $E_{G}(X,\{b_{1},b_{2},b_{3},b'_{2},b'_{3}\})=\emptyset$.
Consequently, $\{b_{3},b_{2},b_{1},b'_{3},b'_{2},w,a_{1},\ldots ,a_{5}\}$ induces a copy of $T_{9}$ in $G$ (see Figure~\ref{fprop5.6-2}).

\medskip
\noindent
\textbf{Case 3:} $k_{0}\leq 4$.

\begin{claim}
\label{cl-prop-M91-1}
We have $N_{\geq 3}(w)=\emptyset $, and $G[N_{2}(w)]$ is $2$-regular.
\end{claim}
\proof
Suppose that $N_{\geq 3}(w)\neq \emptyset $ or $G[N_{2}(w)]$ is not $2$-regular.
Since $|N(u)\cap N(w)|=1$ for $u\in N_{2}(w)$ and $N(u')\cap N(w)=\emptyset $ for $u'\in N_{3}(w)$, there exists $a_{3}\in N_{\geq 2}(w)$ such that $N(a_{3})\cap N_{2}(w)\neq \emptyset $ and $|N(a_{3})\cap N_{\geq 2}(w)|\geq 3$.
Take $a_{2}\in N(a_{3})\cap N_{2}(w)$, and write $N(a_{2})\cap N(w)=\{a_{1}\}$.
Take $a_{4},a'_{4}\in (N(a_{3})\cap N_{\geq 2}(w))-\{a_{2}\}$ with $a_{4}\neq a'_{4}$.

Suppose that one of $a_{4}$ and $a'_{4}$ belongs to $N_{\geq 3}(w)$, say $a_{4}\in N_{\geq 3}(w)$.
Since $|N(a_{4})\cap N(a_{1})|\leq 1$, there exists a vertex $a_{5}\in N(a_{4})-(\{a_{3}\}\cup N(a_{1}))$.
Then $a_{5}\in N_{\geq 2}(w)$.
Since $G$ is $\{C_{3},C_{4}\}$-free, $a_{1}a_{2}a_{3}a_{4}a_{5}$ is an induced path of $G$, which contradicts the assumption of Case~3.
Thus $a_{4},a'_{4}\in N_{2}(w)$.
Write $N(a_{4})\cap N(w)=\{b_{1}\}$ and $N(a'_{4})\cap N(w)=\{b'_{1}\}$.
Since $N(a_{3})\cap N(a_{1})=\{a_{2}\}$ and $N(a_{4})\cap N(a'_{4})=\{a_{3}\}$, three vertices $a_{1},b_{1},b'_{1}$ are distinct.
Take $a'_{3}\in N(a_{2})-\{a_{1},a_{3}\}$.
Note that $a'_{3}\in N_{\geq 2}(w)$.
Since $|N(a'_{3})\cap \{b_{1},b'_{1}\}|\leq |N(a'_{3})\cap N(w)|\leq 1$, we may assume that $a'_{3}b_{1}\notin E(G)$.
Then $b_{1}a_{4}a_{3}a_{2}a'_{3}$ is an induced path of $G$, which contradicts the assumption of Case~3.
\qed

\begin{claim}
\label{cl-prop-M91-2}
Let $C$ be a component of $G[N_{2}(w)]$.
Then $C$ is a cycle with $|V(C)|\equiv 0~({\rm mod~}3)$ and $|N(V(C))\cap N(w)|=3$.
Further, if we write $C=b_{1}b_{2}\cdots b_{3k}b_{1}~(|V(C)|=3k)$, then $N(b_{i})\cap N(w)=N(b_{j})\cap N(w)$ for any $i$ and $j$ with $i\equiv j~({\rm mod~}3)$.
\end{claim}
\proof
By Claim~\ref{cl-prop-M91-1}, $C$ is a cycle.
Let $l=|V(C)|$, and write $C=b_{1}b_{2}\cdots b_{l}b_{1}$ (indices are to be read modulo $l$).
Let $1\leq i\leq l$, and write $N(b_{i})\cap N(w)=\{a\}$.
Since $G$ is $\{C_{3},C_{4}\}$-free, $l\geq 5$, and hence $b_{i}b_{i+1}b_{i+2}b_{i+3}$ is an induced path of $G$.
On the other hand, $ab_{i}b_{i+1}b_{i+2}b_{i+3}$ is not an induced path by the assumption of Case~3.
This forces $ab_{i+3}\in E(G)$, i.e., $N(b_{i+3})\cap N(w)=\{a\}$.
Since $i$ is arbitrary, we obtain the desired conclusion.
\qed

\begin{figure}
\begin{center}
{\unitlength 0.1in%
\begin{picture}(34.5500,13.6500)(6.9500,-22.6500)%
\put(24.0000,-23.3000){\makebox(0,0){$w$}}%
%
\special{sh 1.000}%
\special{ia 2400 2200 50 50 0.0000000 6.2831853}%
\special{pn 8}%
\special{ar 2400 2200 50 50 0.0000000 6.2831853}%
%
\special{pn 8}%
\special{pa 1200 1800}%
\special{pa 3600 1800}%
\special{pa 3600 2000}%
\special{pa 1200 2000}%
\special{pa 1200 1800}%
\special{pa 3600 1800}%
\special{fp}%
%
\special{pn 8}%
\special{pa 700 900}%
\special{pa 4100 900}%
\special{pa 4100 1700}%
\special{pa 700 1700}%
\special{pa 700 900}%
\special{pa 4100 900}%
\special{fp}%
\put(36.5000,-20.0000){\makebox(0,0)[lb]{$N(w)$}}%
\put(41.5000,-17.0000){\makebox(0,0)[lb]{$N_{2}(w)$}}%
%
\special{sh 1.000}%
\special{ia 1200 1200 50 50 0.0000000 6.2831853}%
\special{pn 8}%
\special{ar 1200 1200 50 50 0.0000000 6.2831853}%
%
\special{sh 1.000}%
\special{ia 2000 1200 50 50 0.0000000 6.2831853}%
\special{pn 8}%
\special{ar 2000 1200 50 50 0.0000000 6.2831853}%
%
\special{sh 1.000}%
\special{ia 1200 1400 50 50 0.0000000 6.2831853}%
\special{pn 8}%
\special{ar 1200 1400 50 50 0.0000000 6.2831853}%
%
\special{sh 1.000}%
\special{ia 1600 1600 50 50 0.0000000 6.2831853}%
\special{pn 8}%
\special{ar 1600 1600 50 50 0.0000000 6.2831853}%
%
\special{sh 1.000}%
\special{ia 2000 1400 50 50 0.0000000 6.2831853}%
\special{pn 8}%
\special{ar 2000 1400 50 50 0.0000000 6.2831853}%
%
\special{sh 1.000}%
\special{ia 1600 1900 50 50 0.0000000 6.2831853}%
\special{pn 8}%
\special{ar 1600 1900 50 50 0.0000000 6.2831853}%
\put(14.5000,-19.2000){\makebox(0,0){$a$}}%
%
\special{pn 4}%
\special{pa 1400 1100}%
\special{pa 1200 1200}%
\special{fp}%
\special{pa 1200 1200}%
\special{pa 1200 1400}%
\special{fp}%
\special{pa 1200 1400}%
\special{pa 1600 1600}%
\special{fp}%
\special{pa 1600 1600}%
\special{pa 2000 1400}%
\special{fp}%
\special{pa 2000 1400}%
\special{pa 2000 1200}%
\special{fp}%
\special{pa 2000 1200}%
\special{pa 1800 1100}%
\special{fp}%
%
\special{pn 4}%
\special{sh 1}%
\special{ar 1600 1100 16 16 0 6.2831853}%
\special{sh 1}%
\special{ar 1500 1100 16 16 0 6.2831853}%
\special{sh 1}%
\special{ar 1700 1100 16 16 0 6.2831853}%
\special{sh 1}%
\special{ar 1700 1100 16 16 0 6.2831853}%
\put(9.7000,-12.0000){\makebox(0,0){$b_{3k-1}$}}%
\put(10.4000,-14.0000){\makebox(0,0){$b_{3k}$}}%
\put(16.0000,-14.5000){\makebox(0,0){$b_{1}$}}%
\put(21.4000,-14.0000){\makebox(0,0){$b_{2}$}}%
\put(21.4000,-12.0000){\makebox(0,0){$b_{3}$}}%
%
\special{sh 1.000}%
\special{ia 2805 1200 50 50 0.0000000 6.2831853}%
\special{pn 8}%
\special{ar 2805 1200 50 50 0.0000000 6.2831853}%
%
\special{sh 1.000}%
\special{ia 3605 1200 50 50 0.0000000 6.2831853}%
\special{pn 8}%
\special{ar 3605 1200 50 50 0.0000000 6.2831853}%
%
\special{sh 1.000}%
\special{ia 2805 1400 50 50 0.0000000 6.2831853}%
\special{pn 8}%
\special{ar 2805 1400 50 50 0.0000000 6.2831853}%
%
\special{sh 1.000}%
\special{ia 3205 1600 50 50 0.0000000 6.2831853}%
\special{pn 8}%
\special{ar 3205 1600 50 50 0.0000000 6.2831853}%
%
\special{sh 1.000}%
\special{ia 3605 1400 50 50 0.0000000 6.2831853}%
\special{pn 8}%
\special{ar 3605 1400 50 50 0.0000000 6.2831853}%
%
\special{pn 4}%
\special{pa 3005 1100}%
\special{pa 2805 1200}%
\special{fp}%
\special{pa 2805 1200}%
\special{pa 2805 1400}%
\special{fp}%
\special{pa 2805 1400}%
\special{pa 3205 1600}%
\special{fp}%
\special{pa 3205 1600}%
\special{pa 3605 1400}%
\special{fp}%
\special{pa 3605 1400}%
\special{pa 3605 1200}%
\special{fp}%
\special{pa 3605 1200}%
\special{pa 3405 1100}%
\special{fp}%
%
\special{pn 4}%
\special{sh 1}%
\special{ar 3205 1100 16 16 0 6.2831853}%
\special{sh 1}%
\special{ar 3105 1100 16 16 0 6.2831853}%
\special{sh 1}%
\special{ar 3305 1100 16 16 0 6.2831853}%
\special{sh 1}%
\special{ar 3305 1100 16 16 0 6.2831853}%
\put(25.5500,-12.0000){\makebox(0,0){$b'_{3k'-1}$}}%
\put(26.2500,-14.0000){\makebox(0,0){$b'_{3k'}$}}%
\put(32.0500,-14.5000){\makebox(0,0){$b'_{1}$}}%
\put(37.4500,-14.0000){\makebox(0,0){$b'_{2}$}}%
\put(37.4500,-12.0000){\makebox(0,0){$b'_{3}$}}%
%
\special{sh 1.000}%
\special{ia 3200 1900 50 50 0.0000000 6.2831853}%
\special{pn 8}%
\special{ar 3200 1900 50 50 0.0000000 6.2831853}%
\put(33.5000,-18.9000){\makebox(0,0){$a'$}}%
%
\special{pn 20}%
\special{pa 2400 2200}%
\special{pa 1600 1900}%
\special{fp}%
\special{pa 1600 1900}%
\special{pa 1600 1600}%
\special{fp}%
\special{pa 1600 1600}%
\special{pa 2000 1400}%
\special{fp}%
\special{pa 2000 1400}%
\special{pa 2000 1200}%
\special{fp}%
\special{pa 1600 1600}%
\special{pa 1200 1400}%
\special{fp}%
\special{pa 1200 1400}%
\special{pa 1200 1200}%
\special{fp}%
%
\special{pn 20}%
\special{pa 2400 2200}%
\special{pa 3200 1900}%
\special{fp}%
\special{pa 3200 1900}%
\special{pa 3200 1600}%
\special{fp}%
\special{pa 3200 1600}%
\special{pa 3600 1400}%
\special{fp}%
\special{pa 3600 1400}%
\special{pa 3600 1200}%
\special{fp}%
\end{picture}}%

\caption{An induced copy of $T_{9}$ in Case~3 of Proposition~\ref{prop-M91}}
\label{fprop5.6-3}
\end{center}
\end{figure}
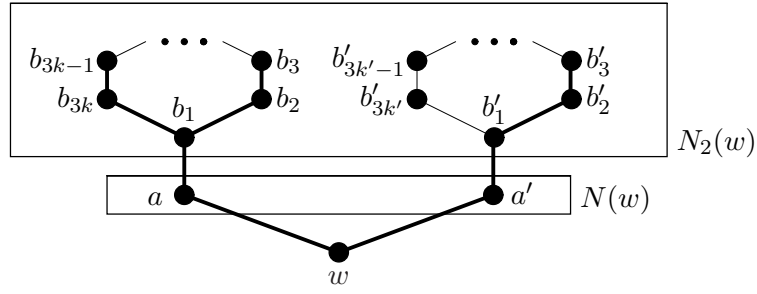

By Claim~\ref{cl-prop-M91-2}, $|N(V(C))\cap N(w)|=3$ for every component $C$ of $G[N_{2}(w)]$.
Since $d_{G}(w)>3$, there are two components $D$ and $D'$ of $G[N_{2}(w)]$ such that $N(V(D))\cap N(w)\neq N(V(D'))\cap N(w)$, and we can take $a\in (N(V(D))-N(V(D')))\cap N(w)$ and $a'\in (N(V(D'))-N(V(D)))\cap N(w)$.
Take $b_{1}\in N(a)\cap V(D)$ and $b'_{1}\in N(a')\cap V(D')$, and write $D=b_{1}b_{2}\cdots b_{3k}b_{1}$ and $D'=b'_{1}b'_{2}\cdots b'_{3k'}b'_{1}$.
In view of Claim~\ref{cl-prop-M91-2}, $\{b_{3},b_{2},b_{1},b_{3k-1},b_{3k},a,w,a',b'_{1},b'_{2},b'_{3}\}$ induces a copy of $T_{9}$ in $G$ (see Figure~\ref{fprop5.6-3}).
\qed

\end{document}